\documentclass[11pt,3p,sort&compress]{elsarticle}
\usepackage[latin9]{inputenc}
\usepackage{geometry}
\usepackage[active]{srcltx}
\usepackage{color}
\usepackage{booktabs}
\usepackage{multirow}
\usepackage{amsmath}
\usepackage{amssymb}
\usepackage{graphicx}
\usepackage[unicode=true,
 bookmarks=true,bookmarksnumbered=false,bookmarksopen=false,
 breaklinks=false,pdfborder={0 0 1},backref=false,colorlinks=true]
 {hyperref}

\makeatletter

\providecommand{\tabularnewline}{\\}

\journal{Comput Methods Appl Mech Eng}
\usepackage{multirow}
\usepackage{upgreek}
\usepackage{paralist}
\usepackage{scalerel}
\newcommand{\CIRC}{\hspace*{0.25pt}\scaleobj{0.5}{\boldsymbol{\circ}}}
\newcommand{\CCIRC}{\hspace*{0.15pt}\scaleobj{0.5}{\boldsymbol{\circ\circ}}}

\usepackage{xfrac}
\usepackage{xcolor}
\usepackage[bottom]{footmisc}

\@ifundefined{showcaptionsetup}{}{%
 \PassOptionsToPackage{caption=false}{subfig}}
\usepackage{subfig}
\makeatother

\begin{document}
\begin{frontmatter} 

\title{High-order implicit time integration scheme based on Padé expansions}

\author{C. Song\corref{cor1}}

\ead{c.song@unsw.edu.au}

\author{S.~Eisenträger}

\address{Centre for Infrastructure Engineering and Safety, School of Civil
and Environmental Engineering, University of New South Wales, Sydney,
NSW 2052, Australia.}

\cortext[cor1]{Corresponding author}
\begin{abstract}
A single-step high-order implicit time integration scheme for the
solution of transient and wave propagation problems is presented.
It is constructed from the Padé expansions of the matrix exponential
solution of a system of first-order ordinary differential equations
formulated in the state-space. A computationally efficient scheme
is developed exploiting the techniques of polynomial factorization
and partial fractions of rational functions, and by decoupling the
solution for the displacement and velocity vectors. An important feature
of the novel algorithm is that no direct inversion of the mass matrix
is required. From the diagonal Padé expansion of order $M$ a time-stepping
scheme of order $2M$ is developed. Here, each elevation of the accuracy
by two orders results in an additional system of real or complex sparse
equations to be solved. These systems are comparable in complexity
to the standard Newmark method, i.e., the effective system matrix
is a linear combination of the static stiffness, damping, and mass
matrices. It is shown that the second-order scheme is equivalent to
Newmark's constant average acceleration method, often also referred
to as trapezoidal rule. The proposed time integrator has been implemented
in \texttt{MATLAB} using the built-in direct linear equation solvers.
In this article, numerical examples featuring nearly one million degrees
of freedom are presented. High-accuracy and efficiency in comparison
with common second-order time integration schemes are observed. The
\texttt{MATLAB}-implementation is available from the authors upon
request or from the GitHub repository (to be added).
\end{abstract}
\end{frontmatter} 

\section{Introduction}

\label{sec:Introduction} The response of structures to time-varying
loads is of utmost importance to many branches of engineering and
science. Due to the complexity of practical problems, analytical solutions
are generally not available for transient analyses and therefore,
numerical methods are usually employed to approximate the structural
response. In finite element analyses (FEA) of transient problems direct
time integration methods are widely used. These methods are based
on two essential ideas: They (i) satisfy the equation of motion at
discrete points in time and (ii) assume a variation of displacements,
velocities, and accelerations within each time step \citep{BookBathe2002}.
A large number of time-stepping methods has been developed over the
last decades and novel approaches are still being proposed continuously.
In commercial finite element software and also in scientific applications
that are directed at studying transient problems the central difference
method (CDM), the family of methods belonging the Newmark's approach
\citep{ArticleNewmark1959}, and the HHT-$\alpha$ method \citep{ArticleHilber1977}
are pre-dominantly utilized. The popularity of these schemes may be
related to their ease of implementation and to the fact they are well-researched
and robust. Still one has to keep in mind that selecting an optimal
time-stepping method for a specific problem is crucially important
in terms of attainable accuracy and efficiency of the transient analysis
\citep{ArticleKim2020b}. The first critical consideration is whether
to apply an explicit or implicit approach. The basic property distinguishing
both types of methods is that in an explicit analysis the equilibrium
conditions are evaluated at time $t$, while in an implicit method
they are considered at time $t\,{+}\,\Delta t$, where $\Delta t$
is the selected time step size. This difference manifests itself in
the numerical properties of the different algorithms. In general,
it can be stated that explicit methods are conditionally stable, i.e.,
a maximum (critical) time step size exists which cannot be exceeded,
while most implicit method are unconditionally stable. This property
is obviously associated with computational cost. In implicit methods,
we typically have to solve a system of equations per each time step,
while simple matrix-vector products are sufficient to advance in time
for truly explicit methods in conjunction with diagonal mass matrices.
Thus, the efficiency of explicit methods critically depends on the
availability of highly accurate mass lumping techniques. At least
to the authors' knowledge, accurate mass lumping schemes have only
been proposed for linear finite elements \citep{ArticleDuczek2019a}
and high-order (tensor product-based) spectral elements \citep{ArticleDuczek2019b}.
Special solutions are available for some selected element types such
as tetrahedral elements \citep{ArticleGeevers2018a,ArticleGeevers_2019,ArticleZhang2019b},
but not for arbitrary order. Otherwise, methods such as diagonal scaling
\citep{ArticleHinton1976} or dual basis functions \citep{ArticleAnitescu2019}
have to be employed, which exhibit a reduced order of convergence.
In this work, we will concentrate implicit methods and therefore,
mass lumping techniques are not essential for an efficient implementation.

In the last 16 years, composite time integration schemes have been
in the focus of the research community. Here, the goal is to harness
the advantages of different time integrators by combining them in
one scheme \citep{ArticleKim2020b}. Arguably the first composite
time-stepping method has been published by Bathe and Baig \citep{ArticleBathe2005},
where one time step is split into two sub-steps using Newmark's constant
average acceleration method in combination with a three-point backward
Euler formulation. In a series of papers, Bathe and co-workers tuned
the numerical properties to have adjustable numerical damping in the
high-frequency range \citep{ArticleBathe2007,ArticleBathe2012,ArticleNoh2013,ArticleNoh2018,ArticleNoh2019a,ArticleNoh2019b,ArticleMalakiyeh2019,ArticleMalakiyeh2021}.
Further composite schemes have been proposed by Kim and Reddy \citep{ArticleKim2017c}
and Kim and Choi \citep{ArticleKim2018a}. A comprehensive review
and an attempt to unify the different formulations has been reported
by Kim in Ref.~\citep{ArticleKim2020b}. By exploiting the idea of
composite time-stepping schemes, an infinite number of methods can
be devised which are tailor-made for different applications. However,
one has to keep in mind that splitting each time step into a certain
number of sub-steps will increase the numerical effort correspondingly
when compared to single-step methods such as Newmark's method \citep{ArticleNewmark1959},
the HHT-$\alpha$ method \citep{ArticleHilber1977}, the generalized-$\alpha$
method \citep{ArticleChung1993}, etc.

Another promising area of research in terms of time-stepping methods
is seen in locally adaptive time integration methods. At this point,
we particularly want to mention the work of Soares \citep{ArticleSoares2017a,ArticleSoares2018b}.
The main idea behind Soares' methods is to have an algorithm including
adaptive parameters that are computed locally. These parameters influence
the properties of the scheme such as numerical (algorithmic) damping.
Based on an oscillatory criterion algorithmic damping can be switched
on and off to damp out spurious waves with minimal effect on the overall
energy of the system. This algorithm is fully automatic and does not
require input from the user. It has to be stressed that the mentioned
adaptivity is not related to the choice of the time step size, which
could be included additionally. In subsequent publications, this approach
has been extended to IMEX (implicit-explicit) schemes \citep{ArticleSoares2019a,ArticleSoares2019b}.
IMEX time-marching techniques aim to exploit the best properties of
both implicit and explicit methods, i.e., in ``stiffer'' domains
(locally refined mesh, heterogeneous material, etc.) an implicit algorithm
is selected enabling a larger critical time step, while in more ``flexible''
regions an explicit method is employed.

The methods discussed above provide second-order accuracy. In the
wide body of literature, there are numerous attempts to construct
high-order schemes. In this context, The final goal is to develop
a unified high-order approach, i.e., high-order spatial \emph{and}
temporal discretizations, and consequently, further research must
be focused on high-order time-stepping methods. A few high-order schemes
are mentioned in the following.

Exploiting a modified continuous Galerkin formulation, Idesman~\citep{ArticleIdesman2007}
introduced the velocity vector as an additional time-dependent variable
to derive a high-order accurate time integrator. The displacement
and velocity vectors are approximated as polynomial functions of order
$N$ and thus, an accuracy of order $2N$ can be achieved. An issue
with this formulation is that the system size grows with $N$. Additionally,
it is stated in Ref.~\citep{ArticleKim2017b} that the method provides
only an $N$th order accurate scheme if implemented as stated in the
original paper by Idesman. Kim and Reddy also propose two high-order
methods based on a modified weighted residual method \citep{ArticleKim2017a,ArticleKim2017b}.
The first approach is based on Lagrange polynomials, while the second
one relies on Hermite polynomials. They have developed a systematic
procedure to derive algorithms that are ($2N-1$)th- and $2N$th-order
accurate if numerical damping is in- or excluded. Here, $N$ denotes
the order of the Lagrange polynomial. In the case of Lagrange polynomials,
Gauss-Lobatto-Legendre points are used for the interpolation process.
The formulae are explicitly provided for schemes up to tenth-order
accuracy. Since Hermite polynomials are C$_{1}$-continuous also the
first derivatives are included in the interpolation. Hence, the order
of accuracy of the derived schemes is $p$ or $p\,{+}\,1$ for a $p$th-order
Hermite polynomial. In this case, all expressions up to eighth-order
are included in the article. These approaches achieve very accurate
results, but the size of the system of simultaneous equations scales
with $N$ or $\sfrac{(p+1)}{2}$. In Ref.~\citep{ArticleKim2019c},
Kim and Lee compare the Lagrange-based method \citep{ArticleKim2017a}
with Fung's algorithm \citep{ArticleFung1997} and arrive at the conclusion
that the numerical characteristics are identical for linear analysis,
but in the case of nonlinear systems a better performance of Kim and
Reddy's approach is observed for long-duration simulations using large
time steps sizes.

Recently, Soares introduced a simple fourth-order accurate scheme
exhibiting the same computational effort as a single-step method of
first- or second-order \citep{ArticleSoares2020b} when neither physical
nor numerical damping is included. In the presence of numerical damping
the method becomes third-order accurate, while it is second-order
accurate if physical damping is taken into account. Other noteworthy
features of this method are that it is truly self-starting, as it
is based on displacement-velocity recursive relations and no computation
of accelerations is required. Due to the fact that the method is only
conditionally stable, its application is limited to problems where
a small time step size is already dictated by the physics of the problem.
An extension of the generalized-$\alpha$ method to high-order accuracy
has be proposed by Behnoudfar et al. \citep{ArticleBehnoudfar2020,ArticleBehnoudfar2021}.
The novel technique retains the advantages of the original method
such as unconditional stability and controllable numerical damping
and increases the order of accuracy by expanding the displacement
vector in a Taylor series. In each time step, $N$ matrix systems
have to be solved consecutively and implicitly. Thereafter, $2N$
auxiliary variables are updated explicitly.

A family of time integration methods has been proposed by making use
of the matrix theory, well-established in mathematics, for the solution
of systems of ordinary differential equations (ODEs)~\citep{Reusch1988,ArticleZhong1994,ArticleFung2008}.
The second-order semi-discrete equations of motion is transformed
into a system of first-order ODEs. The operation requires the inversion
of the mass matrix, which leads to a fully populated matrix if the
mass matrix is not diagonal. The exact solution can be obtained using
the variations-of-constants formula and involves a matrix exponential
function. Zhong and Williams \citep{ArticleZhong1994} derived the
precise time integration method (PIM), which is in the literature
also referred to as an exponential-type time integrator. The matrix
exponential is generally fully populated even if the mass matrix is
diagonal and therefore, it is extremely expensive to compute for large
problems~\citep{Golub1996}. In Ref.~\citep{ArticleZhong1994},
a Taylor series expansion is proposed to obtain the homogeneous solution.
The accuracy of the time-integrator is dependent on the number of
terms that are used. Since the direct computation of the Taylor series
is numerically expensive, a recursive algorithm is proposed to evaluate
the matrix exponential. Considering the inhomogeneous part of the
solution a linear variation of the force per time-step is assumed
limiting the applicability of the algorithm. Over the last decades,
several improvements of this approach have been published \citep{ArticleFung1997,ArticleFung2008,ArticleWang2007,ArticleWang2009}.
In 1997, Fung proposed a general solution in terms of step- and impulse-response
matrices. Based on this approach the system matrices remain symmetric,
but it is concluded that in general the derived high-order methods
are only conditionally stable. An extension of this approach was published
by Fung and Chen in 2008 \citep{ArticleFung2008}, where an additional
Duhamel response matrix was introduced. Wang and Au proposed to use
a Padé series instead of a Taylor series to evaluate the matrix exponential
\citep{ArticleWang2007}. In this approach, no inverse matrix is computed
to treat the inhomogeneous term, but the variation of the force is
limited to linear and sinusoidal functions. In a follow-up article,
Wang and Au recommended to approximate the force in each time step
by Lagrange polynomials defined on Gauss-Lobatto-Chebyshev points,
which provides more flexibility in terms of the excitation signal.
A similar approach is followed by Luan, who derived exponential Runge-Kutta
methods which are also applicable to systems of first-order ODEs \citep{ArticleLuan2021}.
As an example, a fourth and fifth-order accurate scheme are derived.

Techniques to improve the computational efficiency in treating rational
functions such as high-order Padé expansions are important and have
been developed in the past. Wolf \citep{Wolf1991} developed a consistent
lumped-parameter model by formulating a rational function as a partial-fraction
expansion. The technique breaks a complex high-order representation
into a series of simple (first- or second-order) ones. Later, it was
extended to handle complex dynamic soil-structure interaction problems
by Wolf and Paranesso~\citep{Wolf1992}. Birk and Ruge~\citep{BIRK20071152}
further developed this approach to perform dynamic dam-reservoir interaction
analysis. Reusch et al.~\citep{Reusch1988} derived a time integration
scheme for parabolic equations by explicitly factoring the polynomials
of the Padé expansion of the matrix exponential.

Overall, it can be concluded that, despite the excellent accuracy,
the wide-spread application of high-order time integration methods
for structural dynamics is hindered by their significant numerical
overhead. Therefore, our goal is to introduce an efficient high-order
method and compare its performance with established second-order algorithms
that are used both in engineering practice and scientific computations. 

\section{High-order accurate time integration scheme}

\label{sec:HighOrderPade} In this section, the theory needed for
the derivation of the high-order time integration scheme is explained.
Based on the exact solution of the equations of motion an approximation
using a Padé series expansion is employed. This yields a highly accurate
unconditionally stable implicit scheme that can be used for various
problems in structural dynamics and also wave propagation. 

\subsection{Time integration by matrix exponential}

\label{sec:MatrixExponent} The point of departure for the derivation
of any time integration scheme are the equations of motion, which
are semi-discrete system of second-order ODEs in time $t$, commonly
expressed as 
\begin{equation}
\mathbf{M}\ddot{\mathbf{u}}+\mathbf{C}\dot{\mathbf{u}}+\mathbf{K}\mathbf{u}=\mathbf{f}\,.\label{eq:EqOfMotion}
\end{equation}
Here, $\mathbf{M}$, $\mathbf{C}$, and $\mathbf{K}$ denote the mass,
damping, and stiffness matrices, respectively. Regarding the derivation
of a time integration approach, it is of no importance whether these
matrices are obtained from a variant of the finite element method
(FEM), a finite difference method (FDM), or a meshless method, to
name just a few. However, it is fair to assume that the semi-discrete
system of equations is generally sparse and often also symmetric.
In Eq.~\eqref{eq:EqOfMotion}, the displacement vector is represented
as $\mathbf{u}$, while its temporal derivatives are indicated by
an overdot as $\dot{\square}$, which is equivalent to $\square_{,t}$.
Hence, $\dot{\mathbf{u}}$ and $\ddot{\mathbf{u}}$ represent the
velocity and acceleration vectors, respectively. The vector $\mathbf{f}$
stands for the external excitation force. To solve Eq.~\eqref{eq:EqOfMotion},
initial conditions (ICs) at $t=0$ need to be provided 
\begin{equation}
\mathbf{u}(t=0)=\mathbf{u}_{0}\quad\mathrm{and}\quad\dot{\mathbf{u}}(t=0)=\dot{\mathbf{u}}_{0}\label{eq:InitialCondition}
\end{equation}
with the initial displacement vector $\mathbf{u}_{0}$ and the initial
velocity vector $\dot{\mathbf{u}}_{0}$.

Since it is generally intractable (and often impossible) to derive
closed-form solutions for the initial value problem (IVP) expressed
in Eqs.~\eqref{eq:EqOfMotion} and \eqref{eq:InitialCondition},
numerical approaches have to be applied. In this regard, time-stepping
methods are often employed and are very popular in the computational
dynamics community (see the comprehensive review articles by Tamma
et al. \citep{ArticleTamma2000,ArticleTamma2011}). In a numerical
time-stepping algorithm, the overall simulation time is divided into
a finite number of intervals, which is in principle quite similar
to a one-dimensional spatial discretization. Without loss of generality,
it is assumed that the time step $n$ spans the interval $[t_{n-1},t_{n}]$.
Thus, the time step size (increment) is simply defined as $\Delta t\,{=}\,t_{n}-t_{n-1}$,
where $n\,{=}\,1,2,\ldots,n_{\mathrm{S}}$ with the number of time
steps $n_{\mathrm{S}}$. For the sake of a simple and consistent notation,
a (local) dimensionless time variable $s$ is introduced for each
time step. $s$ is defined as $0$ at the beginning of a time step,
while it is equal to $1$ at the end. The time within the time step
$n$ is consequently defined as 
\begin{equation}
t(s)=t_{n-1}+s\Delta t\,,\quad0\le s\le1\,,\label{eq:dimensionlessTime}
\end{equation}
with $t(s\,{=}\,0)\,{=}\,t_{n-1}$ and $t(s\,{=}\,1)\,{=}\,t_{n}$.
Within a time step ($t_{n-1}\,{\le}\,t\,{\le}\,t_{n}$), the velocity
and acceleration vectors are expressed in terms of the dimensionless
time $s$ by applying the chain rule
\begin{subequations}
\begin{align}
\dot{\mathbf{u}} & =\dfrac{1}{\Delta t}\cfrac{\mathrm{d}\mathbf{u}}{\mathrm{d}s}\equiv\dfrac{1}{\Delta t}\overset{\CIRC}{\mathbf{u}}\,,\label{eq:velocityDimensionless}\\
\ddot{\mathbf{u}} & =\dfrac{1}{\Delta t^{2}}\cfrac{\mathrm{d^{2}}\mathbf{u}}{\mathrm{d}s^{2}}\equiv\dfrac{1}{\Delta t^{2}}\overset{\CCIRC}{\mathbf{u}}\,,\label{eq:accelerationDimensionless}
\end{align}
\label{eq:motionDimensionless}
\end{subequations}
where an overhead $\overset{\CIRC}{\square}$ denotes a derivative
with respect to the dimensionless time $s$. Now, we can rewrite the
semi-discrete equations of motion by substituting Eqs.~\eqref{eq:motionDimensionless}
into Eq.~\eqref{eq:EqOfMotion} 
\begin{equation}
\mathbf{M}\overset{\CCIRC}{\mathbf{u}}+\Delta t\mathbf{C}\overset{\CIRC}{\mathbf{u}}+\Delta t^{2}\mathbf{K}\mathbf{u}=\Delta t^{2}\mathbf{f}\,.\label{eq:equationOfMotionDimensionless}
\end{equation}
At the time $t_{n-1}$, i.e., $s\,{=}\,0$, the displacement and velocity
responses are expressed as 
\begin{equation}
\mathbf{u}(s=0)=\mathbf{u}_{n-1}\quad\mathrm{and}\quad\overset{\CIRC}{\mathbf{u}}(s=0)=\Delta t\dot{\mathbf{u}}_{n-1}\,.\label{eq:condition0}
\end{equation}
In the remainder of the theoretical derivation, it is assumed, without
loss of generality, that the matrices $\mathbf{M}$, $\mathbf{C}$,
and $\mathbf{K}$ are constant within a single time step. Furthermore,
we introduce the state-space vector $\mathbf{z}$ defined as 
\begin{equation}
\mathbf{z}=\begin{Bmatrix}\overset{\CIRC}{\mathbf{u}}\\
\mathbf{u}
\end{Bmatrix}\,.\label{eq:stateVariable}
\end{equation}
Using this definition, Eq.~\eqref{eq:equationOfMotionDimensionless}
is transformed into a system of first-order ODEs 
\begin{equation}
\overset{\CIRC}{\mathbf{z}}\equiv\cfrac{\mathrm{d}\mathbf{z}}{\mathrm{d}s}=\mathbf{A}\mathbf{z}+\mathbf{F}\label{eq:ODE1st}
\end{equation}
with the constant coefficient matrix $\mathbf{A}$ 
\begin{equation}
\mathbf{A}=\begin{bmatrix}-\Delta t\mathbf{M}^{-1}\mathbf{C} & -\Delta t^{2}\mathbf{M}^{-1}\mathbf{K}\\
\mathbf{I} & \mathbf{0}
\end{bmatrix}\label{eq:matrixA}
\end{equation}
and the force term 
\begin{equation}
\mathbf{F}=\begin{Bmatrix}\Delta t^{2}\mathbf{M}^{-1}\mathbf{f}\\
\mathbf{0}
\end{Bmatrix}\,.\label{eq:vectorf}
\end{equation}
Here, $\mathbf{I}$ denotes the identity matrix of the same size as
$\mathbf{K}$. Note that when the damping matrix $\mathbf{C}$ vanishes,
i.e., no physical damping is present in the system, all eigenvalues
of $\mathbf{A}$ are imaginary and proportional to the eigenvalues
of $\mathbf{M}^{-1}\mathbf{K}$, which determine the natural/resonant
frequencies. To be precise, the following relation holds 
\begin{equation}
\lambda(\mathbf{A})\,{=}\,\pm\Delta t\lambda(\mathbf{M}^{-1}\mathbf{K})\mathrm{i}\,,
\end{equation}
where $\mathrm{i}$ denotes the imaginary unit being defined as $\mathrm{i}\,{=}\,\sqrt{-1}$.

Using the matrix exponential function (see \ref{sec:Appendix:-Matrix-exponential})
and the variations-of-constants formula, the general solution of Eq.~\eqref{eq:ODE1st}
is expressed as 
\begin{equation}
\mathbf{z}=e^{\mathbf{A}s}\left(\mathbf{c}+\int\limits _{0}^{s}e^{-\mathbf{A}\tau}\mathbf{F}(\tau)\,\mathrm{d}\tau\right)\,,\label{eq:generalSolution}
\end{equation}
where $\mathbf{c}$ is the vector of integration constants. Considering
the displacement and velocity values at the beginning of a time step---see
Eq.~\eqref{eq:condition0}---and the definition of the state-space
vector---see Eq.~\eqref{eq:stateVariable}, the integration constants
are determined by substituting $s\,{=}\,0$ into Eq.~\eqref{eq:generalSolution}
which yields
\begin{equation}
\mathbf{c}=\mathbf{z}_{n-1}=\begin{Bmatrix}\Delta t\dot{\mathbf{u}}_{n-1}\\
\mathbf{u}_{n-1}
\end{Bmatrix}\,.\label{eq:IntegrationConstant}
\end{equation}
Substituting Eq.~\eqref{eq:IntegrationConstant} into Eq.~\eqref{eq:generalSolution},
the solution of Eq.~\eqref{eq:ODE1st} is expressed as 
\begin{equation}
\mathbf{z}=e^{\mathbf{A}s}\mathbf{z}_{n-1}+e^{\mathbf{A}s}\int\limits _{0}^{s}e^{-\mathbf{A}\tau}\mathbf{F}(\tau)\,\mathrm{d}\tau\,.\label{eq:sln_ds}
\end{equation}
The solution at time $t_{n}$, i.e., for $s=1$, is obtained as 
\begin{equation}
\mathbf{z}_{n}=e^{\mathbf{A}}\mathbf{z}_{n-1}+e^{\mathbf{A}}\int\limits _{0}^{1}e^{-\mathbf{A}\tau}\mathbf{F}(\tau)\,\mathrm{d}\tau\,.\label{eq:sln-dn}
\end{equation}

In order to be able to derive a closed-form (analytical) solution
for the integral expression in Eq.~\eqref{eq:sln-dn}, the excitation
force vector $\mathbf{f}$, and thus, $\mathbf{F}$---see Eq.~\eqref{eq:vectorf},
is assumed to be sufficiently smooth within a time step. It is approximated
by a polynomial function in the dimensionless time $s$. For time
step $n$, it can be either expressed as a Taylor series expansion
\begin{equation}
\mathbf{F}_{n}(s)=\sum\limits _{k=0}^{p_{\mathrm{f}}}\cfrac{1}{k!}\cfrac{\mathrm{d}^{k}\mathbf{F}_{\mathrm{m}n}}{\mathrm{d}s^{k}}\;(s-0.5)^{k}=\mathbf{F}_{\mathrm{m}n}+\overset{\CIRC}{\mathbf{F}}_{\mathrm{m}n}(s-0.5)+\cfrac{1}{2}\;\overset{\CCIRC}{\mathbf{F}}_{\mathrm{m}n}(s-0.5)^{2}+\ldots+\cfrac{1}{p_{\mathrm{f}}!}\;\mathbf{F}_{\mathrm{m}n}^{(p_{\mathrm{f}})}(s-0.5)^{p_{\mathrm{f}}}\,,\label{eq:ForceExpansionTaylor}
\end{equation}
or as a polynomial function that is determined by curve-fitting methods,
e.g., least-squares, 
\begin{equation}
\mathbf{F}_{n}(s)=\sum\limits _{k=0}^{p_{\mathrm{f}}}\tilde{\mathbf{F}}_{\mathrm{m}n}^{(k)}(s-0.5)^{k}=\tilde{\mathbf{F}}_{\mathrm{m}n}^{(0)}+\tilde{\mathbf{F}}_{\mathrm{m}n}^{(1)}(s-0.5)+\tilde{\mathbf{F}}_{\mathrm{m}n}^{(2)}(s-0.5)^{2}+\ldots+\tilde{\mathbf{F}}_{\mathrm{m}n}^{(p_{\mathrm{f}})}(s-0.5)^{p_{\mathrm{f}}}\,.\label{eq:ForceExpansion}
\end{equation}
Note that the force-terms in Eqs.~\eqref{eq:ForceExpansionTaylor}
and \eqref{eq:ForceExpansion} are closely related and can be easily
converted into each other. The main difference in both approaches
is that the Taylor series is more suitable for deriving analytical
solutions, while the curve-fitting is straightforward to implement
in numerical algorithms. In Eq.~\eqref{eq:ForceExpansionTaylor},
$\mathbf{F}_{\mathrm{m}n}$ denotes force vector at the middle of
time step $n$, i.e., at $t_{n-\sfrac{1}{2}}\,{=}\,(n-\sfrac{1}{2})\Delta t$,
and the other terms represent the derivatives of the force vector
with respect to the dimensionless time $s$ at the interval midpoint.
These vectors can be easily determined by an analytical differentiation
of the forcing function $\mathbf{F}(s)$. The expansion consists of
$p_{\mathrm{f}}\,{+}\,1$ terms, where $p_{\mathrm{f}}$ is the order
of the polynomial approximation of the forcing function. Based on
an approximation by means of Lagrangian polynomials defined at the
Gau{ß}-Lobatto-Legendre (GLL) points, mapped to the interval $[0,1]$,
the polynomial function within a time step is easily defined (for
more information see \ref{sec:LagrangeInterpolation}). Despite the
simplicity of this scheme, we decided to implement the curve-fitting
approach in our code. Based on a least-squares fit of the polynomial
function given in Eq.~\eqref{eq:ForceExpansion}, the vectors $\tilde{\mathbf{F}}_{\mathrm{m}n}^{(k)}$
are calculated. As the fitting procedure is also based on GLL-points
and the values of the original force vector at those points, the following
relation holds 
\begin{equation}
\tilde{\mathbf{F}}_{\mathrm{m}n}^{(k)}=\cfrac{1}{k!}\cfrac{\mathrm{d}^{k}\mathbf{F}_{\mathrm{m}n}}{\mathrm{d}s^{k}}\,.
\end{equation}
Hence, both approximations are equivalent and the implementation according
to Eq.~\eqref{eq:ForceExpansion} only holds advantages in terms
of programming as mentioned before. The primary intention of introducing
this approximation is to evaluate the integral term in Eq.~\eqref{eq:sln-dn}
analytically. In the following, the integration of the terms $(s-0.5)^{k}$,
where $k\in\{0,1,2,\ldots,p_{\mathrm{f}}\}$, is detailed. The general
solution is derived using integration by parts and can be written
as a recurrence relation 
\begin{equation}
\mathbf{B}_{k}=e^{\mathbf{A}}\int\limits _{0}^{1}(\tau-0.5)^{k}e^{-\mathbf{A}\tau}\mathrm{d}\tau=\mathbf{A}^{-1}\left(k\mathbf{B}_{k-1}+\left(-\cfrac{1}{2}\right)^{k}(e^{\mathbf{A}}-(-1)^{k}\mathbf{I})\right)\quad\;\forall k=0,1,2,\ldots,p_{\mathrm{f}}\,,\label{eq:intergrationExpmI}
\end{equation}
where the solution of the integral for $k$ depends on the solution
for $k\,{-}\,1$. Considering the constant term ($k\,{=}\,0$) which
is denoted as $\mathbf{B}_{0}$, the solution to the integral simplifies
to 
\begin{equation}
\mathbf{B}_{0}=e^{\mathbf{A}}\int\limits _{0}^{1}e^{-\mathbf{A}\tau}\mathrm{d}\tau=\mathbf{A}^{-1}\left(e^{\mathbf{A}}-\mathbf{I}\right)\,,\label{eq:integrationExpm0}
\end{equation}
which is the starting point for the recursion. Substituting Eq.~\eqref{eq:ForceExpansion}
along with Eq.~\eqref{eq:intergrationExpmI} into Eq.~\eqref{eq:sln-dn},
the overall solution is expressed as %
\begin{equation}
\mathbf{z}_{n}=e^{\mathbf{A}}\mathbf{z}_{n-1}+\sum\limits _{k=0}^{p_{\mathrm{f}}}\mathbf{B}_{k}\tilde{\mathbf{F}}_{\mathrm{m}n}^{(k)}\,.\label{eq:timeSteppingExpm}
\end{equation}
The right-hand side of this equation is known at time $t_{n}$ and
therefore, the time-stepping can be easily performed starting at $n\,{=}\,1$
with the prescribed ICs already known at $n\,{=}\,0$. This equation
represents an exact time-stepping scheme if the excitation forces
vary according to a polynomial function within each time step. In
this case, exact means the solution to discretized system is computed
with an accuracy up to machine precision. This does not mean that
the results are physically accurate as the error introduced by the
spatial discretization is not accounted for.

Various algorithms for accurately computing the matrix exponential
function $e^{\mathbf{A}}$ have been devised~\citep{Golub1996}.
The ``scaling and squaring'' algorithm, which is based on a Padé
series expansion of the exponential function, is often employed and
also available in commercial software such as \texttt{MATLAB} (implemented
in the command \texttt{expm(x)}). However, the direct use of Eq.~\eqref{eq:timeSteppingExpm}
for time-stepping is not practical for large-scale problems since
the computational costs increase rapidly with the number of DOFs.
The main reasons are listed in the following: 
\begin{enumerate}
\item The definition of matrix $\mathbf{A}$---see Eq.~\eqref{eq:matrixA}---involves
the matrix products $\mathbf{M}^{-1}\mathbf{C}$ and $\mathbf{M}^{-1}\mathbf{K}$.
Considering the use of a consistent mass matrix formulation in contrast
to a lumped one\footnote{Note that it is not possible to derive a variationally consistent
formulation of a diagonal mass matrix \citep{ArticleDuczek2019b}.
However, if possible, an optimal convergence is not always guaranteed
\citep{ArticleDuczek2019a}.}, the resulting products will be full matrices, which are less efficient
to treat and require significantly more computer memory compared to
the sparse matrices $\mathbf{M}$, $\mathbf{C},$ and $\mathbf{K}$. 
\item The result of computing the matrix exponential function $e^{\mathbf{A}}$
is a fully populated matrix. Its computation involves matrix multiplications
as well as the solution of a dense system of simultaneous equations
and is consequently, highly expensive. 
\end{enumerate}
In summary, due to both the high memory requirements and high demands
on the computational resources, the direct use of this algorithm is
intractable. Therefore, the matrix exponential function needs to be
approximated in a suitable way that guarantees a reasonable accuracy
and can be efficiently implemented. One idea is based on the Padé
series expansion and will be discussed in detail in the remainder
of this section. 

\subsection{Time-stepping using a Padé expansion of the matrix exponential function}

\label{sec:Pade-TimeStepping} To reduce the computational costs,
the matrix exponential in Eq.~\eqref{eq:timeSteppingExpm} can be
approximated by simpler and computationally more efficient functions.
The use of polynomial approximation techniques such as the Taylor
expansion---see Eq.~\eqref{eq:MPFunc-1} in \ref{sec:Appendix:-Matrix-exponential},
will lead to \emph{explicit} time-stepping schemes that are only conditionally
stable~\citep{Gallopoulos1989}, i.e., a critical time increment
exists that must not be exceeded. In contrast, the use of approximation
techniques based on rational functions such as the Padé expansion
(see \ref{sec:Pade-approximation}) will lead to \emph{implicit} algorithms
that can be unconditionally stable. In this case, the size of the
time step is only dependent on the accuracy requirements on the response
history.

We decided to apply the diagonal Padé expansion\footnote{The term \emph{diagonal} Padé series expansion refers to the fact
that the polynomial orders of the numerator $L$ and denominator $M$
are identical. Therefore, it is sufficient to indicate the order by
just one value $M$.}---see Eq.~\eqref{eq:ExpPade}---to approximate the matrix exponential
$e^{\mathbf{A}}$. The diagonal Padé approximation of order $M$ is
expressed as 
\begin{equation}
e^{\mathbf{A}}\approx e_{M}^{\mathbf{A}}=\mathbf{Q}_{M}^{-1}(\mathbf{A})\mathbf{P}_{M}(\mathbf{A})\,,\label{eq:ExpmPade}
\end{equation}
where the polynomials $\mathbf{P}_{M}(\mathbf{A})$ and $\mathbf{Q}_{M}(\mathbf{A})$
\begin{subequations}
\label{eq:Exp_M_PQ} 
\begin{align}
\mathbf{P}_{M}(\mathbf{A}) & =\sum\limits _{m=0}^{M}\cfrac{(2M-m)!}{m!(M-m)!}\mathbf{A}^{m}\,,\label{eq:Exp_M_P}\\
\mathbf{Q}_{M}(\mathbf{A}) & =\sum\limits _{m=0}^{M}\cfrac{(2M-m)!}{m!(M-m)!}(-\mathbf{A})^{m}\,,\label{eq:Exp_M_Q}
\end{align}
\end{subequations}
are obtained from Eqs.~\eqref{eq:Exp_PM} and \eqref{eq:Exp_QM},
respectively. Unless necessary, the order $M$ and the argument $\mathbf{A}$
will be omitted hereafter for simplicity of notation. Pre-multiplying
Eq.~\eqref{eq:timeSteppingExpm} with $\mathbf{Q}$ and using Eq.~\eqref{eq:ExpmPade}
leads to
\begin{equation}
\mathbf{Q}\mathbf{z}_{n}=\mathbf{P}\mathbf{z}_{n-1}+\sum\limits _{k=0}^{p_{\mathrm{f}}}\mathbf{C}_{k}\tilde{\mathbf{F}}_{\mathrm{m}n}^{(k)}\,,\label{eq:PadeStepping}
\end{equation}
where the matrices $\mathbf{C}_{k}$ are introduced and expressed
as polynomials of $\mathbf{A}$ employing Eqs.~\eqref{eq:intergrationExpmI}
and \eqref{eq:ExpmPade}. The general formula to determine $\mathbf{C}_{k}$
is given as 
\begin{equation}
\mathbf{C}_{k}=\mathbf{Q}\mathbf{B}_{k}=\mathbf{A}^{-1}\left(k\mathbf{C}_{k-1}+\left(-\cfrac{1}{2}\right)^{k}(\mathbf{P}-(-1)^{k}\mathbf{Q})\right)\quad\;\forall k=0,1,2,\ldots,p_{\mathrm{f}}\,,\label{eq:PadeStepping_C}
\end{equation}
which can be further simplified for $k\,{=}\,0$: 
\begin{equation}
\mathbf{C}_{0}=\mathbf{Q}\mathbf{B}_{0}=\mathbf{A}^{-1}\left(\mathbf{P}-\mathbf{Q}\right)\,.\label{eq:PadeStepping_C0}
\end{equation}

The time-stepping scheme illustrated in Eq.~\eqref{eq:PadeStepping}
still requires the explicit computation of the matrix $\mathbf{A}$
and therefore, a computationally efficient implementation procedure
needs to be devised, which is described in the next section. 

\subsection{Efficient computational implementation}

\label{sec:Efficient-implementation} To achieve a computationally
efficient implementation\footnote{Remark: For the sake of a fast prototyping of algorithms, the proposed
time-stepping scheme is implemented in the high-level programming
tool \texttt{MATLAB}.}, the matrix polynomial $\mathbf{Q}$ is factorized such that only
terms linear or quadratic in $\mathbf{A}$ need to be evaluated as
suggested in Ref.~\citep{Reusch1988} for the solution of parabolic
differential equations. Therefore, the solution of Eq.~\eqref{eq:PadeStepping}
is obtained by successively solving equations that have been set up
for individual terms. For structural dynamics, the equation of a term
is partitioned into two sets of equations. After decoupling, a system
of equations similar to that of Newmark's method \citep{ArticleNewmark1959}
is obtained. Note that the matrix $\mathbf{A}$---see Eq.~\eqref{eq:matrixA}---is
not explicitly constructed and no inversion of the mass matrix is
required. Moreover, if the stiffness matrix $\mathbf{K}$, the damping
matrix $\mathbf{C}$, and the mass matrix $\mathbf{M}$ are sparse,
the algorithm requires only sparse matrix operations. Instead of directly
solving Eq.~\eqref{eq:PadeStepping}, it is helpful to either add
or subtract $\mathbf{Q}\mathbf{z}_{n-1}$ from both sides of the equation
depending on the order $M$ of the Padé approximation. Hence, for
odd orders of $M$ we will solve 
\begin{equation}
\mathbf{Q}(\mathbf{z}_{n}+\mathbf{z}_{n-1})=(\mathbf{P}+\mathbf{Q})\mathbf{z}_{n-1}+\sum\limits _{k=0}^{p_{\mathrm{f}}}\mathbf{C}_{k}\tilde{\mathbf{F}}_{\mathrm{m}n}^{(k)}\,,\label{eq:PadeStepping2}
\end{equation}
while it is advantageous from the point of view of an efficient computational
implementation to solve 
\begin{equation}
\mathbf{Q}(\mathbf{z}_{n}-\mathbf{z}_{n-1})=(\mathbf{P}-\mathbf{Q})\mathbf{z}_{n-1}+\sum\limits _{k=0}^{p_{\mathrm{f}}}\mathbf{C}_{k}\tilde{\mathbf{F}}_{\mathrm{m}n}^{(k)}\label{eq:PadeStepping3}
\end{equation}
for even orders $M$. In this way, the terms of the highest order
of $\mathbf{A}$ in polynomials $\mathbf{P}$ and $\mathbf{Q}$ cancel---see
Eq.~\eqref{eq:Exp_M_PQ} for their definition---on the right-hand
side of the equations, simplifying the numerical operations needed.
It is easy to see that the same time-stepping algorithm can be used
for both cases with a suitable substitution of variables.

Following Eq.~\eqref{eq:Exp_QM_factor} taken from \ref{sec:Pade-approximation},
the $M$th degree polynomial $\mathbf{Q}$---recall that the dependencies
have been dropped for a simplified notation $\mathbf{Q}_{M}(\mathbf{A})$---can
factorized according to its roots as 
\begin{equation}
\mathbf{Q}=\left(r_{1}\mathbf{I}-\mathbf{A}\right)\left(r_{2}\mathbf{I}-\mathbf{A}\right)\ldots\left(r_{M}\mathbf{I}-\mathbf{A}\right)\,.\label{eq:Qfactor}
\end{equation}
Using Eq.~\eqref{eq:Qfactor}, Eq.~\eqref{eq:PadeStepping2} for
odd values of $M$ and Eq.~\eqref{eq:PadeStepping3} for even values
of $M$ are combined to a single expression 
\begin{equation}
\left(r_{1}\mathbf{I}-\mathbf{A}\right)\left(r_{2}\mathbf{I}-\mathbf{A}\right)\ldots\left(r_{M}\mathbf{I}-\mathbf{A}\right)\hat{\mathbf{z}}_{n}=\mathbf{b}_{n}\quad\textrm{with }\begin{cases}
\hat{\mathbf{z}}_{n}=\mathbf{z}_{n}+\mathbf{z}_{n-1}, & \textrm{when \ensuremath{M} is odd}\\
\hat{\mathbf{z}}_{n}=\mathbf{z}_{n}-\mathbf{z}_{n-1}, & \textrm{when \ensuremath{M} is even}
\end{cases}\,,\label{eq:PadeSteppingEq}
\end{equation}
where the right-hand side is expressed as 
\begin{equation}
\mathbf{b}_{n}=\hat{\mathbf{P}}\mathbf{z}_{n-1}+\sum\limits _{k=0}^{p_{\mathrm{f}}}\mathbf{C}_{k}\tilde{\mathbf{F}}_{\mathrm{m}n}^{(k)}\quad\textrm{with }\begin{cases}
\hat{\mathbf{P}}=\mathbf{P}+\mathbf{Q}, & \textrm{when \ensuremath{M} is odd}\\
\hat{\mathbf{P}}=\mathbf{P}-\mathbf{Q}, & \textrm{when \ensuremath{M} is even}
\end{cases}\,.\label{eq:PadeStepping_b}
\end{equation}
Equation~\eqref{eq:PadeSteppingEq} is reformulated as a system of
equations which are only linear in the matrix $\mathbf{A}$ by solving
for auxiliary variables $\mathbf{z}^{(k)}$, where $k\in\{1,2,\ldots,M-1\}$.
\begin{equation}
\begin{split}\left(r_{1}\mathbf{I}-\mathbf{A}\right)\mathbf{z}^{(1)} & =\mathbf{b}_{n}\,,\\
\left(r_{2}\mathbf{I}-\mathbf{A}\right)\mathbf{z}^{(2)} & =\mathbf{z}^{(1)}\,,\\
\cdots\\
\left(r_{M}\mathbf{I}-\mathbf{A}\right)\hat{\mathbf{z}}_{n} & =\mathbf{z}^{(M-1)}\,.
\end{split}
\label{eq:SteppingSuccessive}
\end{equation}
For the sake of completeness, the definition of the auxiliary variables
is provided at this point 
\begin{equation}
\mathbf{z}^{(k)}=\left(r_{k+1}\mathbf{I}-\mathbf{A}\right)\left(r_{k+2}\mathbf{I}-\mathbf{A}\right)\ldots\left(r_{M}\mathbf{I}-\mathbf{A}\right)\hat{\mathbf{z}}_{n}\quad\;\forall k=1,2,\ldots,M-1\,.
\end{equation}
Equation~\eqref{eq:SteppingSuccessive} can be solved successively
starting from the first line with the known right-hand side $\mathbf{b}_{n}$
and then working through all $M$ equations. Note that the polynomial
$\mathbf{Q}$ may have two types of roots: (i) single real roots and
(ii) pairs of complex conjugate roots. Their treatments are discussed
below in Sects.~\ref{subsec:realRoot} and \ref{subsec:complexRoot},
respectively.

The force-related summation term at the right-hand side of Eq.~\eqref{eq:PadeStepping}
combined with Eq.~\eqref{eq:PadeStepping_C} essentially only involves
the product of the system matrix $\mathbf{A}$ with some vectors.
The operation\footnote{The definitions of the vectors $\mathbf{y}$ and $\mathbf{x}$ is
strictly limited to this paragraph and will not be used in the remainder
of this section.} is denoted as 
\begin{equation}
\mathbf{y}=\mathbf{A}\mathbf{x}\label{eq:Ax_def}
\end{equation}
with partitions conforming to the size of $\mathbf{A}$ 
\begin{equation}
\mathbf{y}=\begin{Bmatrix}\mathbf{y}_{1}\\
\mathbf{y}_{2}
\end{Bmatrix}\qquad\text{and}\qquad\mathbf{x}=\begin{Bmatrix}\mathbf{x}_{1}\\
\mathbf{x}_{2}
\end{Bmatrix}\,.\label{eq:Ax_partition}
\end{equation}
Exploiting the definition of matrix $\mathbf{A},$ given in Eq.~\eqref{eq:matrixA},
the vector $\mathbf{y}$ is obtained from the following expressions
\begin{subequations}
\label{eq:Ax} 
\begin{align}
\mathbf{M}\mathbf{y}_{1} & =-\Delta t\mathbf{C}\mathbf{x}_{1}-\Delta t^{2}\mathbf{K}\mathbf{x}_{2}\,,\label{eq:Ax1}\\
\mathbf{y}_{2} & =\mathbf{x}_{1}\,.\label{eq:Ax2}
\end{align}
\end{subequations}
 This algorithm avoids the explicit construction of the matrix $\mathbf{A}$.
Products involving a higher integer power of $\mathbf{A}$ with a
vector can be straightforwardly computed by applying Eq.~\eqref{eq:Ax}
repeatedly. 

\subsubsection{Real root case}

\label{subsec:realRoot} When a root of the matrix polynomial $\mathbf{Q}$
is a real number, the corresponding line in Eq.~\eqref{eq:SteppingSuccessive}
is denoted as 
\begin{equation}
\left(r\mathbf{I}-\mathbf{A}\right)\mathbf{x}=\mathbf{g}\,,\label{eq:realRootEq}
\end{equation}
where $r$ is the real root. The unknown vector $\mathbf{x}$ is determined
in relation to a given right-hand side denoted by $\mathbf{g}$\footnote{The definitions of the vectors $\mathbf{x}$ and $\mathbf{g}$ is
strictly limited to this paragraph and will not be used in the remainder
of this section.}. The two vectors are partitioned conforming to the size of the matrix
$\mathbf{A}$ as 
\begin{equation}
\mathbf{x}=\begin{Bmatrix}\mathbf{x}_{1}\\
\mathbf{x}_{2}
\end{Bmatrix}\qquad\text{and}\qquad\mathbf{g}=\begin{Bmatrix}\mathbf{g}_{1}\\
\mathbf{g}_{2}
\end{Bmatrix}\,.\label{eq:realRootPartition}
\end{equation}
Exploiting the definition of matrix $\mathbf{A}$ given in Eq.~\eqref{eq:matrixA}
and substituting Eq.~\eqref{eq:realRootPartition} into Eq.~\eqref{eq:realRootEq}
we arrive at 
\begin{equation}
\begin{bmatrix}r\mathbf{I}+\Delta t\mathbf{M}^{-1}\mathbf{C} & \Delta t^{2}\mathbf{M}^{-1}\mathbf{K}\\
-\mathbf{I} & r\mathbf{I}
\end{bmatrix}\begin{Bmatrix}\mathbf{x}_{1}\\
\mathbf{x}_{2}
\end{Bmatrix}=\begin{Bmatrix}\mathbf{g}_{1}\\
\mathbf{g}_{2}
\end{Bmatrix}\,.\label{eq:realRootEq2}
\end{equation}
Pre-multiplying the first row block with $r\mathbf{M}$ yields 
\begin{subequations}
\label{eq:realRootEq3} 
\begin{align}
\left(r^{2}\mathbf{M}+r\Delta t\mathbf{C}\right)\mathbf{x}_{1}+r\Delta t^{2}\mathbf{K}\mathbf{x}_{2} & =r\mathbf{M}\mathbf{g}_{1}\,,\label{eq:realRootEq3a}\\
-\mathbf{x}_{1}+r\mathbf{x}_{2} & =\mathbf{g}_{2}\,.\label{eq:realRootEq3b}
\end{align}
\end{subequations}
 In the next step, Eq.~\eqref{eq:realRootEq3a} is rearranged with
respect to $\mathbf{x}_{2}$ 
\begin{equation}
\mathbf{x}_{2}=\cfrac{1}{r}(\mathbf{x}_{1}+\mathbf{g}_{2})\,.\label{eq:realRootEq4}
\end{equation}
Substituting Eq.~\eqref{eq:realRootEq4} into Eq.~\eqref{eq:realRootEq3a}
and simplifying the resulting expression leads to 
\begin{equation}
\left(r^{2}\mathbf{M}+r\Delta t\mathbf{C}+\Delta t^{2}\mathbf{K}\right)\mathbf{x}_{1}=r\mathbf{M}\mathbf{g}_{1}-\Delta t^{2}\mathbf{K}\mathbf{g}_{2}\,.\label{eq:realRootEq5}
\end{equation}
After solving Eq.~\eqref{eq:realRootEq5} for $\mathbf{x}_{1}$ and
determining $\mathbf{x}_{2}$ using Eq.~\eqref{eq:realRootEq4},
the overall solution $\mathbf{x}$ of Eq.~\eqref{eq:realRootEq}
is obtained.

Note that Eq.~\eqref{eq:realRootEq5} is similar to the Newmark time-stepping
scheme, sharing several helpful properties such as the fact that (i)
no inverse of the mass matrix is required and (ii) the coefficient
matrix is symmetric and positive definite. Due to the sparse nature
of the stiffness, damping, and mass matrices resulting from most spatial
discretization methods, established sparse matrix algorithms can be
exploited. 

\subsubsection{Complex conjugate roots case}

\label{subsec:complexRoot} When the roots of the matrix polynomial
$\mathbf{Q}$ contain a pair of complex conjugate roots the treatment
of the equations needs to be adjusted. According to Ref.~\citep{Reusch1988},
it is advantageous to consider the pair of complex conjugate roots
together. Thus, the equation to be solved is expressed as 
\begin{equation}
\left(r\mathbf{I}-\mathbf{A}\right)\left(\bar{r}\mathbf{I}-\mathbf{A}\right)\mathbf{x}=\mathbf{g}\,,\label{eq:ComplexRootEq}
\end{equation}
where $r$ and $\bar{r}$ denote the pair of complex conjugate roots.
Both the given right-hand side $\mathbf{g}$ and the unknown vector
$\mathbf{x}$ are real-valued\footnote{The definitions of the vectors $\mathbf{x}$, $\mathbf{y}$ and $\mathbf{g}$
is strictly limited to this paragraph and will not be used in the
remainder of this section.}. Mathematically, the solution can be simply written as 
\begin{equation}
\mathbf{x}=\left[\left(r\mathbf{I}-\mathbf{A}\right)\left(\bar{r}\mathbf{I}-\mathbf{A}\right)\right]^{-1}\mathbf{g}\,.\label{eq:ComplexRootEq2}
\end{equation}
Expressing the inversion as partial fractions---see Eq.~\eqref{eq:PartialFraction},
the solution is formulated as 
\begin{equation}
\mathbf{x}=\cfrac{-1}{2\mathrm{Im}(r)\mathrm{i}}\left[\left(r\mathbf{I}-\mathbf{A}\right)^{-1}-\left(\bar{r}\mathbf{I}-\mathbf{A}\right)^{-1}\right]\mathbf{g}\,.\label{eq:ComplexRootPartialFraction}
\end{equation}
Introducing the auxiliary vector $\mathbf{y}$ 
\begin{equation}
\mathbf{y}=\left(r\mathbf{I}-\mathbf{A}\right)^{-1}\mathbf{g}\label{eq:ComplexRoot_yDef}
\end{equation}
with its complex conjugate 
\begin{equation}
\bar{\mathbf{y}}=\left(\bar{r}\mathbf{I}-\mathbf{A}\right)^{-1}\mathbf{g\,,}\label{eq:ComplexRoot_yConjugate}
\end{equation}
Eq.~\eqref{eq:ComplexRootPartialFraction} is rewritten as 
\begin{equation}
\mathbf{x}=\cfrac{-1}{2\mathrm{Im}(r)\mathrm{i}}\left(\mathbf{y}-\bar{\mathbf{y}}\right)=-\cfrac{\mathrm{Im}(\mathbf{y})}{\mathrm{Im}(r)}\,.\label{eq:ComplexRoot-slnx}
\end{equation}
The auxiliary vector $\mathbf{y}$ is obtained by solving 
\begin{equation}
\left(r\mathbf{I}-\mathbf{A}\right)\mathbf{y}=\mathbf{g}\,.\label{eq:ComplexRootEqy}
\end{equation}
Note that Eq.~\eqref{eq:ComplexRootEqy} which is used to determine
$\mathbf{y}$ is identical in its mathematical structure to Eq.~\eqref{eq:realRootEq}
which is used to calculate $\mathbf{x}$. Partitioning the vector
$\mathbf{y}$ in the same way 
\begin{equation}
\mathbf{y}=\begin{Bmatrix}\mathbf{y}_{1}\\
\mathbf{y}_{2}
\end{Bmatrix}\label{eq:ComplexRoot_sln_y}
\end{equation}
and following the derivation from Eqs.~\eqref{eq:realRootEq2} to
\eqref{eq:realRootEq5}, we can easily calculate the subvector $\mathbf{y}_{1}$
by solving 
\begin{equation}
\left(r^{2}\mathbf{M}+r\Delta t\mathbf{C}+\Delta t^{2}\mathbf{K}\right)\mathbf{y}_{1}=r\mathbf{M}\mathbf{g}_{1}-\Delta t^{2}\mathbf{K}\mathbf{g}_{2}\,,\label{eq:ComplexRoot_sln_y1}
\end{equation}
while the second subvector $\mathbf{y}_{2}$ is obtained from 
\begin{equation}
\mathbf{y}_{2}=\cfrac{1}{r}\left(\mathbf{y}_{1}+\mathbf{g}_{2}\right)\,.\label{eq:ComplexRoot_sln_y2}
\end{equation}
Equation~\eqref{eq:ComplexRoot_sln_y1} is again similar to the Newmark
time-stepping scheme. Since $r$ is a complex number the coefficient
matrix $\left(r^{2}\mathbf{M}+r\Delta t\mathbf{C}+\Delta t^{2}\mathbf{K}\right)$
of Eq.~\eqref{eq:ComplexRoot_sln_y1} is symmetric but not Hermitian.

For later use, the matrix-vector product $\mathbf{A}\mathbf{x}$ is
obtained utilizing Eqs.~\eqref{eq:ComplexRoot-slnx} and \eqref{eq:ComplexRootEqy}
and keeping in mind that the vector $\mathbf{g}$ is real-valued and
thus, its imaginary part is zero 
\begin{equation}
\mathbf{A}\mathbf{x}=-\cfrac{\mathrm{Im}(\mathbf{A}\mathbf{y})}{\mathrm{Im}(r)}=-\cfrac{\mathrm{Im}(r\mathbf{y}-\mathbf{g})}{\mathrm{Im}(r)}=-\cfrac{\mathrm{Im}(r\mathbf{y})}{\mathrm{Im}(r)}\,.\label{eq:ComplexRootAx}
\end{equation}
We can obtain the product $\mathbf{A}^{2}\mathbf{x}$ in a similar
way using Eq.~\eqref{eq:ComplexRootEq} 
\begin{equation}
\mathbf{A}^{2}\mathbf{x}=\mathbf{g}+2\mathrm{Re}(r)\mathbf{A}\mathbf{x}-r\bar{r}\mathbf{x}\label{eq:ComplexRootAAx}
\end{equation}
and consequently, no direct computation of the matrix-vector product
is required. 

\subsection{Implementation of the novel algorithm for second-, fourth-, sixth-,
and eighth-order accurate methods}

\label{subsec:Time-stepping-schemes} The equations necessary for
the computational implementation of the present time-stepping scheme
are summarized below for orders $M\,{=}\,1$ to $M\,{=}\,4$. The
diagonal Padé approximation of order $M$ can be easily derived from
the formulae provided in Eqs.~\eqref{eq:ExpmPade} and \eqref{eq:Exp_M_PQ}. 

\subsubsection{Second-order accurate time-stepping scheme}

\label{subsec:Order11} The second-order accurate implicit time integration
scheme is derived from a diagonal Padé expansion of order $M\,{=}\,1$.
According to the definition in Eq.~\eqref{eq:Exp_M_PQ}, the numerator
and denominator polynomials of the rational function are written as
\begin{subequations}
\label{eq:Exp_M1-PQ} 
\begin{align}
\mathbf{P} & =2\mathbf{I}+\mathbf{A}\,,\label{eq:Exp_M1-P}\\
\mathbf{Q} & =2\mathbf{I}-\mathbf{A}\,.\label{eq:Exp_M1-Q}
\end{align}
\end{subequations}
 The root of the polynomial $\mathbf{Q}$ is 
\begin{equation}
r=2\label{eq:Exp_M1-r}
\end{equation}
and therefore, the case of a single real root, discussed in Sect.~\ref{subsec:realRoot},
is encountered. Based on the results from Eq.~\eqref{eq:Exp_M1-PQ},
Eq.~\eqref{eq:PadeStepping_C} can be evaluated and the values for
$\mathbf{C}_{0}$ and $\mathbf{C}_{1}$ are determined 
\begin{subequations}
\label{eq:PadeStepping_C_M1} 
\begin{align}
\mathbf{C}_{0} & =2\mathbf{I}\,,\label{eq:PadeStepping_C0_M1}\\
\mathbf{C}_{1} & =\mathbf{0}\,.\label{eq:PadeStepping_C1_M1}
\end{align}
\end{subequations}
 Here, it is assumed that the force varies linearly with time in each
time step, i.e., $p_{\mathrm{f}}\,{=}\,1$. This assumption is justified
as it sufficient to ensure optimal convergence of the time-stepping
scheme and therefore, a high-order approximation of the actual force
is not required and would only lead to an unnecessary numerical overhead\footnote{Remark: The force term $\mathbf{F}_{n}$ has to be approximated with
the same polynomial order as the denominator and numerator polynomials
$\mathbf{Q}$ and $\mathbf{P}$ to ensure optimal convergence. Thus,
the specific choice $p_{\mathrm{f}}\,{=}\,M\,{=}\,L$ is recommended.}. In the next step, the expressions obtained in Eqs.~\eqref{eq:Exp_M1-PQ}--\eqref{eq:PadeStepping_C_M1}
are substituted into the time-stepping scheme formulated in Eq.~\eqref{eq:PadeSteppingEq}.
This results in 
\begin{equation}
\left(2\mathbf{I}-\mathbf{A}\right)\left(\mathbf{z}_{n}+\mathbf{z}_{n-1}\right)=4\mathbf{z}_{n-1}+2\mathbf{F}_{\mathrm{m}n}\,.\label{eq:Exp_M1_stepping}
\end{equation}
By following the solution procedure, discussed in detail in Sect.~\ref{subsec:realRoot},
we can determine the auxiliary vector $\hat{\mathbf{z}}_{n}\,{=}\,(\mathbf{z}_{n}+\mathbf{z}_{n-1})$\footnote{Note that this is the definition of the auxiliary vector for odd orders
$M$---see Eq.~\eqref{eq:PadeSteppingEq}.} and hence, the state-space vector $\mathbf{z}_{n}$ at time $t_{n}$
is also known.

The present time integration scheme of order $M\,{=}\,1$ and $p_{\mathrm{f}}\,{=}\,1$
is equivalent to the average constant acceleration scheme of the Newmark
family of time integration methods with the two Newmark parameters
set to $\gamma=1/2$ and $\beta=1/4$ (commonly also referred to as
trapezoidal rule). To illustrate the equivalence of both schemes,
the basic idea is to evaluate the semi-discrete equations of motion
at times $t_{n}$ and $t_{n-1}$ and average the results that are
obtained by Newmark's algorithm. First, we can identify the following
correspondences by comparing Eq.~\eqref{eq:Exp_M1_stepping} to Eqs.~\eqref{eq:realRootEq}
and \eqref{eq:realRootPartition}
\begin{equation}
\begin{split}\mathbf{x}_{1} & =\overset{\CIRC}{\mathbf{u}}_{n}+\overset{\CIRC}{\mathbf{u}}_{n-1}\,,\\
\mathbf{x}_{2} & =\mathbf{u}_{n}+\mathbf{u}_{n-1}\,,\\
\mathbf{g}_{1} & =4\overset{\CIRC}{\mathbf{u}}_{n-1}+2\Delta t^{2}\mathbf{M}^{-1}\mathbf{f}_{\text{m}n}\,,\\
\mathbf{g}_{2} & =4\mathbf{u}_{n-1}\,.
\end{split}
\label{eq:Exp_M1_terms}
\end{equation}
For the definition of the forcing terms$\mathbf{f}$ (related to the
original second-order ODE) and $\mathbf{F}$(state-space formulation)
please refer to Eq.~\eqref{eq:vectorf}. Consequently, Eq.~\eqref{eq:realRootEq5}
is expressed as 
\begin{equation}
\left(4\mathbf{M}+2\Delta t\mathbf{C}+\Delta t^{2}\mathbf{K}\right)\left(\overset{\CIRC}{\mathbf{u}}_{n}+\overset{\CIRC}{\mathbf{u}}_{n-1}\right)=4\Delta t^{2}\mathbf{f}_{\text{m}n}+8\mathbf{M}\overset{\CIRC}{\mathbf{u}}_{n-1}-4\Delta t^{2}\mathbf{K}\mathbf{u}_{n-1}\label{eq:Exp_M1_stepping2}
\end{equation}
from which the velocity $\overset{\CIRC}{\mathbf{u}}_{n}$ is determined.
In the last step, Eq.~\eqref{eq:ComplexRoot_sln_y2} is solved, which
yields the displacement $\mathbf{u}_{n}$ at time $t_{n}$ 
\begin{equation}
\mathbf{u}_{n}=\mathbf{u}_{n-1}+\cfrac{1}{2}\left(\overset{\CIRC}{\mathbf{u}}_{n}+\overset{\CIRC}{\mathbf{u}}_{n-1}\right)\,.\label{eq:Exp_M1_stepping3}
\end{equation}

In order to show its equivalence with the novel time-stepping scheme
proposed above, the Newmark's constant average acceleration method
is formulated at two time instances $t_{n}$ and $t_{n-1}$. We recall
that the expressions for the velocity and displacement vectors at
time $t_{n}$ given in terms of the dimensionless time $s$ are expressed
as
\begin{subequations}
\label{eq:Newmark_uv} 
\begin{align}
\overset{\CIRC}{\mathbf{u}}_{n} & =\overset{\CIRC}{\mathbf{u}}_{n-1}+0.5\overset{\CCIRC}{\mathbf{u}}_{n-1}+0.5\overset{\CCIRC}{\mathbf{u}}_{n}\,,\label{eq:Newmark_u}\\
\overset{}{\mathbf{u}}_{n} & =\overset{}{\mathbf{u}}_{n-1}+\overset{\CIRC}{\mathbf{u}}_{n-1}+0.25\overset{\CCIRC}{\mathbf{u}}_{n-1}+0.25\overset{\CCIRC}{\mathbf{u}}_{n}\,.\label{eq:Newmark_v}
\end{align}
\end{subequations}
 Evaluating Eq.~\eqref{eq:equationOfMotionDimensionless} at the
time $t_{n}$ gives 
\begin{equation}
\mathbf{M}\overset{\CCIRC}{\mathbf{u}}_{n}+\Delta t\mathbf{C}\overset{\CIRC}{\mathbf{u}}_{n}+\Delta t^{2}\mathbf{K}\mathbf{u}_{n}=\Delta t^{2}\mathbf{f}_{n}\,.\label{eq:equationOfMotionDimensionless_n}
\end{equation}
In the next step, Eqs.~\eqref{eq:Newmark_uv} are substituted into
Eq.~\eqref{eq:equationOfMotionDimensionless_n} and all terms related
to the acceleration at time $t_{n}$ are moved to the left-hand side,
while all other terms are added/subtracted to/from the right-hand
side yielding 
\begin{equation}
\left(\mathbf{M}+0.5\Delta t\mathbf{C}+0.25\Delta t^{2}\mathbf{K}\right)\overset{\CCIRC}{\mathbf{u}}_{n}=\Delta t^{2}\mathbf{f}_{n}-\Delta t\mathbf{C}\left(\overset{\CIRC}{\mathbf{u}}_{n-1}+0.5\overset{\CCIRC}{\mathbf{u}}_{n-1}\right)-\Delta t^{2}\mathbf{K}\left(\mathbf{u}_{n-1}+\overset{\CIRC}{\mathbf{u}}_{n-1}+0.25\overset{\CCIRC}{\mathbf{u}}_{n-1}\right)\,.\label{eq:Newmark_Stepn}
\end{equation}
Expressing Eq.~\eqref{eq:Newmark_uv} in terms of the time $t_{n-1}$
leads to 
\begin{subequations}
\label{eq:Newmark_uv-1} 
\begin{align}
\overset{\CIRC}{\mathbf{u}}_{n-1} & =\overset{\CIRC}{\mathbf{u}}_{n-2}+0.5\overset{\CCIRC}{\mathbf{u}}_{n-2}+0.5\overset{\CCIRC}{\mathbf{u}}_{n-1}\,,\label{eq:Newmark_u-1}\\
\overset{}{\mathbf{u}}_{n-1} & =\overset{}{\mathbf{u}}_{n-2}+\overset{\CIRC}{\mathbf{u}}_{n-2}+0.25\overset{\CCIRC}{\mathbf{u}}_{n-2}+0.25\overset{\CCIRC}{\mathbf{u}}_{n-1}\,.\label{eq:Newmark_v-1}
\end{align}
\end{subequations}
 This time, Eq.~\eqref{eq:equationOfMotionDimensionless} is evaluated
at the time $t_{n-1}$ 
\begin{equation}
\mathbf{M}\overset{\CCIRC}{\mathbf{u}}_{n-1}+\Delta t\mathbf{C}\overset{\CIRC}{\mathbf{u}}_{n-1}+\Delta t^{2}\mathbf{K}\mathbf{u}_{n-1}=\Delta t^{2}\mathbf{f}_{n-1}\,.\label{eq:equationOfMotionDimensionless_n-1}
\end{equation}
Substituting Eqs.~\eqref{eq:Newmark_uv-1} into Eq.~\eqref{eq:equationOfMotionDimensionless_n-1}
yields an expression analogous to Eq.~\eqref{eq:Newmark_Stepn} at
time $t_{n-1}$ 
\begin{equation}
\left(\mathbf{M}+0.5\Delta t\mathbf{C}+0.25\Delta t^{2}\mathbf{K}\right)\overset{\CCIRC}{\mathbf{u}}_{n-1}=\Delta t^{2}\mathbf{f}_{n-1}-\Delta t\mathbf{C}\underbrace{\left(\overset{\CIRC}{\mathbf{u}}_{n-1}-0.5\overset{\CCIRC}{\mathbf{u}}_{n-1}\right)}_{=\;\overset{\CIRC}{\mathbf{u}}_{n-2}\;+\;0.5\overset{\CCIRC}{\mathbf{u}}_{n-2}}-\Delta t^{2}\mathbf{K}\underbrace{\left(\mathbf{u}_{n-1}-0.25\overset{\CCIRC}{\mathbf{u}}_{n-1}\right)}_{=\;\overset{}{\mathbf{u}}_{n-2}\;+\;\overset{\CIRC}{\mathbf{u}}_{n-2}\;+\;0.25\overset{\CCIRC}{\mathbf{u}}_{n-2}}\,.\label{eq:Newmark_Stepn1}
\end{equation}
Averaging Eqs.~\eqref{eq:Newmark_Stepn} and \eqref{eq:Newmark_Stepn1}
and eliminating the average acceleration using Eq.~\eqref{eq:Newmark_v}
yields 
\begin{equation}
\left(\mathbf{M}+0.5\Delta t\mathbf{C}+0.25\Delta t^{2}\mathbf{K}\right)\underbrace{\left(\overset{\CIRC}{\mathbf{u}}_{n}-\overset{\CIRC}{\mathbf{u}}_{n-1}\right)}_{=\;\sfrac{1}{2}(\overset{\CCIRC}{\mathbf{u}}_{n-1}\;+\;\overset{\CCIRC}{\mathbf{u}}_{n})}=\Delta t^{2}\mathbf{f}_{\text{m}n}-\Delta t\mathbf{C}\overset{\CIRC}{\mathbf{u}}_{n-1}-\Delta t^{2}\mathbf{K}\left(\mathbf{u}_{n-1}+0.5\overset{\CIRC}{\mathbf{u}}_{n-1}\right)\,,\label{eq:Newmark_Stepn-1}
\end{equation}
where the effective force term is defined as $\mathbf{f}_{\text{m}n}=(\mathbf{f}_{n-1}+\mathbf{f}_{n})/2$.
In the next step, Eq.~\eqref{eq:Newmark_Stepn-1} is further manipulated
by adding the term 2$(\mathbf{M}+0.5\Delta t\mathbf{C}+0.25\Delta t^{2}\mathbf{K})\overset{\CIRC}{\mathbf{u}}_{n-1}$
on both sides of Eq.~\eqref{eq:Newmark_Stepn-1}. This operation
leads to 
\begin{equation}
(\mathbf{M}+0.5\Delta t\mathbf{C}+0.25\Delta t^{2}\mathbf{K})(\overset{\CIRC}{\mathbf{u}}_{n}+\overset{\CIRC}{\mathbf{u}}_{n-1})=\Delta t^{2}\mathbf{f}_{\text{m}n}+2\mathbf{M}\overset{\CIRC}{\mathbf{u}}_{n-1}-\Delta t^{2}\mathbf{K}\mathbf{u}_{n-1}\,,\label{eq:Newmark_Stepn-2}
\end{equation}
which is identical to Eq.~\eqref{eq:Exp_M1_stepping2} divided by
4, highlighting that the average constant acceleration method is mathematically
identical to the novel time integration scheme for a diagonal Padé
expansion of order $M\,{=}\,1$ of the matrix exponential and a linear
approximation of the force vector within a time step, i.e., $p_{\mathrm{f}}\,{=}\,1$.

\subsubsection{Fourth-order accurate time-stepping scheme}

\label{subsec:Order22} The fourth-order accurate implicit time integration
scheme is derived from a diagonal Padé expansion of order $M\,{=}\,2$.
According to the definition in Eq.~\eqref{eq:Exp_M_PQ}, the numerator
and denominator polynomials of the rational function are written as
\begin{subequations}
\label{eq:Exp_M2-PQ} 
\begin{align}
\mathbf{P} & =12\mathbf{I}+6\mathbf{A}+\mathbf{A}^{2}\,,\label{eq:Exp_M2-P}\\
\mathbf{Q} & =12\mathbf{I}-6\mathbf{A}+\mathbf{A}^{2}\,.\label{eq:Exp_M2-Q}
\end{align}
\end{subequations}
 The polynomial $\mathbf{Q}$ has the roots 
\begin{equation}
r_{1}=3+\sqrt{3}\mathrm{i}\qquad\text{and}\qquad\bar{r}_{1}=3-\sqrt{3}\mathrm{i}\,.\label{eq:Exp_M2-r}
\end{equation}
Therefore, the case of a pair of complex conjugate roots, discussed
in Sect.~\ref{subsec:complexRoot}, is encountered. Based on the
results from Eq.~\eqref{eq:Exp_M2-PQ}, Eq.~\eqref{eq:PadeStepping_C}
can be evaluated and the values for $\mathbf{C}_{0}$, $\mathbf{C}_{1}$,
and $\mathbf{C}_{2}$ are determined 
\begin{subequations}
\label{eq:PadeStepping_C_M2} 
\begin{align}
\mathbf{C}_{0} & =12\mathbf{I}\,,\label{eq:PadeStepping_C0_M2}\\
\mathbf{C}_{1} & =-\mathbf{A}\,,\label{eq:PadeStepping_C1_M2}\\
\mathbf{C}_{2} & =\mathbf{I}\,.\label{eq:PadeStepping_C2_M2}
\end{align}
\end{subequations}
 Here, it is assumed that the force varies according to a quadratic
polynomial with time in each time step, i.e., $p_{\mathrm{f}}\,{=}\,2$.
In the next step, the expressions obtained in Eqs.~\eqref{eq:Exp_M2-PQ}--\eqref{eq:PadeStepping_C_M2}
are substituted into the time-stepping scheme formulated in Eq.~\eqref{eq:PadeSteppingEq}.
This results in 
\begin{equation}
\left(r_{1}\mathbf{I}-\mathbf{A}\right)\left(\bar{r}_{1}\mathbf{I}-\mathbf{A}\right)\left(\mathbf{z}_{n}-\mathbf{z}_{n-1}\right)=\mathbf{b}_{n}\label{eq:Exp_M2_stepping}
\end{equation}
with the right-hand side vector defined according to Eq.~\eqref{eq:PadeStepping_b}
\begin{equation}
\mathbf{b}_{n}=12\mathbf{A}\mathbf{z}_{n-1}+\sum\limits _{k=0}^{2}\mathbf{C}_{k}\tilde{\mathbf{F}}_{\mathrm{m}n}^{(k)}\,.\label{eq:Exp_M2_stepping_bn}
\end{equation}
Time-stepping is performed utilizing the procedure discussed in Sect.~\ref{subsec:complexRoot}
for a pair of complex conjugate roots. In each time step, the products
$\mathbf{A}(\mathbf{z}_{n}-\mathbf{z}_{n-1})$ and $\mathbf{A}\mathbf{z}_{n}$
are calculated employing the expression derived in Eq.~\eqref{eq:ComplexRootAx}.
The result of the product $\mathbf{A}\mathbf{z}_{n}$ at the time
$t_{n}$ will be stored for the next time step $n\,{+}\,1$, where
it can be employed in Eq.~\eqref{eq:Exp_M2_stepping_bn}---as ``$\mathbf{A}\mathbf{z}_{n-1}$''---to
calculate the right-hand side of Eq.~\eqref{eq:Exp_M2_stepping}.
The additional products including the matrix $\mathbf{A}$ (implicitly
contained in the definition of matrices $\mathbf{C}_{k}$) and the
force terms $\tilde{\mathbf{F}}_{\mathrm{m}n}^{(k)}$ are computed
based on the procedure derived in Eq.~\eqref{eq:Ax}. 

\subsubsection{Sixth-order accurate time-stepping scheme}

\label{subsec:Order33} The sixth-order accurate implicit time integration
scheme is derived from a diagonal Padé expansion of order $M\,{=}\,3$.
According to the definition in Eq.~\eqref{eq:Exp_M_PQ}, the numerator
and denominator polynomials of the rational function are written as
\begin{subequations}
\label{eq:Exp_M3-PQ} 
\begin{align}
\mathbf{P} & =120\mathbf{I}+60\mathbf{A}+12\mathbf{A}^{2}+\mathbf{A}^{3}\,,\label{eq:Exp_M3-P}\\
\mathbf{Q} & =120\mathbf{I}-60\mathbf{A}+12\mathbf{A}^{2}-\mathbf{A}^{3}\,.\label{eq:Exp_M3-Q}
\end{align}
\end{subequations}
 The polynomial $\mathbf{Q}$ has a single real root and a pair of
complex conjugate roots 
\begin{equation}
\begin{split}r_{1} & =4.644370709252171\,,\\
r_{2} & =3.677814645373910+3.508761919567443\mathrm{i}\,,\qquad\text{and}\\
\bar{r}_{2} & =3.677814645373910-3.508761919567443\mathrm{i}\,.
\end{split}
\label{eq:Exp_M3-r}
\end{equation}
Therefore, a combination of the two cases discussed in Sects.~\ref{subsec:realRoot}
and \ref{subsec:complexRoot}, is encountered. Using Eq.~\eqref{eq:Exp_M3-PQ},
Eq.~\eqref{eq:PadeStepping_C} can be evaluated and the values for
$\mathbf{C}_{0}$ to $\mathbf{C}_{3}$ are determined 
\begin{subequations}
\label{eq:PadeStepping_C_M3} 
\begin{align}
\mathbf{C}_{0} & =120\mathbf{I}+2\mathbf{A}^{2}\,,\label{eq:PadeStepping_C0_M3}\\
\mathbf{C}_{1} & =-10\mathbf{A}\,,\label{eq:PadeStepping_C1_M3}\\
\mathbf{C}_{2} & =10\mathbf{I}+\cfrac{1}{2}\mathbf{A}^{2}\,,\label{eq:PadeStepping_C2_M3}\\
\mathbf{C}_{3} & =-\cfrac{3}{2}\mathbf{A}\,.\label{eq:PadeStepping_C3_M3}
\end{align}
\end{subequations}
 Here, it is assumed that the force varies according to a cubic polynomial
with time in each time step, i.e., $p_{\mathrm{f}}\,{=}\,3$. In the
next step, the expressions obtained in Eqs.~\eqref{eq:Exp_M3-PQ}--\eqref{eq:PadeStepping_C_M3}
are substituted into the time-stepping scheme formulated in Eq.~\eqref{eq:PadeSteppingEq}.
This results in 
\begin{equation}
\left(r_{1}\mathbf{I}-\mathbf{A}\right)\left(r_{2}\mathbf{I}-\mathbf{A}\right)\left(\bar{r}_{2}\mathbf{I}-\mathbf{A}\right)\left(\mathbf{z}_{n}+\mathbf{z}_{n-1}\right)=\mathbf{b}_{n}\label{eq:Exp_M3_stepping}
\end{equation}
with the right-hand side vector defined according to Eq.~\eqref{eq:PadeStepping_b}
\begin{equation}
\mathbf{b}_{n}=\left(240\mathbf{I}+24\mathbf{A}^{2}\right)\mathbf{z}_{n-1}+\sum\limits _{k=0}^{3}\mathbf{C}_{k}\tilde{\mathbf{F}}_{\mathrm{m}n}^{(k)}\,.\label{eq:Exp_M3_stepping_bn}
\end{equation}
Equation~\eqref{eq:Exp_M3_stepping} is solved according to the procedure
detailed in Eq.~\eqref{eq:SteppingSuccessive} and hence, it is rewritten
as 
\begin{equation}
\begin{split}\left(r_{1}\mathbf{I}-\mathbf{A}\right)\mathbf{z}^{(1)} & =\mathbf{b}_{n}\,,\\
\left(r_{2}\mathbf{I}-\mathbf{A}\right)\left(\bar{r}_{2}\mathbf{I}-\mathbf{A}\right)\left(\mathbf{z}_{n}+\mathbf{z}_{n-1}\right) & =\mathbf{z}^{(1)}\,.
\end{split}
\label{eq:SteppingSuccessive-M3}
\end{equation}
The time-stepping is performed by solving the two equations successively
using the procedures explained in Sect~\ref{subsec:realRoot} and
\ref{subsec:complexRoot}. The matrix--vector products $\mathbf{A}\mathbf{z}_{n}$
and $\mathbf{A}^{2}\mathbf{z}_{n}$ are calculated via Eqs.~\eqref{eq:ComplexRootAx}
and \eqref{eq:ComplexRootAAx}, respectively. The obtained values
are again stored at the end of each time step $n$ to be re-used to
evaluate the right-hand side of Eq.~\eqref{eq:Exp_M3_stepping_bn}
for the next time step $n\,{+}\,1$. 

\subsubsection{Eighth-order accurate time-stepping scheme}

\label{subsec:Order44} The eighth-order accurate implicit time integration
scheme is derived from a diagonal Padé expansion of order $M\,{=}\,4$.
According to the definition in Eq.~\eqref{eq:Exp_M_PQ}, the numerator
and denominator polynomials of the rational function are written as
\begin{subequations}
\label{eq:Exp_M4-PQ} 
\begin{align}
\mathbf{P} & =1680\mathbf{I}+840\mathbf{A}+180\mathbf{A}^{2}+20\mathbf{A}^{3}+\mathbf{A}^{4}\,,\label{eq:Exp_M4-P}\\
\mathbf{Q} & =1680\mathbf{I}-840\mathbf{A}+180\mathbf{A}^{2}-20\mathbf{A}^{3}+\mathbf{A}^{4}\,.\label{eq:Exp_M4-Q}
\end{align}
\end{subequations}
 The polynomial $\mathbf{Q}$ has two pairs of complex conjugate roots
\begin{equation}
\begin{split}r_{1} & =4.207578794359259+5.314836083713504\mathrm{i}\,,\\
\bar{r}_{1} & =4.207578794359259-5.314836083713504\mathrm{i}\,,\\
r_{2} & =5.792421205640749+1.734468257869007\mathrm{i}\,,\qquad\text{and}\\
\bar{r}_{2} & =5.792421205640749-1.734468257869007\mathrm{i}\,.
\end{split}
\label{eq:Exp_M4-r}
\end{equation}
Therefore, the case of a pair of complex conjugate roots, discussed
in Sect.~\ref{subsec:complexRoot}, is encountered. Based on the
results from Eq.~\eqref{eq:Exp_M4-PQ}, Eq.~\eqref{eq:PadeStepping_C}
can be evaluated and the values for $\mathbf{C}_{0}$ to $\mathbf{C}_{4}$
are determined 
\begin{subequations}
\label{eq:PadeStepping_C_M4} 
\begin{align}
\mathbf{C}_{0} & =1680\mathbf{I}+40\mathbf{A}^{2}\,,\label{eq:PadeStepping_C0_M4}\\
\mathbf{C}_{1} & =-140\mathbf{A}-\mathbf{A}^{3}\,,\label{eq:PadeStepping_C1_M4}\\
\mathbf{C}_{2} & =140\mathbf{I}+8\mathbf{A}^{2}\,,\label{eq:PadeStepping_C2_M4}\\
\mathbf{C}_{3} & =-21\mathbf{A}-\cfrac{1}{4}\mathbf{A}^{3}\,,\label{eq:PadeStepping_C3_M4}\\
\mathbf{C}_{4} & =21\mathbf{I}+\cfrac{3}{2}\mathbf{A}^{2}\,.\label{eq:PadeStepping_C4_M4}
\end{align}
\end{subequations}
 Here, it is assumed that the force varies according to a quartic
polynomial with time in each time step, i.e., $p_{\mathrm{f}}\,{=}\,4$.
In the next step, the expressions obtained in Eqs.~\eqref{eq:Exp_M4-PQ}--\eqref{eq:PadeStepping_C_M4}
are substituted into the time-stepping scheme formulated in Eq.~\eqref{eq:PadeSteppingEq}.
This results in 
\begin{equation}
\left(r_{1}\mathbf{I}-\mathbf{A}\right)\left(\bar{r}_{1}\mathbf{I}-\mathbf{A}\right)\left(r_{2}\mathbf{I}-\mathbf{A}\right)\left(\bar{r}_{2}\mathbf{I}-\mathbf{A}\right)\left(\mathbf{z}_{n}-\mathbf{z}_{n-1}\right)=\mathbf{b}_{n}\label{eq:Exp_M4_stepping}
\end{equation}
with the right-hand side vector defined according to Eq.~\eqref{eq:PadeStepping_b}
\begin{equation}
\mathbf{b}_{n}=\left(1680\mathbf{A}+40\mathbf{A}^{3}\right)\mathbf{z}_{n-1}+\sum\limits _{k=0}^{4}\mathbf{C}_{k}\tilde{\mathbf{F}}_{\mathrm{m}n}^{(k)}\,.\label{eq:Exp_M4_stepping_bn}
\end{equation}
Equation~\eqref{eq:Exp_M4_stepping} is solved according to the procedure
detailed in Eq.~\eqref{eq:SteppingSuccessive} and hence, it is rewritten
as 
\begin{equation}
\begin{split}\left(r_{1}\mathbf{I}-\mathbf{A}\right)\left(\bar{r}_{1}\mathbf{I}-\mathbf{A}\right)\mathbf{z}^{(2)} & =\mathbf{b}_{n}\,,\\
\left(r_{2}\mathbf{I}-\mathbf{A}\right)\left(\bar{r}_{2}\mathbf{I}-\mathbf{A}\right)\left(\mathbf{z}_{n}-\mathbf{z}_{n-1}\right) & =\mathbf{z}^{(2)}\,.
\end{split}
\label{eq:SteppingSuccessive-M4}
\end{equation}
Time-stepping is performed by solving the two equations successively
using the procedure discussed in Sect.~\ref{subsec:complexRoot}.
The matrix--vector products $\mathbf{A}\mathbf{z}_{n}$ and $\mathbf{A}^{2}\mathbf{z}_{n}$
are calculated via Eqs.~\eqref{eq:ComplexRootAx} and \eqref{eq:ComplexRootAAx},
respectively. The obtained values are again stored at the end of each
time step $n$ to be re-used to evaluate the right-hand side of Eq.~\eqref{eq:Exp_M3_stepping_bn}
for the next time step $n\,{+}\,1$. In order to compute the matrix-vector
product $\mathbf{A}^{3}\mathbf{z}_{n}$, Eq.~\eqref{eq:SteppingSuccessive-M4}
is rewritten as 
\begin{equation}
\begin{split}\left(r_{1}\mathbf{I}-\mathbf{A}\right)\mathbf{z}^{(1)} & =\mathbf{b}_{n}\,,\\
\left(\bar{r}_{1}\mathbf{I}-\mathbf{A}\right)\left(r_{2}\mathbf{I}-\mathbf{A}\right)\left(\bar{r}_{2}\mathbf{I}-\mathbf{A}\right)\left(\mathbf{z}_{n}-\mathbf{z}_{n-1}\right) & =\mathbf{z}^{(1)}\,.
\end{split}
\label{eq:SteppingSuccessive-M4-AAAx1}
\end{equation}
The first equation which is used to determine $\mathbf{z}^{(1)}$
is also contained in the first equation of Eq.~\eqref{eq:SteppingSuccessive-M4},
which is required to calculate $\mathbf{z}^{(2)}$. Note that we make
use of the definition of $\mathbf{z}^{(1)}$ already indicated in
Eq.~\eqref{eq:SteppingSuccessive} and repeated at this point for
the sake of convenience 
\begin{equation}
\mathbf{z}^{(1)}=\left(\bar{r}_{1}\mathbf{I}-\mathbf{A}\right)\mathbf{z}^{(2)}=\bar{r}_{1}\mathbf{z}^{(2)}-\mathbf{A}\mathbf{z}^{(2)}\,.\label{eq:SteppingSuccessive-M4-z1z2}
\end{equation}
The solution of $\mathbf{z}^{(1)}$ is obtained using the procedure
suggested in Eq.~\eqref{eq:ComplexRootEqy} which is included in
the process of solving for $\mathbf{z}^{(2)}$. Consequently, no additional
computations are needed. Using the expression for $\mathbf{z}^{(2)}$
given in the second equation of Eq.~\eqref{eq:SteppingSuccessive-M4}
and substituting it for the second term in Eq.~\eqref{eq:SteppingSuccessive-M4-z1z2},
we can obtain a different variant of the second equation of Eq.~\eqref{eq:SteppingSuccessive-M4-AAAx1}
\begin{equation}
\bar{r}_{1}\mathbf{z}^{(2)}-\mathbf{A}\left(r_{2}\mathbf{I}-\mathbf{A}\right)\left(\bar{r}_{2}\mathbf{I}-\mathbf{A}\right)\left(\mathbf{z}_{n}-\mathbf{z}_{n-1}\right)=\mathbf{z}^{(1)}\,.\label{eq:SteppingSuccessive-M4-AAAx2}
\end{equation}
Expanding the second term and rearranging leads to 
\begin{equation}
\mathbf{A}^{3}\mathbf{z}_{n}=\mathbf{A}^{3}\mathbf{z}_{n-1}+\bar{r}_{1}\mathbf{z}^{(2)}-\mathbf{z}^{(1)}-\left(r_{2}\bar{r}_{2}\mathbf{A}-2\mathrm{Re}(r_{2})\mathbf{A}^{2}\right)\left(\mathbf{z}_{n}-\mathbf{z}_{n-1}\right)\,.\label{eq:SteppingSuccessive-M4-AAAx}
\end{equation}
Note that the evaluation of the matrix--vector product $\mathbf{A}^{3}\mathbf{z}_{n}$
is linear in the state-space and auxiliary vectors.

Employing the methodology discussed in this section allows us to derive
implicit and unconditionally stable time-stepping schemes of arbitrary
order of accuracy depending only on the order $M$ of the diagonal
Padé expansion of the matrix exponential function and the order $p_{\mathrm{f}}$
of the polynomial approximation of the force vector within each time
step. In contrast to other high-order time integration methods, the
inverse of the mass matrix is not required and the size of the system
of equations does not grow when elevating the order of accuracy. 

\section{Numerical examples}

\label{sec:Examples} In this section, the performance of the proposed
time integration scheme is investigated by means of several numerical
examples ranging from simple academic problems to more complex structures
of practical relevance. Here, we are particularly interested in the
accuracy and numerical costs of the novel scheme in comparison to
established time-stepping methods that are widely used in commercial
software packages. This provides an indication whether the increased
computational requirements are justified and advantages in terms of
efficiency can be leveraged.

The codes have been written in \texttt{MATLAB} version 9.8.0.1538580
(R2020a) Update 6 and the elapsed time is measured using the commands
\texttt{tic} and \texttt{toc}. Hence, we do not measure the actual
CPU time, but the physical time (wall-clock time) is recorded. Note
that only the time integration itself is timed, since the set-up of
the stiffness, damping, and mass matrices as well as the incorporation
of boundary conditions (Dirichlet and Neumann) are independent processes
(equal for all different time integrators) and not part of the actual
transient analysis.

Most of the computer time of the proposed high-order schemes is spent
on the solution of the systems of linear algebraic equations~\eqref{eq:realRootEq5}
for real roots and \eqref{eq:ComplexRoot_sln_y1} for pairs of complex
roots. In the current version, only direct solvers have been used.
In our \texttt{MATLAB}-implementation, the symmetric positive-definite
system given by Eq.~\eqref{eq:realRootEq5} (for the real root case)
is solved by means of the \texttt{decomposition}-command, which creates
reusable matrix decompositions (LU, LDL, Cholesky, QR, etc.) depending
on the properties of the input matrix. The performance was found to
be better compared to a direct call to the \texttt{chol} with permutation-
or \texttt{lu}-commands. The system of equations stated by Eq.~\eqref{eq:ComplexRoot_sln_y1},
obtained from a complex root, is symmetric but not Hermitian. The
\texttt{lu}-command for general matrices, which factorizes full or
sparse input matrices into an upper triangular matrix and a permuted
lower triangular matrix, is used. The option \texttt{`vector'} is
used to store the row and column permutation matrices. It is evident
that the comparisons on computer time reported in this paper are limited
to the above solvers. Iterative solvers and solvers that take advantage
of the symmetric structures of complex matrices will be investigated
and reported in future communications.

All simulation for this article (if not stated otherwise) are run
on a standard desktop workstation with the following specifications:
Precision 3630-Tower Workstation; Intel(R) Core(TM) i7-8700 CPU @
3.20$\,$GHz; 64$\,$GB DDR4 (4$\times$16$\,$GB), 2666 MHz; Intel(R)
UHD Graphics 630. 

\subsection{Single degree of freedom systems}

\label{sec:SDOF} In the wide body of literature, single-degree-of-freedom
(SDOF) systems are often deployed to study the accuracy and stability
properties of a time integration method. Studying SDOF systems is
of great importance as each multi-degree-of-freedom (MDOF) problem
can (theoretically) be transformed into a $N$-SDOF problems, where
$N$ denotes the number of DOFs. This can be achieved exploiting the
modal properties of a structure by employing mode decomposition techniques
\citep{BookChopra2019}. Moreover, since the numerical solution is
independent of the spatial discretization, the effect of the temporal
discretization can be isolated. Therefore, different cases considering
varying ICs, damping parameters, and external excitations are conveniently
investigated. The particular examples are taken from Ref.~\citep{ArticleKim2019d}
and will be discussed in detail in the remainder of this section.

\begin{figure}[t!]
\centering \subfloat[Time-domain signal $f_{1}(t)$\label{subfig:f1}]{\includegraphics{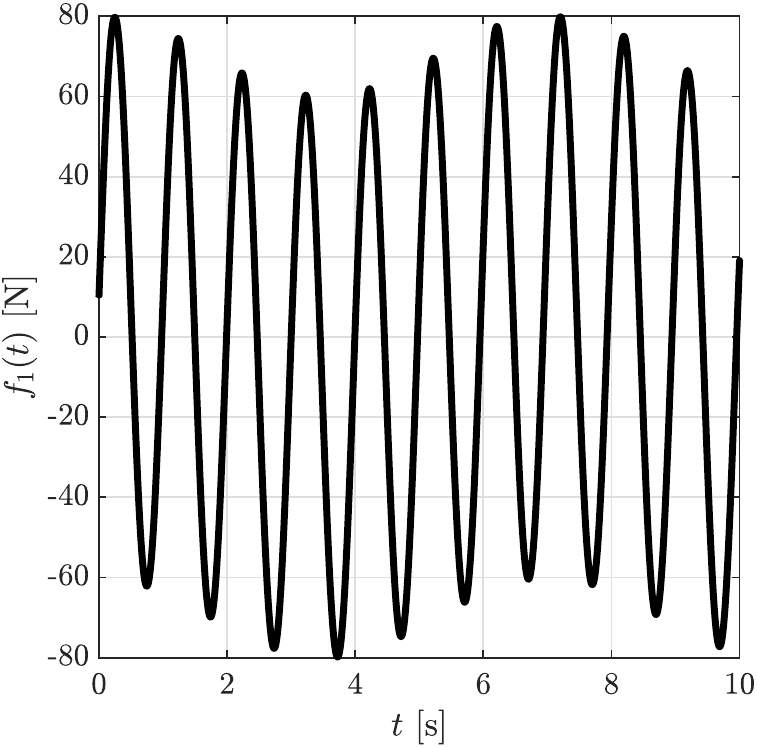}

}\hfill{}\subfloat[Time-domain signal $f_{2}(t)$\label{subfig:f2}]{\includegraphics{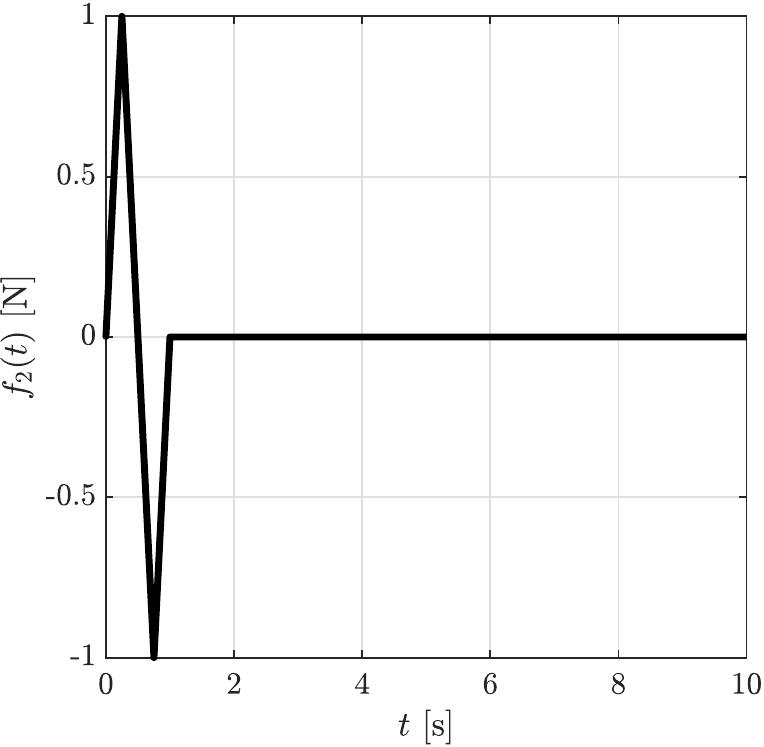}

}\\
 \subfloat[Frequency-domain signal $\hat{F}_{1}(f)$\label{subfig:F1}]{\includegraphics{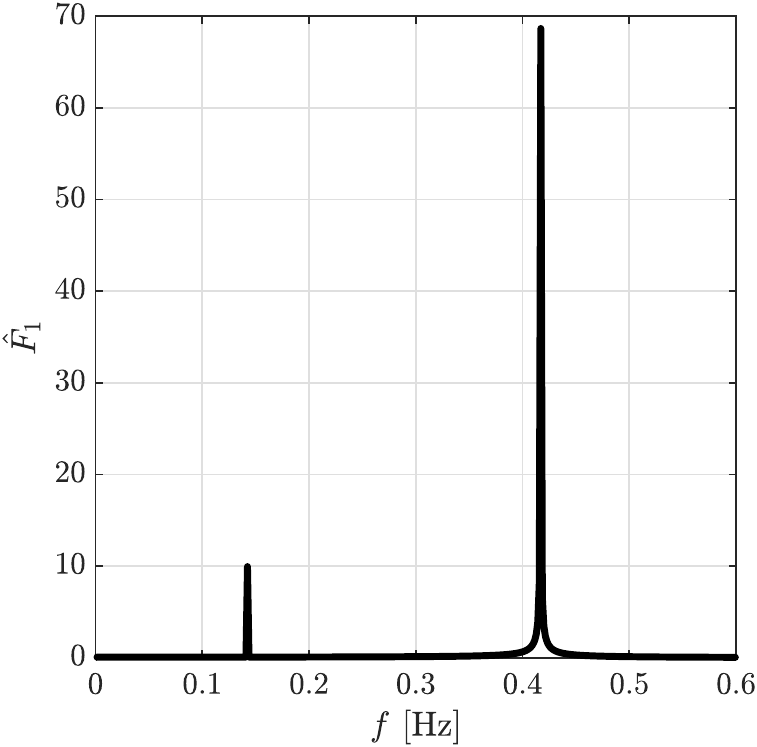}

}\hfill{}\subfloat[Frequency-domain signal $\hat{F}_{2}(f)$\label{subfig:F2}]{\includegraphics{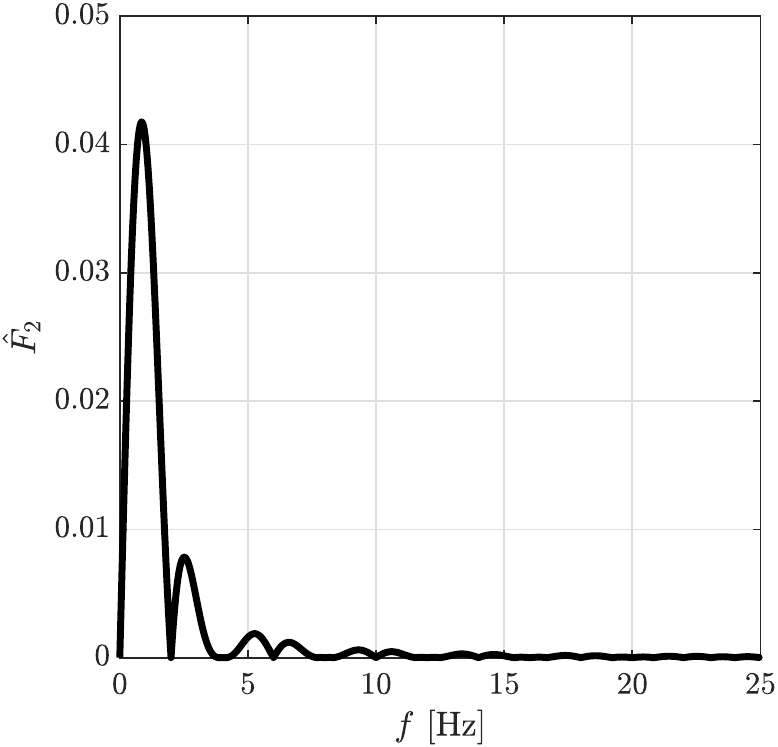}

}\caption{Excitation signal in the time and frequency domain---SDOF system.
\label{fig:ExcitationSignal}}
\end{figure}

The point of departure for the analysis of the proposed time integration
scheme is a generic SDOF system which can be written as 
\begin{equation}
\ddot{u}(t)+2\zeta\omega_{\mathrm{n}}\dot{u}(t)+\omega_{\mathrm{n}}^{2}u(t)=p(t)\,,\label{eq:ODE_SDOF}
\end{equation}
with the initial displacement and velocity 
\begin{equation}
u(t=0)=u_{0}\qquad\mathrm{and}\qquad\dot{u}(t=0)=\dot{u}_{0}\,.
\end{equation}
Here, $u(t)$ denotes the displacement, $\zeta$ stands for the damping
ratio, $\omega_{\mathrm{n}}$ is the natural (angular) frequency of
the structure, and $p(t)$ represents the external excitation force.
The stiffness of the system can be computed (if needed) from the mass
and the natural frequency of the structure as $k\,{=}\,\omega_{\mathrm{n}}^{2}m$,
where the mass is assumed to be 1$\,$kg for our examples. The damping
coefficient $c$ which is often used instead of $\zeta$ is defined
as $c\,{=}\,2\zeta\omega_{\mathrm{n}}m$. As mentioned before, different
cases covering a wide spectrum of applications are studied and the
selected parameters are listed in Table~\ref{tab:SDOF_parameters}.
The time-dependent amplitudes of the external excitation force are
either given as a harmonic function 
\begin{equation}
f_{1}(t)=\hat{f}_{1}\cos\left(\omega_{1}t\right)+\hat{f}_{2}\sin\left(\omega_{2}t\right)\,,\label{eq:f1}
\end{equation}
with $\hat{f}_{1}\,{=}\,10\,$N, $\hat{f}_{2}\,{=}\,70\,$N, $\omega_{1}\,{=}\,\sfrac{2\sqrt{5}}{5}\,$rad/s,
and $\omega_{2}\,{=}\,2\sqrt{10}\,$rad/s or as a piece-wise linear
one 
\begin{equation}
f_{2}(t)=\begin{cases}
4t\,, & 0\,\mathrm{s}\le t<0.25\,\mathrm{s}\\
-4t+2\,, & 0.25\,\mathrm{s}\le t<0.75\,\mathrm{s}\\
4t-4\,, & 0.75\,\mathrm{s}\le t<1\,\mathrm{s}\\
0\,, & t\ge1\,\mathrm{s}
\end{cases}\,.\label{eq:f2}
\end{equation}
The dynamic response of the SDOF system is analyzed for a certain
time interval $[0,t_{\mathrm{sim}}]$, where $t_{\mathrm{sim}}$ is
set to 10$\,$s in this section. The time history as well as the frequency
content of the two signals are depicted in Fig.~\ref{fig:ExcitationSignal}.
Here, we easily observe the two distinct peaks in the frequency spectrum
for the harmonic excitation with the sine- and cosine-terms, while
the spectrum of the triangular impulses contains much higher frequency
components and hence, can be considered as a broad-band excitation.
To obtain the frequency spectra, the discrete Fourier transform (DFT)
algorithm has been employed and the signal has been truncated after
$t_{\mathrm{sim}}$. For the application of the DFT algorithms, the
amplitude has been set to zero for $t\,{>}\,t_{\mathrm{sim}}$, which
explains why two extended peaks are observed in Fig.~\ref{subfig:F1}
instead of two discrete lines. 
\begin{table}[t!]
\centering \caption{Parameters for the different single degree of freedom systems. \label{tab:SDOF_parameters}}
\begin{tabular}{cccccc}
\toprule 
\multirow{1}{*}{Case} & Natural freq.  & Damping ratio  & Ext. excitation  & Init. disp.  & Init. velo.\tabularnewline
 & {[}rad/s{]}  & {[}--{]} & {[}N{]}  & {[}m{]}  & {[}m/s{]}\tabularnewline
\midrule 
\#1  & $2\pi$  & 0  & 0  & 1  & 0\tabularnewline
\#2  & $2\pi$  & 0  & 0  & 0  & $2\pi$\tabularnewline
\#3  & $2\pi$  & 0.05  & 0  & 1  & 0\tabularnewline
\#4  & $2\pi$  & 0.05  & 0  & 0  & $2\pi$\tabularnewline
\#5  & $2\pi$  & 0  & $f_{1}(t)$  & 2  & $\pi/3$\tabularnewline
\#6  & $2\pi$  & 0  & $f_{2}(t)$  & 2  & $\pi/3$\tabularnewline
\bottomrule
\end{tabular}
\end{table}

An advantage of studying simple SDOF systems is that it is possible
to derive analytical solutions in closed-form and consequently, a
suitable error measure is easily defined. In the following, the error
in displacements based on the $L_{2}$-norm ($\epsilon_{L_{2}}$)
is used to assess the accuracy of the time integration algorithms
\begin{equation}
\epsilon_{L_{2}}=\cfrac{\int\limits _{t=0s}^{t_{\mathrm{sim}}}(u_{\mathrm{ref}}(t)-u_{\mathrm{num}}(t))^{2}\;\mathrm{d}t}{\int\limits _{t=0s}^{t_{\mathrm{sim}}}u_{\mathrm{ref}}(t)^{2}\;\mathrm{d}t}\times100[\%]\,,\label{eq:ErrorL2}
\end{equation}
where $u_{\mathrm{ref}}(t)$ denotes the analytical (theoretical)
solution and $u_{\mathrm{num}}(t)$ represents the numerical solution
obtained utilizing a specific time-stepping algorithm.

At this point, let us recall the analytical solution for the free
vibration response of undamped and damped SDOF systems \citep{BookChopra2019}.
In this particular case, no external excitation force is acting on
the structure and therefore, the ICs determine the structural dynamic
behavior. The solution to the second-order ODE---Eq.~\eqref{eq:ODE_SDOF}---for
$\zeta\,{=}\,0$ is 
\begin{equation}
u(t)=u_{0}\cos(\omega_{\mathrm{n}}t)+\cfrac{\dot{u}_{0}}{\omega_{\mathrm{n}}}\,\sin(\omega_{\mathrm{n}}t)\,,\label{eq:FreeVibrationsUndamped}
\end{equation}
whereas in the damped case, i.e., $\zeta\ne0$, the theoretical solution
takes the following form 
\begin{equation}
u(t)=\exp(-\zeta\omega_{\mathrm{n}}t)\left(u_{0}\cos(\omega_{\mathrm{D}}t)+\cfrac{\dot{u}_{0}+\zeta\omega_{\mathrm{D}}u_{0}}{\omega_{\mathrm{D}}}\,\sin(\omega_{\mathrm{D}}t)\right)\,,\label{eq:FreeVibrationsDamped}
\end{equation}
where $\omega_{\mathrm{D}}$ is the natural (angular) frequency of
the damped systems, defined as $\omega_{\mathrm{D}}\,{=}\,\omega_{\mathrm{n}}\sqrt{1-\zeta^{2}}$.
Considering an external excitation, the solution is a bit more complex
since now not only the complementary (transient response) but also
the particular part (steady-state response) of the solution have to
be taken into account. The latter one is solely related to the non-homogeneous
ODE and therefore, determined by the excitation function. The general
solution for a harmonic excitation of an undamped system featuring
both sine- and cosine-terms with two distinct frequencies is given
in closed-form as 
\begin{equation}
\begin{split}u(t)= & \left(u_{0}-\cfrac{\hat{f}_{1}}{k}\,\cfrac{1}{1-\left(\sfrac{\omega_{1}}{\omega_{\mathrm{n}}}\right)^{2}}\right)\cos(\omega_{\mathrm{n}}t)+\left(\cfrac{\dot{u}_{0}}{\omega_{\mathrm{n}}}-\cfrac{\hat{f}_{2}}{k}\,\cfrac{\sfrac{\omega_{2}}{\omega_{\mathrm{n}}}}{1-\left(\sfrac{\omega_{2}}{\omega_{\mathrm{n}}}\right)^{2}}\right)\sin(\omega_{\mathrm{n}}t)+\\
 & \left(\cfrac{\hat{f}_{1}}{k}\,\cfrac{1}{1-\left(\sfrac{\omega_{1}}{\omega_{\mathrm{n}}}\right)^{2}}\right)\cos(\omega_{1}t)+\left(\cfrac{\hat{f}_{2}}{k}\,\cfrac{1}{1-\left(\sfrac{\omega_{2}}{\omega_{\mathrm{n}}}\right)^{2}}\right)\sin(\omega_{2}t)\,.
\end{split}
\label{eq:FrorcedVibrationsDamped}
\end{equation}
In order to assess the accuracy of the proposed time integrator in
relation to established methods, Newmark's method \citep{ArticleNewmark1959}
and Bathe's method \citep{ArticleBathe2005,ArticleNoh2019b} are included
in the analysis. Both methods are second-order accurate, whereas Newmark's
method is a single-step scheme and Bathe's method is a composite scheme
consisting of two sub-steps. Due to the simplicity of the investigated
systems, only the accuracy is considered, while the computational
costs are essentially negligible.

\begin{figure}[p]
\centering \subfloat{\includegraphics{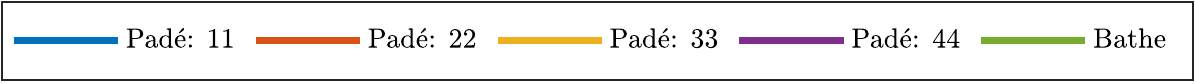}}\\
 \setcounter{subfigure}{0} \subfloat[Case \#1\label{subfig:SDOF_1}]{\includegraphics{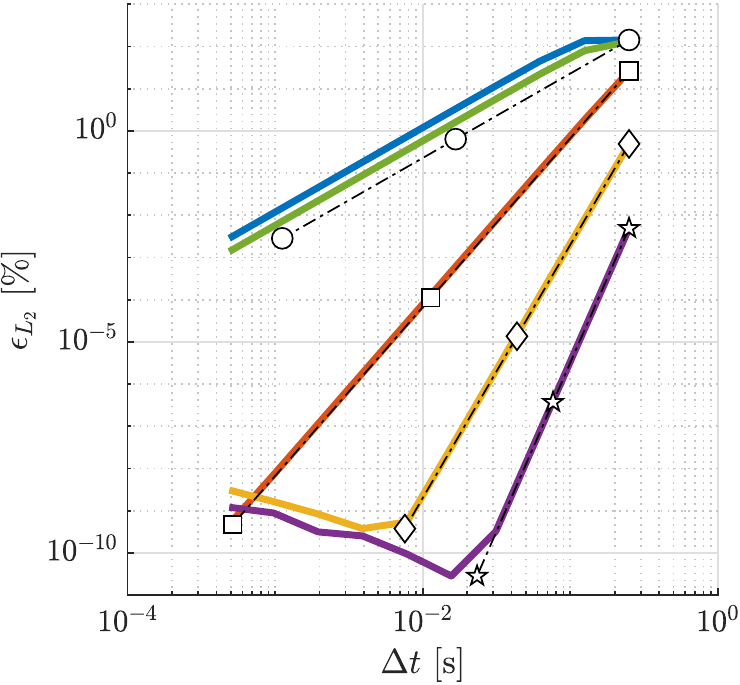}

}\hfill{}\subfloat[Case \#2\label{subfig:SDOF_2}]{\includegraphics{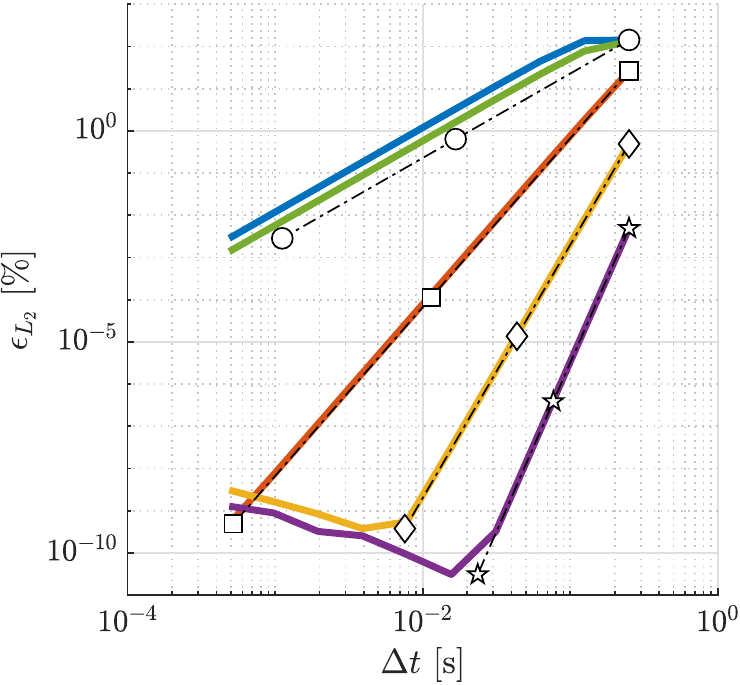}

}\\
 \subfloat[Case \#3\label{subfig:SDOF_3}]{\includegraphics{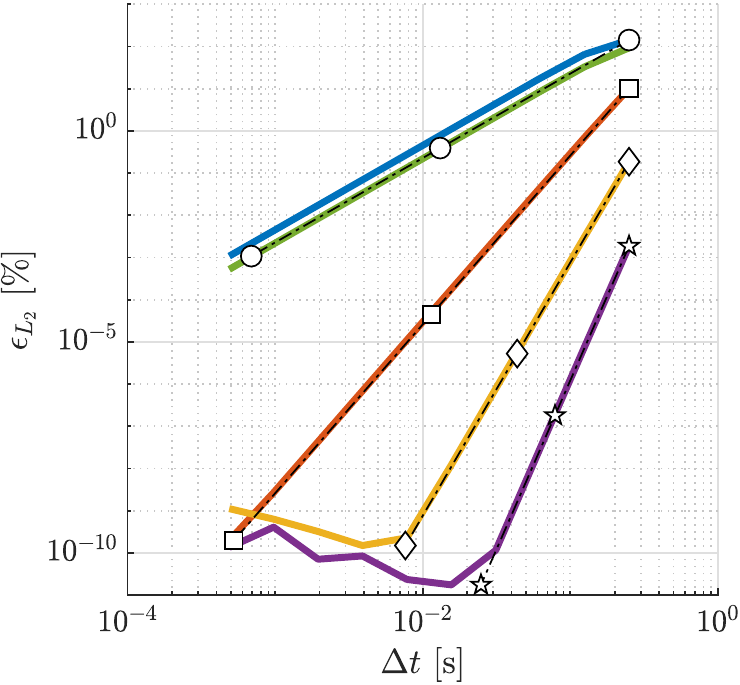}

}\hfill{}\subfloat[Case \#4\label{subfig:SDOF_4}]{\includegraphics{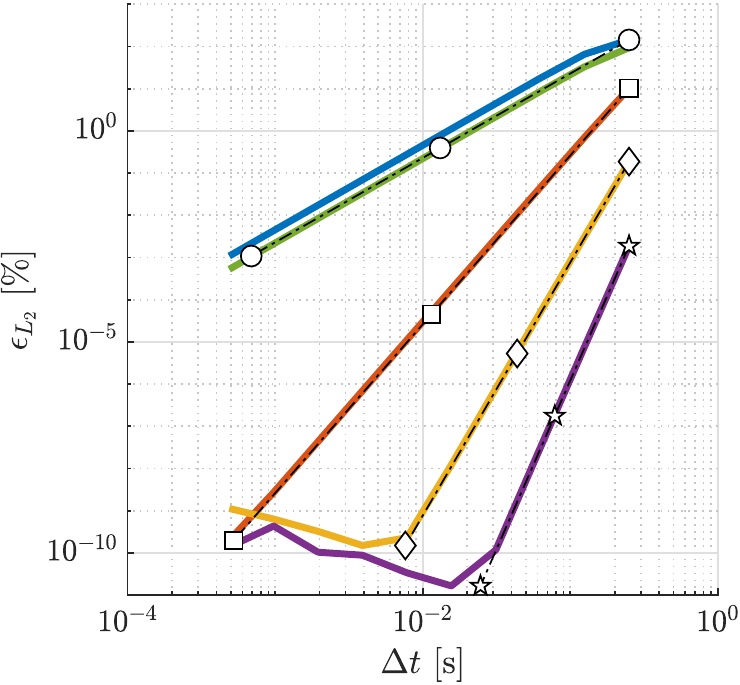}

}\\
 \subfloat[Case \#5\label{subfig:SDOF_5}]{\includegraphics{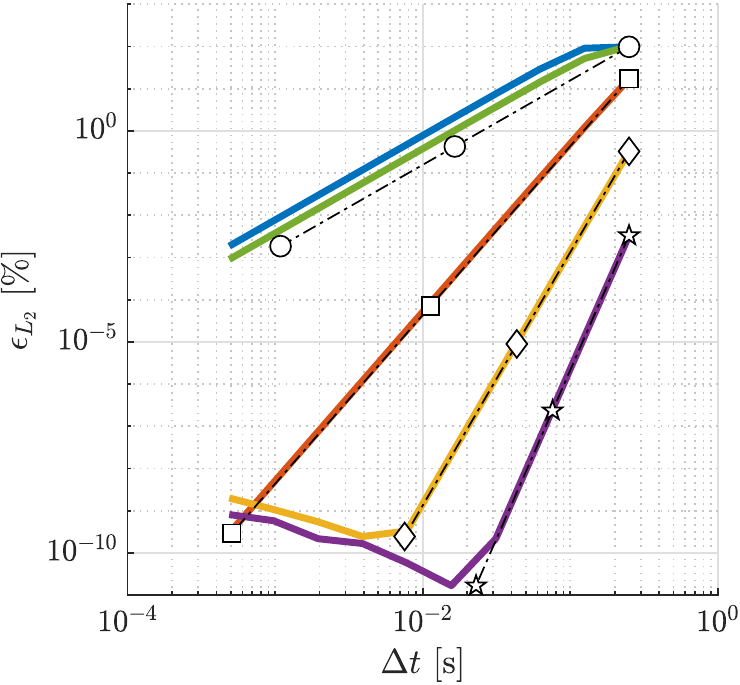}

}\hfill{}\subfloat[Case \#6\label{subfig:SDOF_6}]{\includegraphics{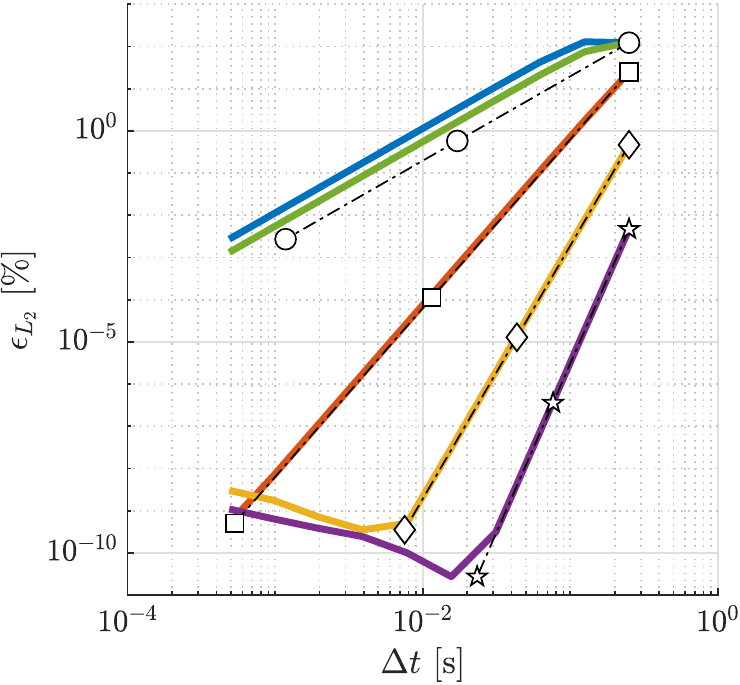}

}\caption{Displacement error in the $L_{2}$-norm for the different configurations
of the single degree of freedom system. The dash-dotted lines indicate
the optimal rates of convergence corresponding to slopes of 2 (circle),
4 (square), 6 (diamond), and 8 (pentagram), respectively. \label{fig:SDOF_Error}}
\end{figure}

The displacement error in the $L_{2}$-norm for the different set-ups
of the SDOF system is plotted in Fig.~\ref{fig:SDOF_Error}. The
thin dash-dotted lines are included to illustrate the theoretically
optimal convergence behavior. We clearly observe that the proposed
time-stepping method is converging with the optimal rate for all examples
until the error reaches a plateau at roughly $1\times10^{-9}\%$.
Independent of the ICs, the presence of physical damping, or the external
loading functions, high rates of convergence are achieved which highlights
the superior performance of the novel method. The error plateau is
clearly related to round-off errors that are inevitable in the numerical
implementation of the algorithm. 

As mentioned before, Bathe's method is included in the comparison
as a reference solution for a more recently developed algorithm. It
is an unconditionally stable composite time integration scheme consisting
of two sub-steps. In the first sub-step the constant average acceleration
method is employed, while in the second on a three-point backward
difference method is utilized. This combination yields a lower period
elongation compared to Newmark's method, while introducing a certain
degree of numerical damping. More details on the performance of this
method can be found in the pertinent literature \citep{ArticleBathe2005,ArticleBathe2012,ArticleNoh2018,ArticleMalakiyeh2019,ArticleNoh2019a,ArticleNoh2019b}.
Bathe's method has been intensively tested and is also implemented
in the commercial finite element software ADINA. Based on the characteristics
of Bathe's method it is expected that it exhibits similar convergence
properties compared to the present scheme of order $\mathcal{O}(1,1)$.
Due to the use of two sub-steps, which in trun means that the computational
effort is roughly doubled per time step, the overall error is slightly
lower for Bathe's method. Please keep in mind that for a fair comparison
in terms of accuracy and computational effort, the time step should
be doubled for Bathe's time-stepping scheme. 

\subsection{Single degree of freedom system: Period elongation and amplitude
error}

\label{sec:SDOF_PE_AE} The main goal of every time-stepping scheme
is to achieve a high-quality approximation of the actual dynamic response
of the structure that is being investigated. To this end, the time
increment must selected such that the maximum frequency $f_{\mathrm{max}}$
of interest is well resolved. As a rule of thumb, it is often stated
that ten increments per smallest period of interest $T_{\mathrm{min}}\,{=}\,\sfrac{1}{f_{\mathrm{max}}}$
are sufficient\footnote{Considering the Nyquist--Shannon sampling theorem, the absolute minimum
is two time steps per smallest period $\Delta t_{\mathrm{max}}\,{=}\,\sfrac{T_{\mathrm{min}}}{2}$.
For lower sampling rates the signal cannot be correctly reconstructed.
Please note that in high-order methods we have $p\,{+}\,1$ sampling
points per time step.} to achieve reasonably accurate results \citep{BookBathe2002} 
\begin{equation}
\Delta t=\cfrac{T_{\mathrm{min}}}{10}\,.\label{eq:DeltaT10}
\end{equation}
However, considering problems where highly precise numerical solutions
are required and the time integration method is only second-order
accurate (which applies to most established methods such as Newmark's
method \citep{ArticleNewmark1959}, the HHT-$\alpha$ method \citep{ArticleHilber1977},
the generalized-$\alpha$ method \citep{ArticleChung1993}, Bathe's
method \citep{ArticleBathe2005}, etc.), this estimate of the time
step size is generally not sufficient. From experience, approximately
100 sampling points per smallest period are recommended for high-fidelity
simulations \citep{ArticleTschoeke2018,ArticleDuczek2019b}.

The situation is obviously different when studying high-order accurate
time-stepping methods such as the one proposed in this article. Here,
larger time increments are possible due to the increased accuracy
of the algorithm. Hence, it is in principle possible to apply similar
refinement strategies in both space and time. The spatial \emph{h}-refinement
corresponds to a decrease in the time step size $\Delta t$ in the
time domain, while an increase in the polynomial order of the shape
functions is equivalent to an elevation of the order of the time integrator.
Thus, to achieve a highly accurate solution with the least numerical
costs a holistic \emph{hp}-refinement strategy might be developed
that is not only limited to the spatial discretization but also exploited
for the temporal one. This, however, is out of the scope of the current
contribution, where we want to introduce the algorithm and discuss
its fundamental properties. Therefore, three important quantities,
i.e., stability, amplitude error (AE), and period elongation (PE),
will be briefly discussed in the remainder of this section. To this
end, we consider the SDOF system described by case \#1 for which the
analytical solution is known: $u(t)\,{=}\,\cos(\omega_{\mathrm{n}}t)$.
The simulation is run for 10,000 periods of vibration to study the
effect of long-term numerical analysis. For a complete analysis, different
IVPs with different ICs, damping, and loading functions would have
to be considered. However, the defining numerical properties of the
time integration scheme can be investigated by only solving the problem
stated above~\citep{BookBathe2002}. For a general comparison, the
results obtained using Newmark's method (constant average acceleration),
the HHT-$\alpha$ method ($\alpha\,{=}\,-0.05;$default setting in
ABAQUS/standard), the generalized-$\alpha$ method ($\rho_{\infty}\,{=}\,0.8$),
and Bathe's method ($\gamma\,{=}\,(2-\sqrt{2})\Delta t$; splitting
ratio) are included in the discussion. 

\subsubsection{Stability}

\label{sec:Stability} In many problems of practical interest, the
response of a structure is dominated by a limited range of frequencies,
while the contribution of higher modes is essentially negligible.
Hence, it is not meaningful to assign a time step size $\Delta t$
depending on the smallest period of the numerical system. Instead,
we are only interested in the maximum frequency with a significant
contribution to the structural response. This means, the higher modes
are integrated with a time increment that not sufficient to fully
resolve them. Hence, the question arises how the time-stepping scheme
can handle large values of the ratio $\sfrac{\Delta t}{\tilde{T}_{\mathrm{min}}}$,
where $\tilde{T}_{\mathrm{min}}$ denotes the smallest (numerical)
period of the system under investigation. This leads us to the topic
of stability of a time-stepping scheme.

According to Ref.~\citep{Gallopoulos1989}, Padé-based time integrators
are unconditionally stable since the spectral radius of the rational
approximation---see Eq.~\eqref{eq:ExpmPade}---of the matrix exponential---see
Eq.~\eqref{eq:timeSteppingExpm}---is of unit value, which means
that amplitude decay is not observed in these methods. Unconditional
stability is a very useful property to have as the error in the simulation
of a dynamic problem will not diverge independent of the choice of
the time step size $\Delta t$. Note that if the value for $\Delta t$
is chosen incorrectly, i.e., too large, the solution can still be
arbitrarily wrong. 

\subsubsection{Amplitude error}

\label{sec:AmplitudeError} Please note that a rigorous study of the
amplification matrix for the proposed approach is out of the scope
of the current contribution and will be included in future communications
on a high-order implicit scheme with controllable numerical damping.
Nonetheless, a simplified error measure referred to as amplitude error
(AE) is introduced at this point to conclusively illustrate the performance
of the novel time integrator. To this end, the amplitude of the displacement
response is evaluated at multiplies of the actual period of vibration
$T_{\mathrm{n}}\,{=}\,\sfrac{2\pi}{\omega_{\mathrm{n}}}$ and a relative
error is computed over all periods 
\begin{equation}
\mathrm{AE}=\cfrac{1}{N}\sum\limits _{k=1}^{N}\cfrac{||u_{\mathrm{ref}}(kT_{\mathrm{n}})-u_{\mathrm{num}}(kT_{\mathrm{n}})||}{u_{\mathrm{ref}}(kT_{\mathrm{n}})}\times100[\%]\,.\label{eq:AE}
\end{equation}

\begin{figure}[t!]
\centering \subfloat{\includegraphics{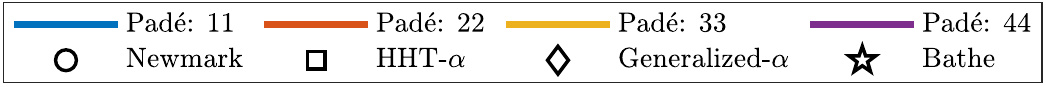}}\\
 \setcounter{subfigure}{0} \subfloat[Relative amplitude error in semi-log scale\label{subfig:SDOF_AElin}]{\includegraphics{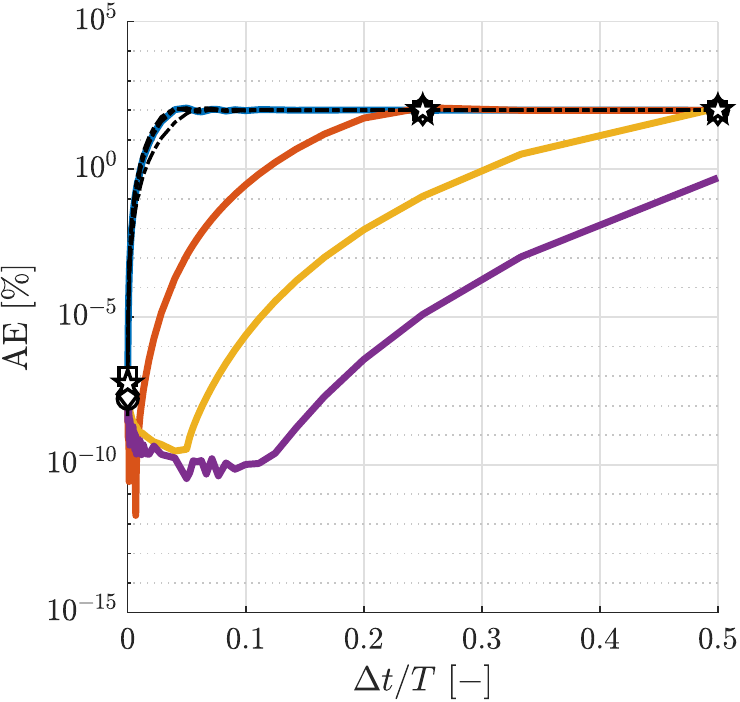}

}\hfill{}\subfloat[Relative amplitude error in log-log scale\label{subfig:SDOF_AE}]{\includegraphics{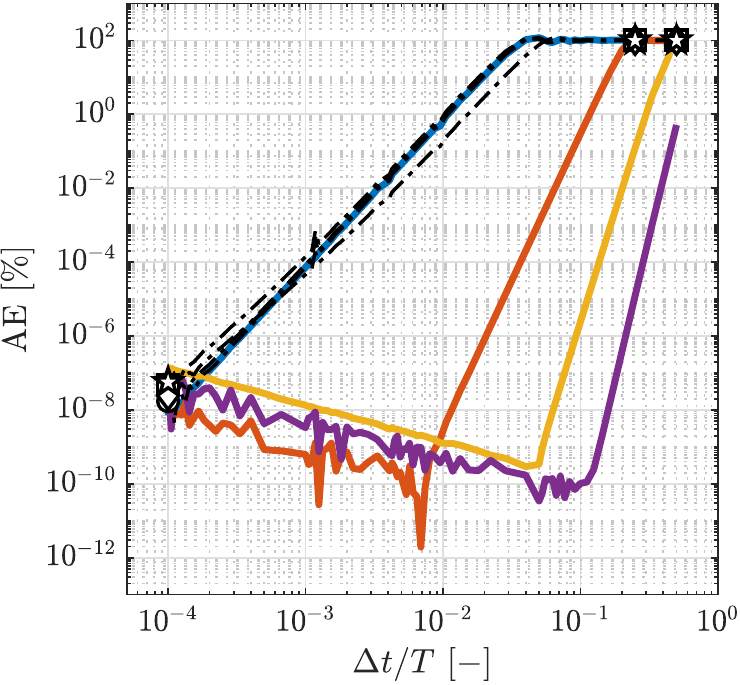}

}

\subfloat[Period elongation\label{subfig:SDOF_PElin}]{\includegraphics{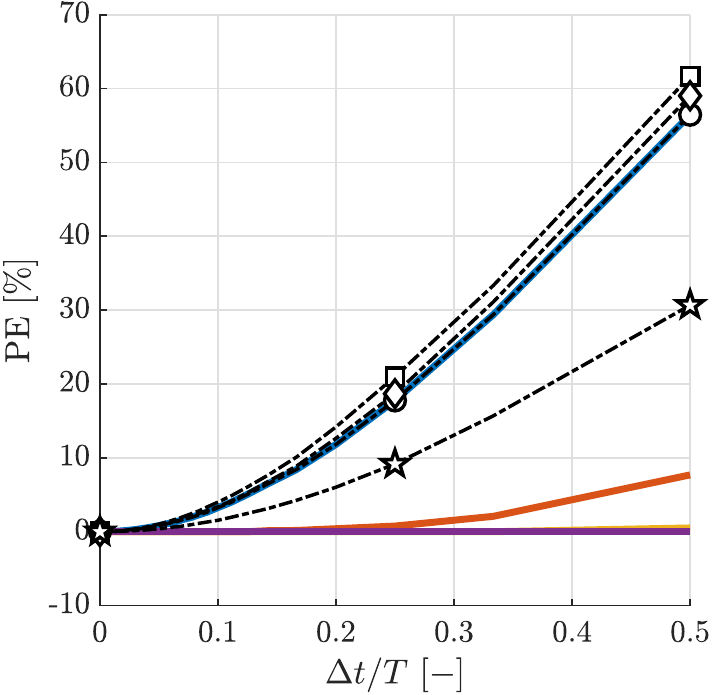}

} \hfill{}\subfloat[Period elongation in log-log scale\label{subfig:SDOF_PE}]{\includegraphics{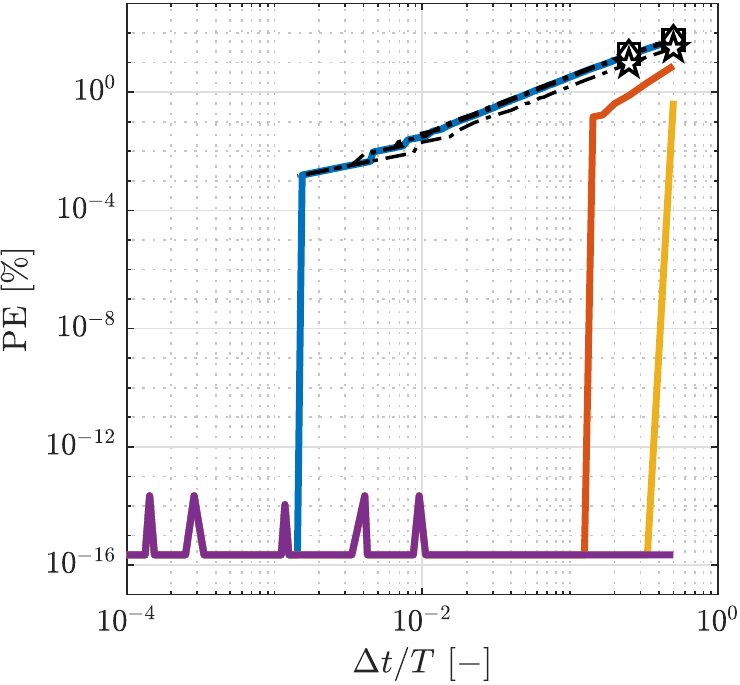}

}\caption{Amplitude error and period elongation after 10,000 periods of vibration---Comparison
of the novel scheme with established methods. \label{fig:SDOF_PE_AE}}
\end{figure}

The results of this particular analysis are depicted in Figs.~\ref{subfig:SDOF_AElin}
and \ref{subfig:SDOF_AE} in linear and logarithmic scale, respectively.
The numerical results highlight that the novel scheme shows a very
low amplitude error over time even if long-term simulations are investigated.
It can be observed that depending on the chosen order already rather
large time steps yield exceptionally accurate results in terms of
the recovery of the amplitude of the displacement response. To put
this into perspective, recall that in engineering analysis an error
of roughly 1\% is often deemed acceptable. Taking this measure, we
see that the eighth-order accurate scheme provides accurate results
already for a time step $\Delta t\,{=}\,\sfrac{T_{\mathrm{n}}}{2}$,
which is the absolute minimum increment that is required to theoretically
resolve the vibration for low-order schemes\footnote{In high-order schemes even larger time steps can be used as implicitly
more sampling points are included in the analysis. This can be seen
for example in the interpolation of the force vector per time step.}. Naturally, smaller time steps are required considering lower approximation
orders. Considering the sixth-order accurate algorithm a time step
of $\Delta t\,{<}\,\sfrac{T_{\mathrm{n}}}{3}$ is needed, while in
the fourth-order case the time increment needs to be reduced further
to $\Delta t\,{<}\,\sfrac{T_{\mathrm{n}}}{8}$. Recall that established
time integration methods are often only second-order accurate and
therefore, a significantly smaller time step is required for those
techniques. From Fig.~\ref{subfig:SDOF_AE}, we infer that $\Delta t\,{<}\,\sfrac{T_{\mathrm{n}}}{85}$
is needed to achieve the prescribed error threshold of 1\%. Thus,
by increasing the order of accuracy of the time marching algorithm,
the required time step can be significantly reduced. This behavior
is expected and corresponds to what we observe when increasing the
polynomial order of the spatial discretization. It has to be kept
in mind, however, that by increasing the order of accuracy, we also
increase the numerical costs of the method which will be discussed
in more detail using more complex examples featuring a moderate number
of DOFs. 

Since the amplitude error uses the displacement response evaluated
at multiples of the natural period of the system under investigation
it essentially combines two types of errors: (i) the actual error
in the maximum displacement and (ii) the period elongation (PE) of
the numerical algorithm. Consequently, we will have a closer look
at the PE error in the following subsection. 

\subsubsection{Period elongation}

\label{sec:PeriodElongation} The error cause by period elongation
leads to the computational phenomenon referred to as numerical dispersion,
i.e., the numerical solution seems to approximate a system with a
different natural frequency. This shift can be both positive (as is
the case for Newmark's constant average acceleration method) or negative
(as is the case in the central difference method). That is to say,
in the numerical results either a shift to a higher or lower frequency
will be apparent.

The error in the approximation of the natural period $T_{\mathrm{n}}$
is assessed by determining the peak values of the displacement response
and saving the corresponding time values in a vector. Computing the
differences between subsequent components of the time vector gives
us the values for the natural period as approximated by the numerical
solution. 
\begin{equation}
\mathrm{PE}=\cfrac{1}{N}\sum\limits _{k=1}^{N}\cfrac{t\left(u_{\mathrm{max,num}}^{(k)}\right)-t\left(u_{\mathrm{max,ref}}^{(k)}\right)}{t\left(u_{\mathrm{max,ref}}^{(k)}\right)}\times100[\%]\,,
\end{equation}
where $u_{\mathrm{max,num}}^{(k)}$ and $u_{\mathrm{max,ref}}^{(k)}$
are the maximum amplitude values for the k-th vibration. The value
that is calculated here is an average period elongation over 10,000
vibration periods as mentioned before.

The numerical results of the PE analysis are depicted in Figs.~\ref{subfig:SDOF_PElin}
and \ref{subfig:SDOF_PE} in linear and logarithmic scale, respectively.
Due to the accuracy of the proposed time-stepping method, even after
10,000 periods of vibration there is basically no phase shift detectable
and therefore, the error is set to the value of machine precision
\texttt{eps} to evaluate the logarithm. In the initial stage, we observe
a monotonous decrease of the PE-error before a sudden jump occurs,
telling us that for the chosen simulation time no phase shift is present
in the numerical results. Overall, the curves highlight that the proposed
scheme exhibits very small values of PE or in other words, a negligible
numerical dispersion is introduced. Thus, the present method is very
attractive for long-term simulations with high accuracy requirements.
We observe that the error in the natural period is basically negligible
over the entire range of time increments if a method of order six
or eight is chosen. Even at order four the shift in frequency is significantly
lower compared to all established second-order accurate methods. For
the maximum time step $\Delta t\,{=}\,\sfrac{T_{\mathrm{n}}}{2}$,
the proposed algorithm exhibits PE-values of approximately 0\%, 0.5\%,
7.7\%, and 56.5\%. In terms of the PE-error, the worst performance
is exhibited by the HHT-$\alpha$ method with 61.7\%, while Bathe's
method shows the best performance for all second-order approaches
with 30.7\%.

In summary, the numerical results demonstrate that any time-stepping
method can be used and is very accurate as long as the ratio of time
step to natural period (or period of interest) $\sfrac{\Delta t}{T_{\mathrm{n}}}$
is below 0.01. In cases where this ratio is larger, the investigated
methods show distinct difference that need to be taken into account
when solving dynamical problems.

\subsection{Dynamic behavior of a rectangular domain}

\label{sec:2D_Rod} 
\begin{figure}[!b]
\centering \includegraphics{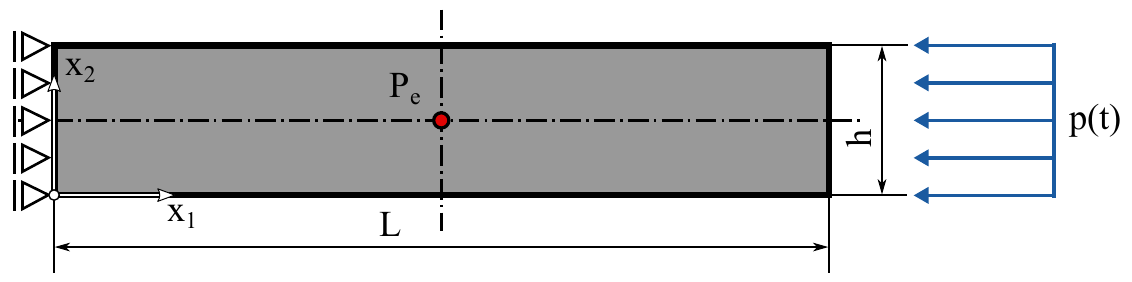} \caption{Homogeneous rectangular rod under plane stress conditions. \label{fig:ModelRod}}
\end{figure}

In this section, the dynamic behavior of a two-dimensional
rectangular body in plane stress conditions is considered \citep{ArticleSoares2015}.
In principle, the structure acts as a one-dimensional rod with the
following material properties: Young's modulus $E\,{=}\,100\,$Pa,
Poisson's ratio $\nu\,{=}\,0.0$, and mass density $\rho\,{=}\,1\,$kg/m$^{3}$.
The dimensions of the rectangular domain are $L\,{=}\,1\,$m and $h\,{=}\,0.2\,$m
(see Fig.~\ref{fig:ModelRod}). A pressure load is applied at the
right edge of the rod, while displacements in horizontal direction
are fixed at its left edge. In contrast to the examples reported in
the body of literature \citep{ArticleSoares2020b,ArticleNoh2019a,ArticleKim2020b},
the loading is not applied as a Heaviside function, even though an
analytical solution is readily available due to the (mathematical)
simplicity of the excitation. Instead, a time-dependent amplitude
function $p(t)$ in form of a sine-burst is used 

\begin{equation}
p(t)=P_{0}\sin(2\pi f_{\mathrm{ex}}t)\exp^{-\cfrac{1}{2}\left(\cfrac{t-t_{0}}{\tau}\right)^{2}}\,,\label{eq:BurstExcitation}
\end{equation}
where $f_{\mathrm{ex}}$ is the center frequency of the excitation
signal, $P_{0}$ is the amplitude, and for the sake of simplicity,
the variables $\tau$ and $t_{0}$ are chosen as functions of the
period of excitation $T\,{=}\,\sfrac{1}{f_{\mathrm{ex}}}$ as 
\[
\tau=T\qquad\text{and}\qquad t_{0}=4T\,.
\]
Hence, the excitation signal is defined by a single parameter which
is $f_{\mathrm{ex}}$. Note that the chosen excitation function makes
the derivation of a analytical solution based on Duhamel's integral
very difficult and therefore, the reference results for this example
are based on a numerical overkill solution of the system. An analytical
solution for the displacement response due to a trigonometric excitation
functions can be found in Ref.~\citep{ArticleGravenkamp2020}, while
the solution for a Heaviside excitation are provided in Ref.~\citep{ArticleSoares2015}.

\begin{figure}[!t]
\centering \subfloat[Time-domain signal $f(t)$]{\includegraphics{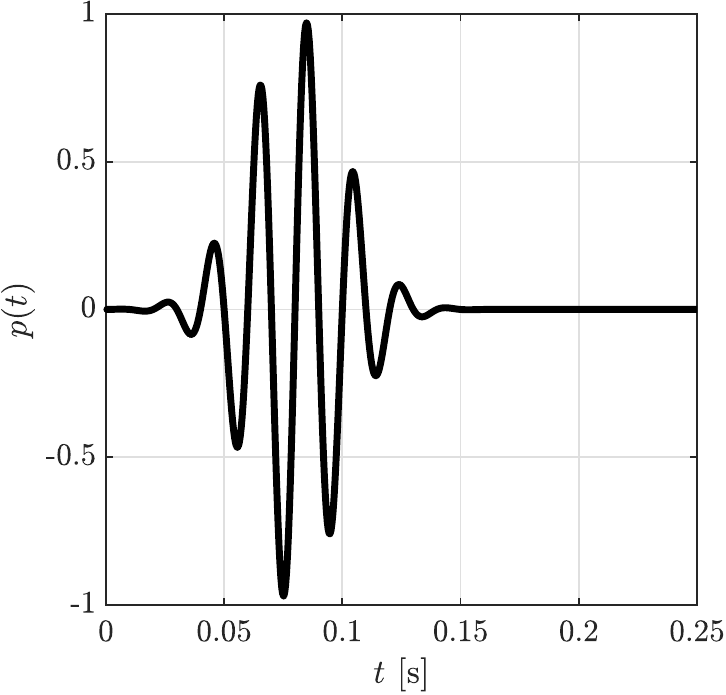}

}\hfill{}\subfloat[Frequency-domain signal $\hat{F}(f)$]{\includegraphics{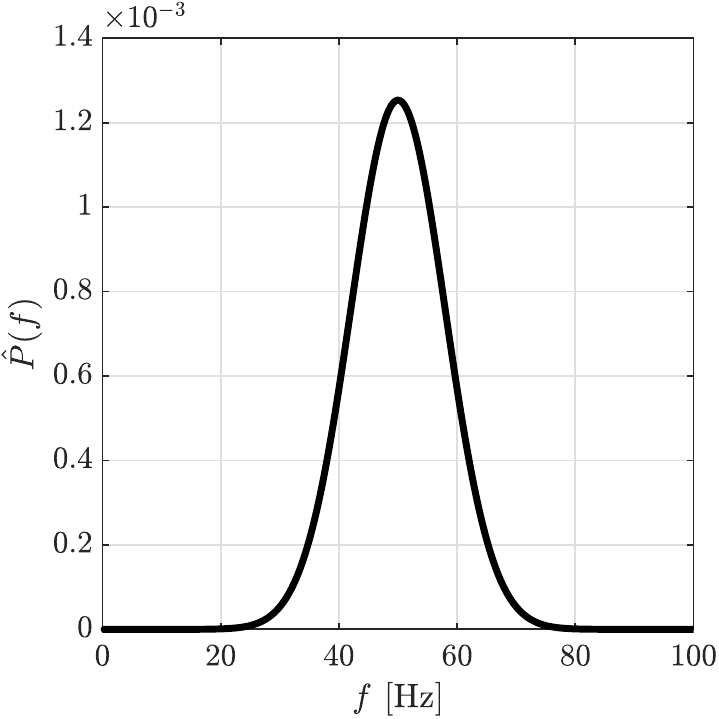}

}\caption{Excitation signal in the time and frequency domain---Rod. \label{fig:ExcitationSignalRod}}
\end{figure}

In Fig.~\ref{fig:ExcitationSignalRod}, both the time history and
the frequency spectrum of the loading function are depicted. For the
current analysis, the center frequency $f_{\mathrm{ex}}$ is set to
$50\,$Hz, while the amplitude of the loading function $P_{0}$ is
taken as $1\,$N. The maximum frequency of interest $f_{\mathrm{max}}$
is determined at the threshold when the amplitude in the frequency-spectrum
is constantly below 1\% of its maximum value. In our specific example,
$f_{\mathrm{max}}$ is 74.125$\,$Hz. To obtain a reasonable temporal
resolution it is recommended (at least for commonly used second-order
accurate time-stepping schemes) to employ at least 20 time steps per
period of vibration determined by the maximum frequency of interest.
Thus, an upper limit for the time step size would typically be $\Delta t_{\mathrm{max}}\,{=}\,\sfrac{1}{(20f_{\mathrm{max}})}\,{=}\,6.75\times10^{-4}\,$s.
However, due to the use of a high-order time-stepping scheme the maximum
time increment can be significantly increased to $\Delta t_{\mathrm{max}}\,{=}\,1.0\times10^{-2}\,$s.

In the current contribution, the main goal is to study the performance
(accuracy and efficiency) of the proposed time integration scheme.
At this point, we are not interested in including effects due to the
spatial discretization in our analysis. Therefore, we limit our discussions
on meshes consisting of bi-linear (4-node) finite elements. These,
elements are implemented in virtually every commercial finite element
code and are widely used in structural dynamics. We are well aware
of the fact that, in general, high-order finite elements offer advantages
in terms of higher convergence rates \citep{ArticleWillberg2012}
and if they are based on Lagrangian shape functions and non-equidistant
nodal distributions (e.g., the spectral element method (SEM) \citep{BookKarniadakis2005,BookPozrikidis2014})
the mass matrix can easily be diagonalized \citep{ArticleDuczek2019b}\footnote{Remark: Only in the spectral element method a variationally consistent
mass lumping procedure which retains the theoretical optimal rates
of convergence can be devised.}. However, they have only found limited application outside academia
and therefore, our analysis is at least for now restricted to elements
with linear shape functions.

It is well-known that a high-order accurate time integration scheme
is numerically more expensive on a per time step basis compared to
low-order schemes and therefore, the higher costs need to be amortized
by more accurate solutions and consequently, the possibility to use
significantly larger time steps. Overall, a solution with a prescribed
error threshold needs to be achieved with the least possible computational
effort. In this section, all analyses are compared with the constant
average acceleration method of the Newmark family \citep{ArticleNewmark1959}.
This established time-stepping scheme serves as a benchmark for assessing
the quality of the numerical solution in terms of accuracy and the
computational costs.

Considering the spatial discretization eight different meshes consisting
of square-shaped (bi-linear) finite elements have been set up. In
the coarsest mesh, the element size has been set to $0.1\,$m, i.e.,
20 elements (10$\times$2) are created. For each new discretization,
the element size is halved and thus, the number of elements is increased
by a factor of four. For the finest mesh, we obtain 327,680 elements
(1,280$\times$256). That is to say, in terms of the number of degrees
of freedom ($n_{\mathrm{DOF}}$), we cover a range from 66 to 658,434
DOFs, giving us a good insight into the performance of the present
time integrator for small to medium-sized systems. Note that in the
current implementation a direct solver is used in each time step,
i.e., in linear dynamics the decomposition of the system matrix is
pre-computed such that only forward and backward substitution steps
need to be performed. For even larger systems, it is worth exchanging
the direct solver with an iterative solver to facilitate a parallel
implementation of the time-integrator on high-performance clusters
\citep{ArticleZhang2021}.

Because the Poisson's ratio is set to zero and therefore, the structure
essentially behaves like a one-dimensional rod, the longitudinal (pressure)
wave velocity, denoted as $c_{\mathrm{L}}$ is simply defined as 
\begin{equation}
c_{\mathrm{L}}=\sqrt{\cfrac{E}{\rho}}\,.\label{eq:cL_rod}
\end{equation}
We have to keep in mind, however, that this expression for the wave
velocity only holds for a one-dimensional structure. According to
Eq.~\eqref{eq:cL_rod}, the wave velocity is $10\,\sfrac{\mathrm{m}}{\mathrm{s}}$
in our example. The corresponding wavelength $\lambda$ is related
to the frequency and the wave velocity by 
\begin{equation}
\lambda=\cfrac{c_{\mathrm{L}}}{f_{\mathrm{ex}}}\,.\label{eq:Wavelength}
\end{equation}
For the given values, the wavelength is $\lambda\,{=}\,0.2\,$m. Consequently,
the spatial resolution of the wave packets varies between 2 and 256
elements per wavelength. As a rule of thumb, 10 elements with linear
shape functions are always recommended although it has been shown
in Ref.~\citep{ArticleWillberg2012} that this is an absolute minimum,
and the error is still well above the 1\% threshold.

For each simulation, the so-called Courant-Friedrichs-Lewy (CFL) number
can be determined, which provides a necessary condition for the convergence
of a numerical solution to partial differential equations (PDEs) \citep{ArticleCourant1928}.
Although, it is closely related to the numerical analysis of explicit
time integration schemes, we would like to mention it as a useful
measure to assess the interplay between spatial and temporal discretizations.
It is commonly defined as 
\begin{equation}
\mathrm{CFL}=\cfrac{c\,\Delta t}{\Delta x}\label{eq:CFL}
\end{equation}
where $c$ is the characteristic wave velocity of the medium and $\Delta x$
denotes the element size (in a structured grid). Hence, the value
of the CFL-number tells us through how many elements a wave can travel
per each time step. This definition is, however, related to finite
elements featuring linear shape functions. In general, a CFL-number
below the value of one means that a wave travels less than the distance
between two neighboring nodes. 

\subsubsection{Analysis of the computational costs}

\label{sec:CostsRod} 
\begin{figure}[!b]
\centering \includegraphics{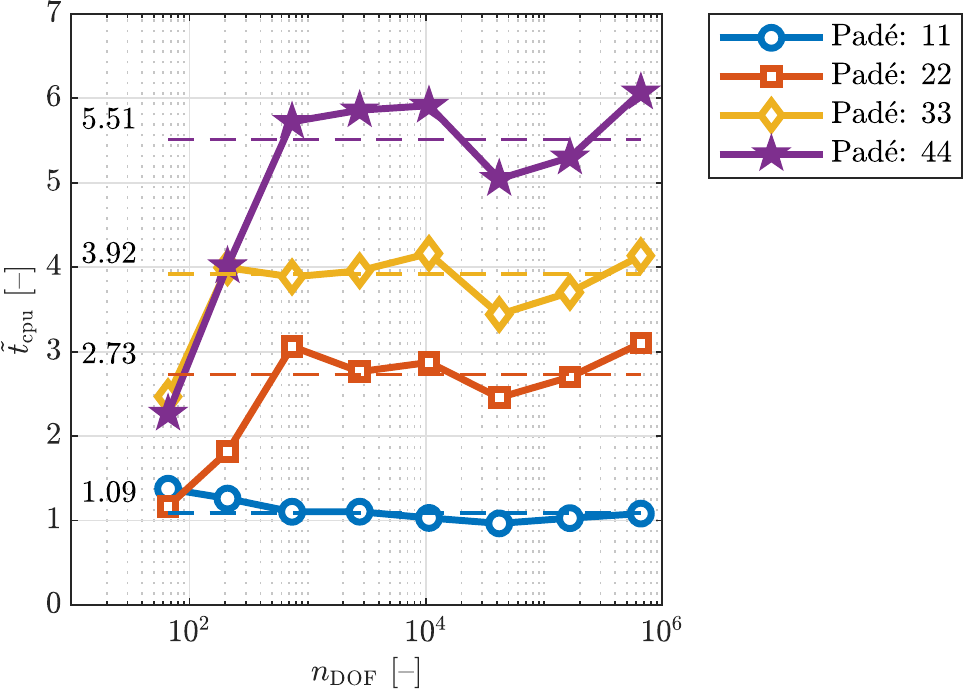} \caption{Computational time normalized with respect to Newmark's constant average
acceleration method---Rod. \label{fig:tnormRod}}
\end{figure}

As mentioned before, it is only natural that
high-order accurate time integration schemes are invariably more costly
per time step compared to conventionally used second-order schemes.
Therefore, it is worthwhile to investigate the required computational
time. In the following, we will normalize the computational times
obtained with the present high-order scheme with respect to the values
obtained employing Newmark's constant average acceleration method.
The results of this analysis are depicted in Fig.~\ref{fig:tnormRod}.
The simulations are run for the eight different discretizations with
a time step size of $\Delta t\,{=}\,1.5625\times10^{-4}\,$s, which
results in 6,400 time steps for a simulation time of $t_{\mathrm{sim}}\,{=}\,1\,$s.

As derived in Sect.~\ref{subsec:Order11}, the present scheme of
order $\mathcal{O}(1,1)$ is mathematically identical to Newmark's
constant average acceleration method and therefore, the computational
times of both approaches should be (ideally) identical. In Fig.~\ref{fig:tnormRod},
we observe that a median value of 1.09 is reached for the normalized
computational time $\tilde{t}_{\mathrm{cpu}}$ over a wide range of
number of DOFs. This slight difference from the theoretically expected
value of 1.0 is related to the different implementations regarding
the solution procedure introduced in Sect.~\ref{sec:HighOrderPade},
compared to standard implementations of Newmark's method as detailed
in Ref.~\citep{BookBathe2002}. Overall, this is of no concern and
shows again the good agreement. Considering the present time integration
method of orders $\mathcal{O}(2,2)$, $\mathcal{O}(3,3)$, and $\mathcal{O}(4,4)$
which are fourth-, sixth-, and eighths-order accurate, respectively,
an increase in the computational costs is noted. This increase is,
however, very moderate, and median values of 2.73, 3.92, and 5.51
are achieved for the different schemes. That is to say, for this particular
example the eighth-order accurate scheme is less than six times more
costly than a standard second-order accurate implicit method of the
Newmark family. This value is very promising considering the gained
accuracy. 

\subsubsection{Analysis of the accuracy of the time-integrator}

\label{sec:AccuracyRod} Besides the computational costs, the attainable
accuracy is another important aspect when investigation time-stepping
methods. Here, the convergence with respect to the exact solution
of the discretized problem is studied. We are only interested in assessing
the error introduced by the temporal discretization and neglect additional
sources of errors such as the spatial discretization and incorrect
mathematical assumptions on the material behavior.

\begin{figure}[!b]
\centering \includegraphics{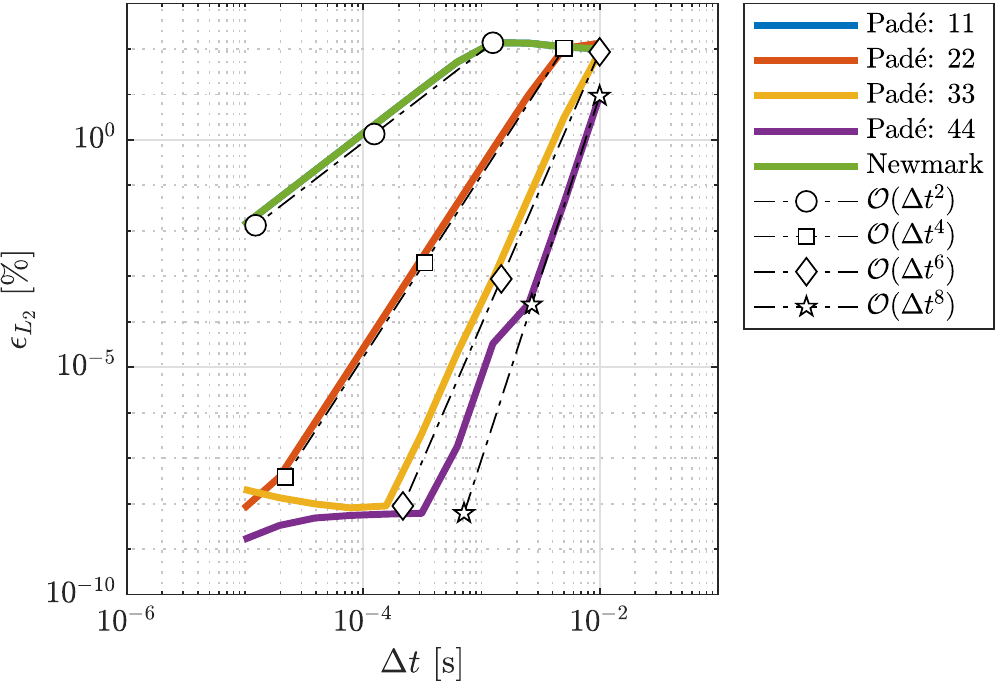} \caption{Displacement error in the $L_{2}$-norm for the homogeneous rectangular
rod. The dash-dotted lines indicate the optimal rates of convergence
corresponding to slopes of 2 (circle), 4 (square), 6 (diamond), and
8 (pentagram), respectively. \label{fig:ErrorRod}}
\end{figure}

In order to assess the accuracy the displacement error in the $L_{2}$-norm
($\epsilon_{L_{2}}$)---see Eq.~\eqref{eq:ErrorL2}---is computed.
In Fig.~\ref{fig:ErrorRod}, the error is depicted for a spatial
discretization using 1,280 (bi-linear) finite elements (80$\times$16).
The convergence curves are virtually identical for the other discretizations
and therefore, only this example is plotted. Similar to the simple
SDOF-systems, we observe optimal convergence until an error plateau
is reached. As theoretically predicted, the present integrator of
order $\mathcal{O}(1,1)$ and Newmark's constant average acceleration
method yield identical results. Due to the higher rates of convergence
a similar accuracy can be easily reached using much larger time step
sizes when the novel approach is employed. Considering an error threshold
of roughly 1\%, which is acceptable in most engineering applications,
a time increment of $\Delta t\,{=}\,8.4\times10^{-5}\,$s is required
for Newmark's constant average acceleration method, while for the
high-order schemes significantly larger values are acceptable. Considering
the present scheme of order $\mathcal{O}(2,2)$, $\Delta t\,{=}\,1.4\times10^{-3}\,$s
is sufficient, while the time increment can be further increased for
orders $\mathcal{O}(3,3)$ and $\mathcal{O}(4,4)$ to $\Delta t\,{=}\,4.1\times10^{-3}\,$s
and $\Delta t\,{=}\,7.6\times10^{-3}\,$s, respectively. That is to
say, the time step size can be increased by factors of 17, 49, and
90 for the novel time-integrator. The results discussed in the current
subsection need to be combined with the findings of the previous subsection
to arrive at a first conclusion regarding the viability of our high-order
time-stepping method. To this end, we recall that the increased computational
costs per time step are depicted in Fig.~\ref{fig:tnormRod}. Here,
it has been observed that the costs increase by a factor of less than
3, 4, and 6 for the fourth-, sixth- and eighth-order accurate schemes.
This is much lower compared the amount of time steps that can be saved.
For this example, a speed-up of roughly 6, 12, and 15 can be achieved.
However, note that this value is problem dependent and might also
be related to the properties of the numerical method that is employed
for the spatial discretization.

\begin{figure}[!t]
\centering \includegraphics{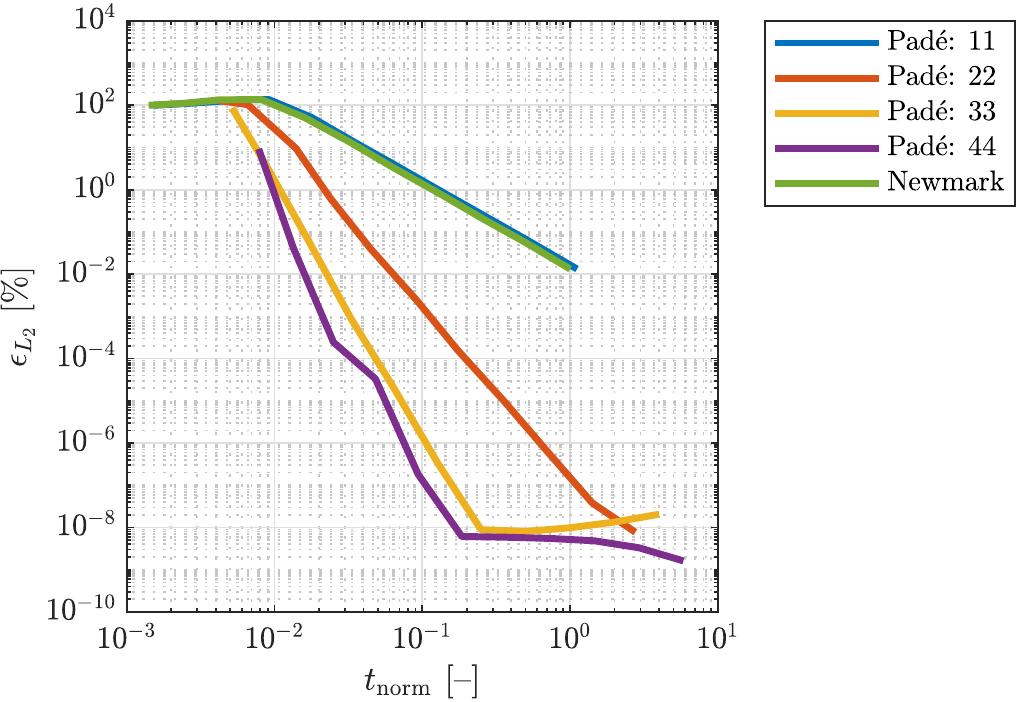} \caption{Displacement error in the $L_{2}$-norm for the different configurations
over the normalized computational time---Rod. \label{fig:ErrorTcpuRod}}
\end{figure}

The findings obtained in this section are summarized in Fig.~\ref{fig:ErrorTcpuRod},
where we plot the attainable accuracy over the normalized computational
time for a spatial discretization of 1,280 (bi-linear) finite elements
(80$\times$16). Here, the simulation time is divided by the maximum
value for Newmark's constant average acceleration method, i.e., we
take the computational time corresponding to Newmark's method with
a time step size of $\Delta t_{\mathrm{max}}\,{=}\,9.765625\times10^{-6}\,$s
as a reference. This figure is in principle a combination of Figs.~\ref{fig:tnormRod}
and \ref{fig:ErrorRod}, illustrating again the advantages of employing
a high-order time integration scheme. We clearly observe that for
a given simulation time the accuracy significantly increases when
using the proposed method (add a vertical line to the graph with the
targeted simulation time), while a prescribed error threshold can
be reached in a shorter time (add a horizontal line to the graph with
the targeted error value). Note that achieving this kind of efficiency
is not trivial for high-order schemes, which are often too costly
for moderately large systems. The results reported in this section
for the relatively coarse mesh are also valid for the finer meshes
that have been investigated for this example. Since the results are
virtually identical it is not meaningful to provide all results and
therefore, only selected graphs are included in the article.

\begin{figure}[!b]
\centering \includegraphics{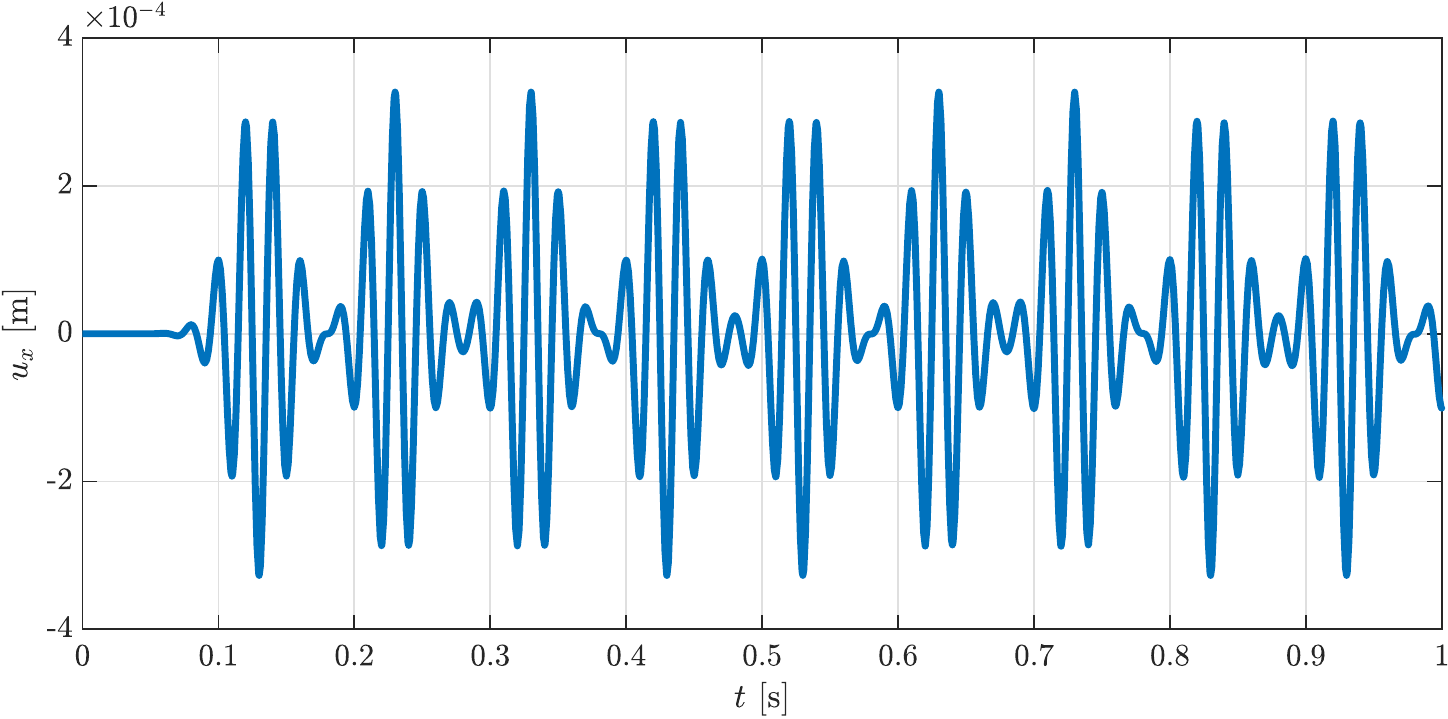} \caption{Displacement in $x$-direction at the observation point $P_{\mathrm{e}}$---Rod.
\label{fig:Rod_UX}}
\end{figure}

For the sake of completeness, we include the displacement response
in $x$-direction at point $P_{\mathrm{e}}$ located at (0.5$\,$m,
0.1$\,$m), see Fig.~\ref{fig:Rod_UX}. The multiple reflections
at the left and right edges are clearly seen, showing how the waves
travels back and forth in the rod-like structure. 

\subsection{Lamb's problems}

\label{sec:Lamb} 
\begin{figure}[!b]
\centering \includegraphics{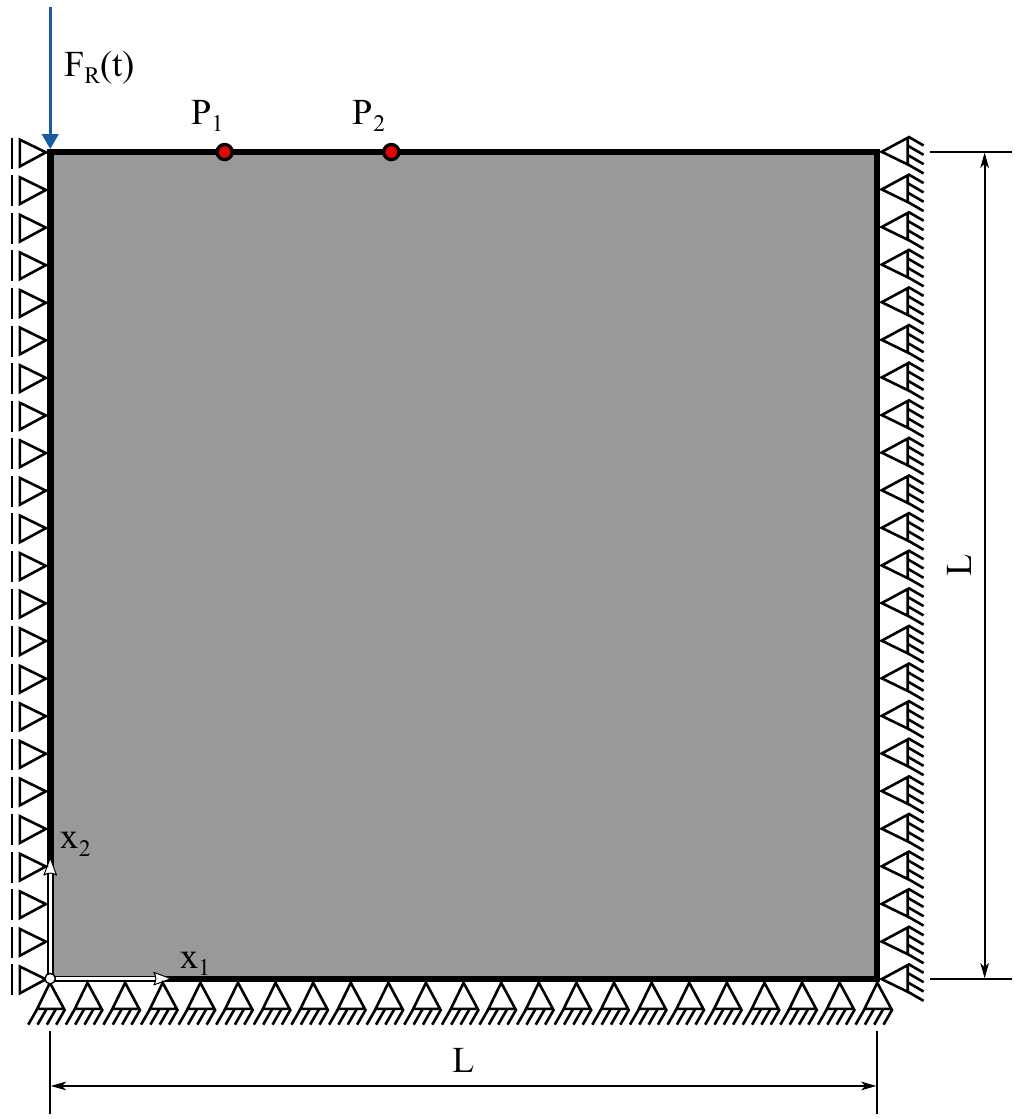} \caption{\textquotedblleft Semi-infinite\textquotedblright{} elastic domain
under plane strain conditions. \label{fig:ModelLamb}}
\end{figure}

In this section, we consider Lamb's problem, in
which wave propagation in a semi-infinite elastic domain are considered
under plane strain conditions. The problem set-up and the material
properties, reported in the following, are taken from Refs.~\citep{ArticleHam2012,ArticleNoh2013}.
The geometry and loading conditions are depicted in Fig.~\ref{fig:ModelLamb}.
The material properties are chosen such that the P-wave velocity ($c_{\mathrm{L}}$)
is 3,200$\,$m/s, the S-wave velocity ($c_{\mathrm{T}}$) is 1,847.5$\,$m/s,
and the Rayleigh wave velocity ($c_{\mathrm{R}}$) is 1,671$\,$m/s,
while the mass density ($\rho$) is 2,200$\,$kg/m$^{3}$. The dimension
of the computational domain $L$ is set to 3,200$\,$m which ensures
that no unwanted reflections from the boundaries of the structure
are present in the displacement response ($t_{\mathrm{sim}}\,{=}\,1\,$s).

Based on the wave velocities and the mass density, Lamé's constants
$\tilde{\lambda}$ and $\tilde{\mu}$ (shear modulus) can be determined
\begin{alignat}{3}
c_{\mathrm{L}} & =\sqrt{\cfrac{\tilde{\lambda}+\tilde{\mu}}{\rho}}\qquad &  & \longrightarrow\qquad\tilde{\lambda} &  & =c_{\mathrm{L}}^{2}\rho-\tilde{\mu}\,,\\
c_{\mathrm{T}} & =\sqrt{\cfrac{\tilde{\mu}}{\rho}}\qquad &  & \longrightarrow\qquad\tilde{\mu} &  & =c_{\mathrm{T}}^{2}\rho\,.
\end{alignat}
To obtain the typical engineering constants Young's modulus $E$ and
Poisson's ratio $\nu$ the following conversion is required 
\begin{align}
E & =\cfrac{\tilde{\mu}(3\tilde{\lambda}+2\tilde{\mu})}{\tilde{\lambda}+\tilde{\mu}}\,,\\
\nu & =\cfrac{\tilde{\lambda}}{2(\tilde{\lambda}+\tilde{\mu})}\,.
\end{align}
Accordingly, the Young's modulus is 20$\,$GPa and Poisson's ratio
is 0.33. These values can now be used to set up the numerical model
for the wave propagation analysis.

\begin{figure}[!b]
\centering \subfloat[Time-domain signal $f(t)$]{\includegraphics{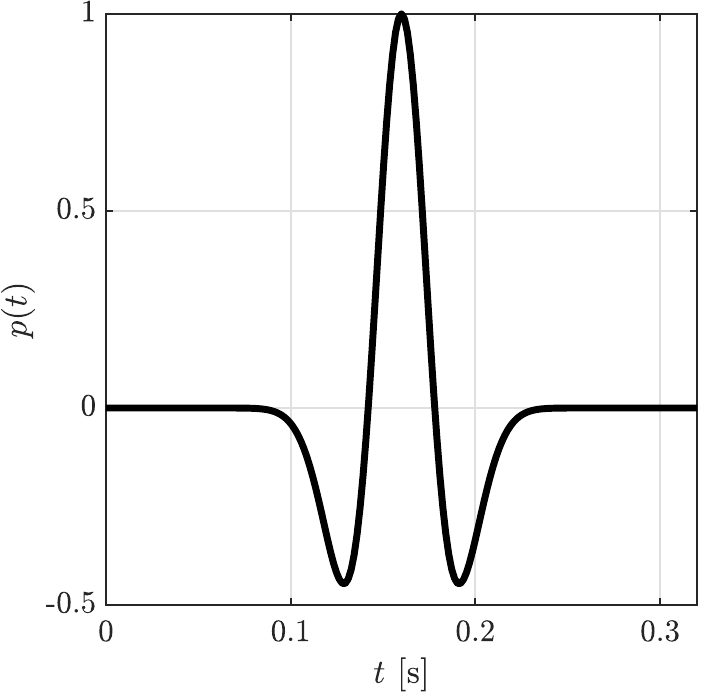}

}\hfill{}\subfloat[Frequency-domain signal $\hat{F}(f)$]{\includegraphics{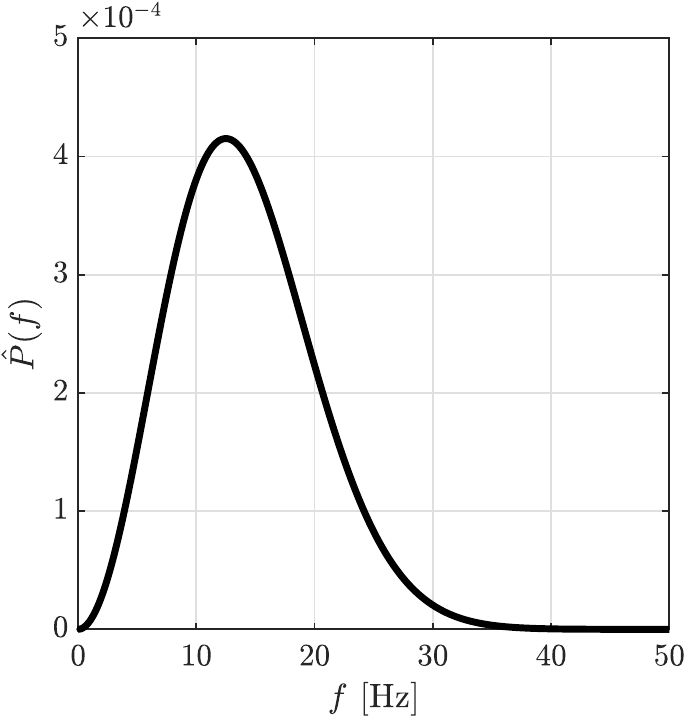}

}\caption{Excitation signal in the time and frequency domain---Lamb's problem.
\label{fig:ExcitationSignalLamb}}
\end{figure}

The wave packets are excited by means of a concentrated point load
located at the upper-left corner of the domain. The time-dependent
excitation follows a Ricker wavelet function 
\begin{equation}
F_{\mathrm{R}}(t)=F_{0}\left(1-\Psi(t)\right)\exp^{-\sfrac{1}{2}\Psi(t)}\qquad\text{with}\qquad\Psi(t)=2\left(\pi f_{\mathrm{ex}}(t-t_{0})\right)^{2}\,,
\end{equation}
where $f_{\mathrm{ex}}$ is the center frequency of the excitation,
$F_{0}$ is the prescribed amplitude, and $t_{0}$ is a time parameter.
In this example, the following values are chosen: $f_{\mathrm{ex}}\,{=}\,12.5\,$Hz,
$F_{0}\,{=}\,100\,$N, and $t_{0}\,{=}\,\sfrac{2}{f_{\mathrm{ex}}}$
and the time- as well as frequency-domain signal are depicted in Fig.~\ref{fig:ExcitationSignalLamb}.
The maximum frequency of interest $f_{\mathrm{max}}$ is determined
at the threshold when the amplitude in the frequency-spectrum is constantly
below 1\% of its maximum value. In our specific example, $f_{\mathrm{max}}$
is 34.54$\,$Hz and therefore, a suitable time steps size to start
the investigation is $\Delta t_{\mathrm{max}}\,{=}\,\sfrac{1}{2f_{\mathrm{max}}}\,{\approx}\,1.5\times10^{-2}\,$s.

The Dirichlet boundary conditions are chosen according to Refs.~\citep{ArticleHam2012,ArticleNoh2013},
i.e., the bottom and right edges of the domain are fixed (clamped),
while symmetry boundary conditions are applied to the left edge. The
top edge is assumed to be free and therefore, no additional Dirichlet
or non-homogeneous Neumann boundary conditions are applied.

Considering the spatial discretization, we use bi-linear (4-node)
quadrilateral elements, and the element size is selected as 5$\,$m.
Thus, 409,600 elements are created which corresponds to 821,762 DOFs.
The spatial discretization is set up with respect to the Rayleigh
wavelength being the wave type with the lowest wave velocity. According
to Eq.~\eqref{eq:Wavelength}, this also means that its wavelength
is the smallest, making it the one defining the element size. For
a maximum considered frequency of $f_{\mathrm{max}}\,{=}\,34.54\,$Hz,
we get $\lambda_{\mathrm{R}}\,{=}\,48.38\,$m and thus, roughly ten
elements per wavelength are employed. As mentioned before and shown
in Ref.~\citep{ArticleWillberg2012}, this is an absolute minimum
and would result in an error of above 1\%. However, at this point,
we are not interested in obtaining highly accurate results with a
minimal error due to the spatial discretization, but we want to determine
the error contributed by the temporal discretization and the chosen
time-integrator. Therefore, the spatial resolution is sufficient for
our purposes.

\begin{figure}[p]
\centering \subfloat[$t\,{=}\,0.2508\,$s\label{subfig:Lamb_1}]{\includegraphics[width=0.425\textwidth]{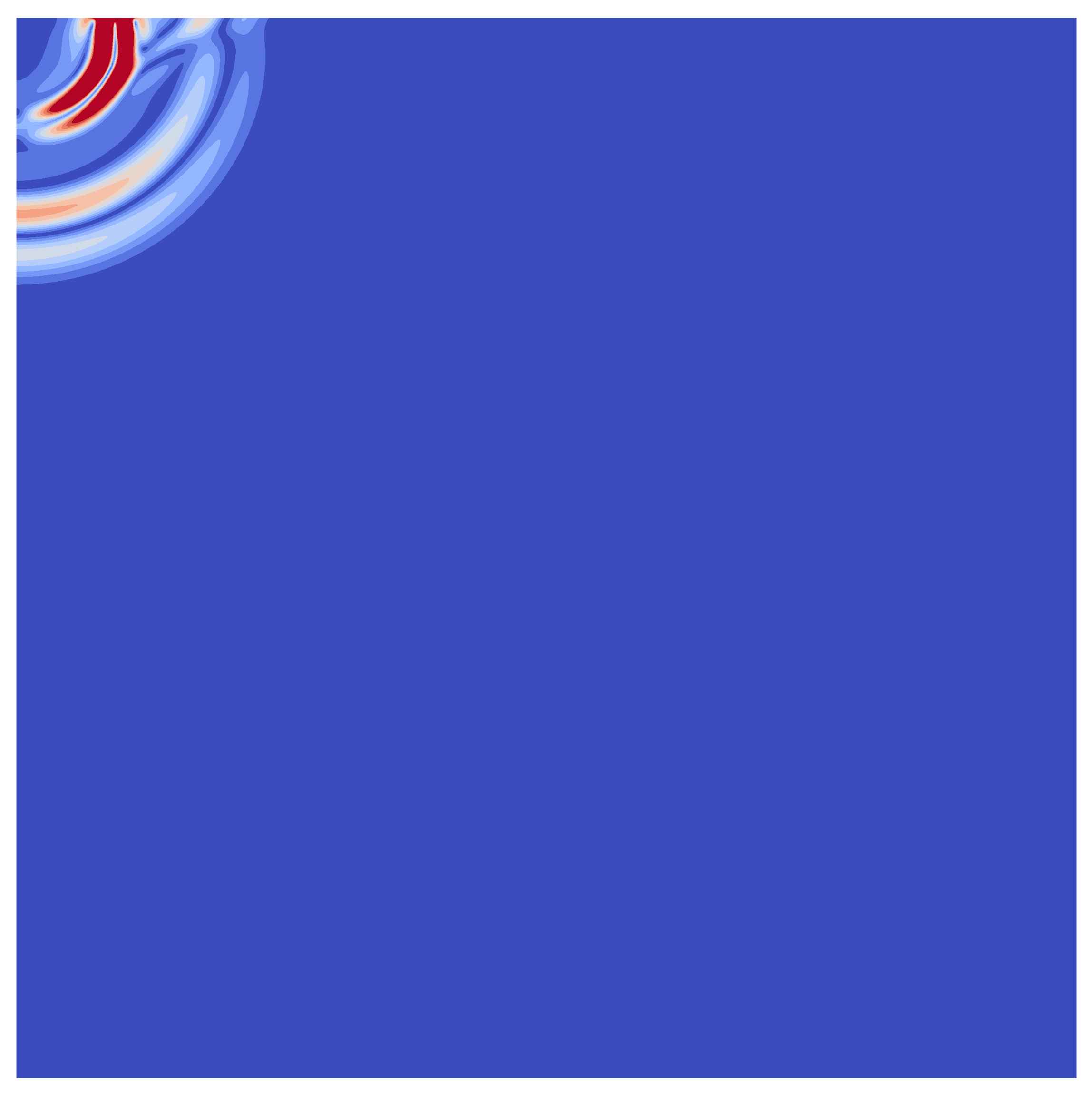}

}\hfill{}\subfloat[$t\,{=}\,0.3511\,$s\label{subfig:Lamb_2}]{\includegraphics[width=0.425\textwidth]{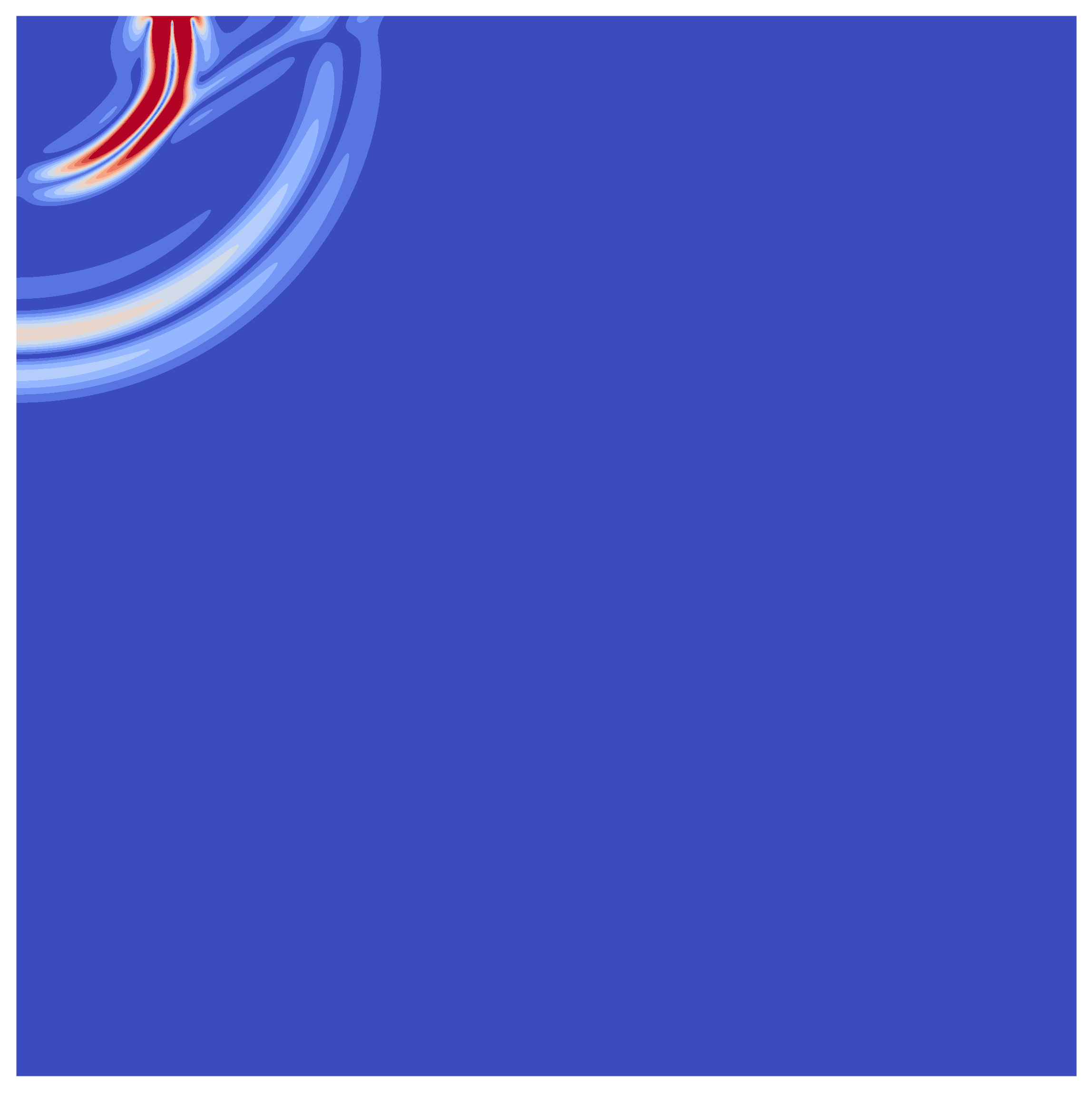}

}\\
 \subfloat[$t\,{=}\,0.4514\,$s\label{subfig:Lamb_3}]{\includegraphics[width=0.425\textwidth]{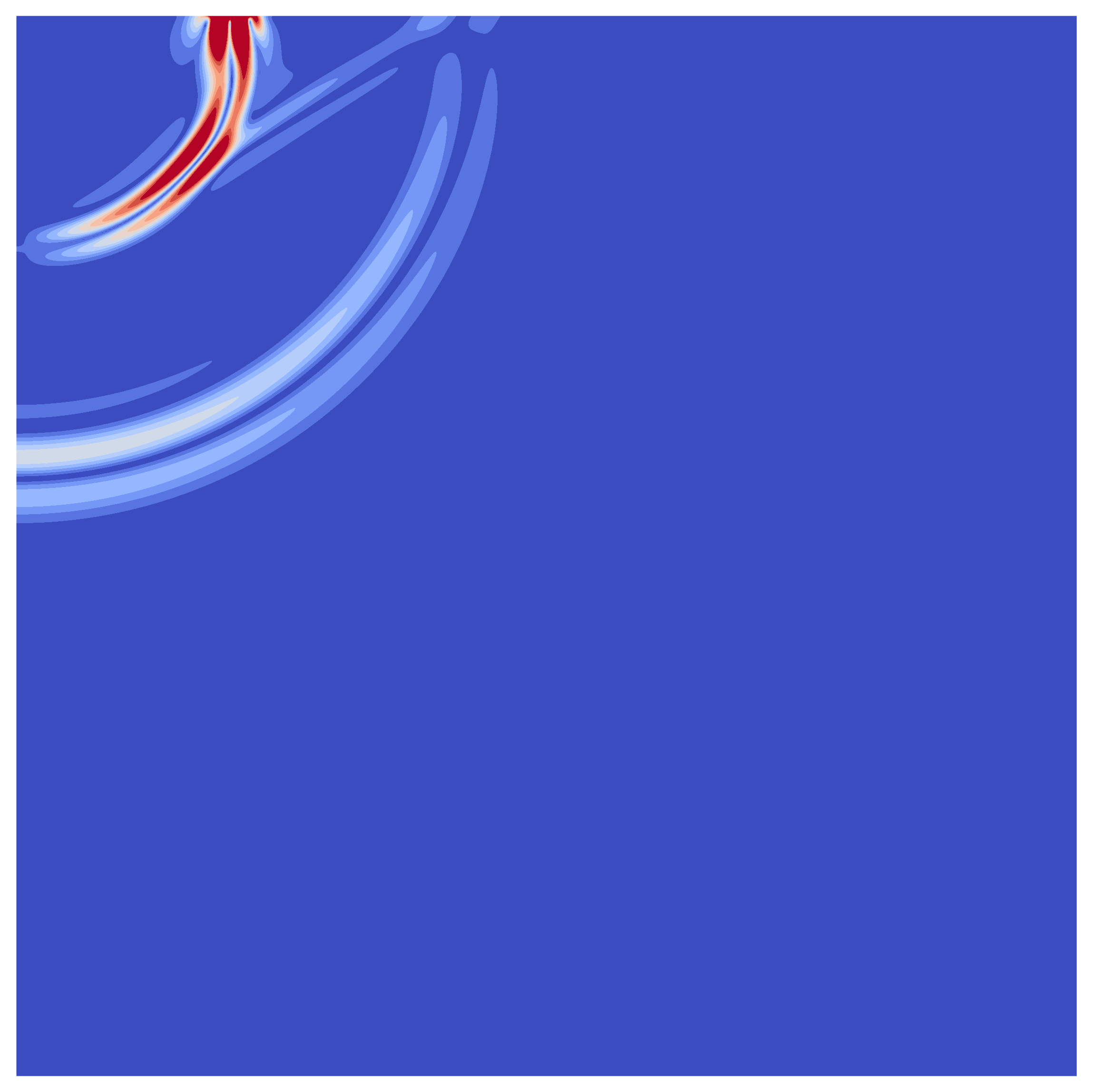}

}\hfill{}\subfloat[$t\,{=}\,0.7523\,$s\label{subfig:Lamb_4}]{\includegraphics[width=0.425\textwidth]{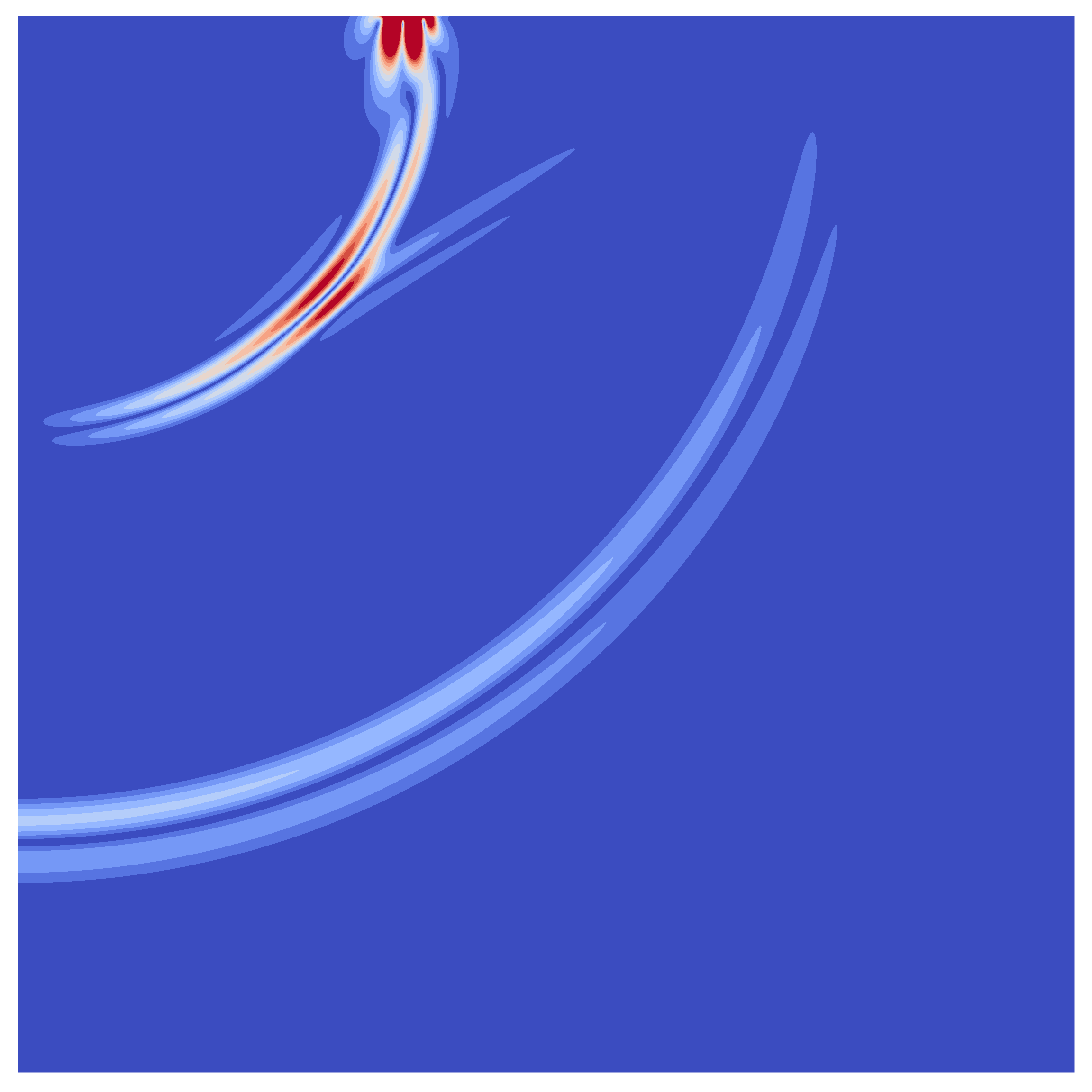}

}\\
 \subfloat[$t\,{=}\,0.8527\,$s\label{subfig:Lamb_5}]{\includegraphics[width=0.425\textwidth]{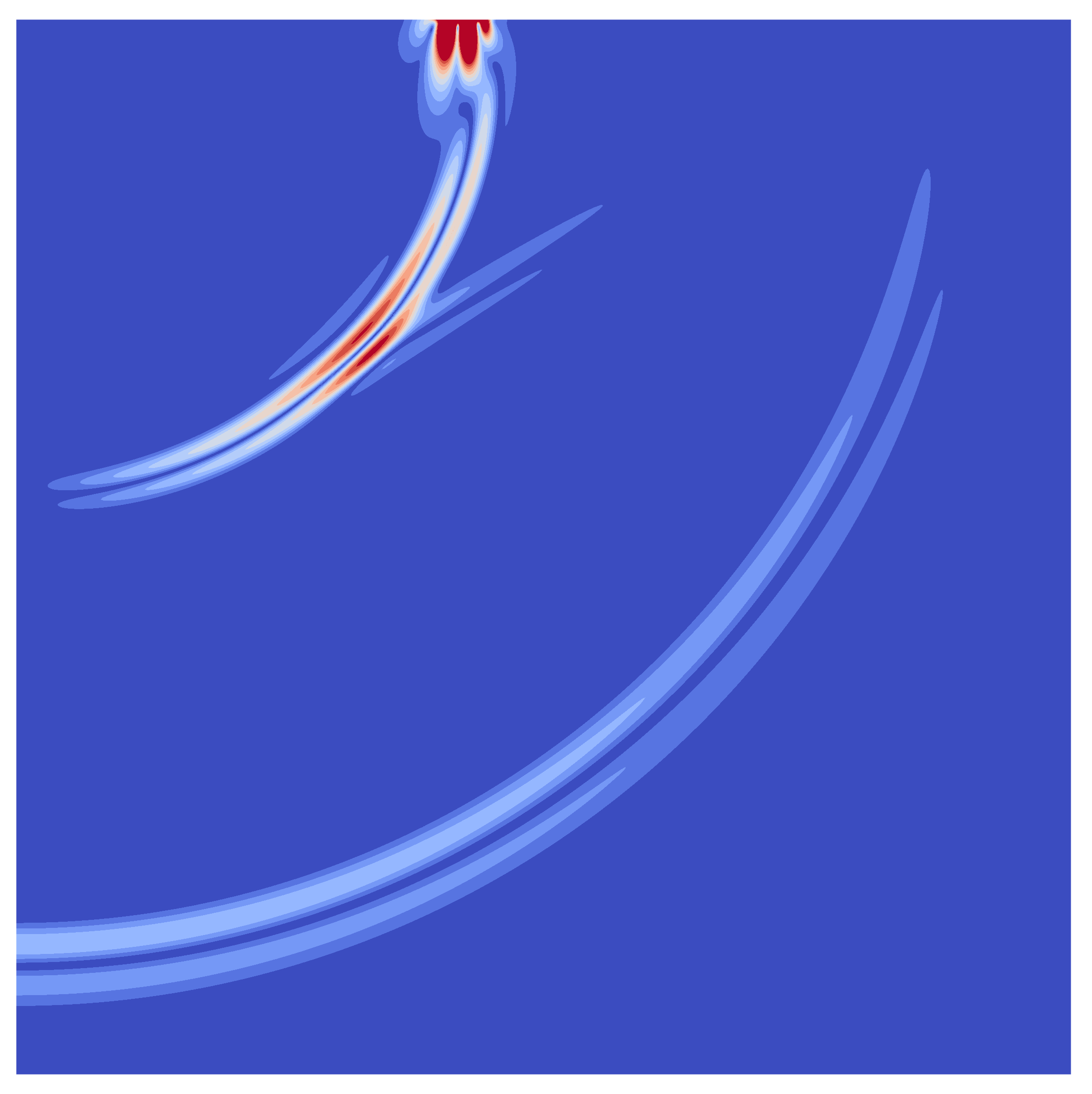}

}\hfill{}\subfloat[$t\,{=}\,0.9530\,$s\label{subfig:Lamb_6}]{\includegraphics[width=0.425\textwidth]{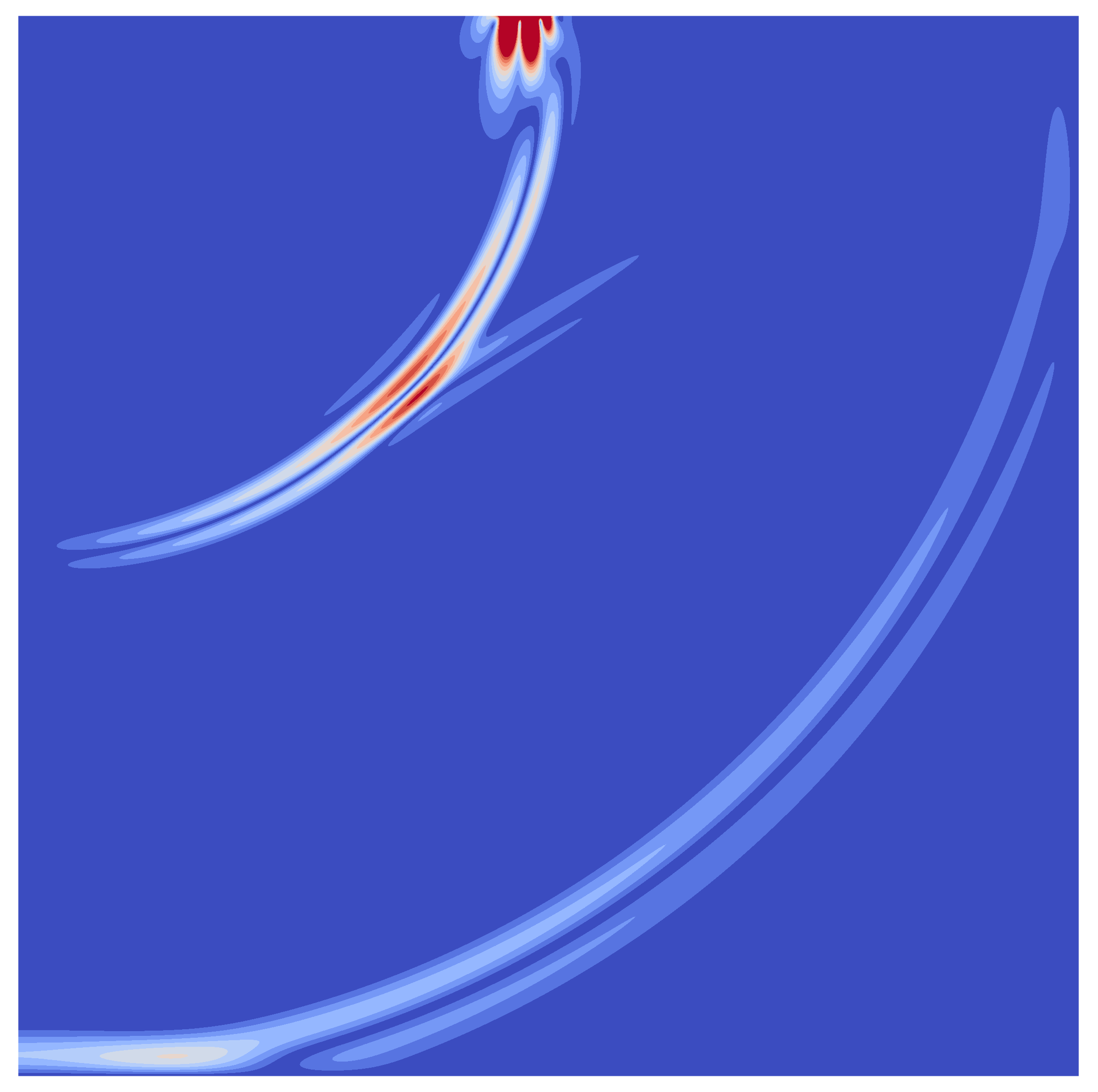}

}\caption{Displacement field ($u_{\mathrm{mag}}$) for Lamb's problem at different
time instances ($\Delta t\,{=}\,2.34375\times10^{-4}\,$s; Present
scheme of order $\mathcal{O}$(4,4). In order to observe not only
the dominant Rayleigh wave the maximum color value of the displacement
is set to $1\times10^{-5}\,$m.\label{fig:DispLamb}}
\end{figure}

\begin{figure}[!t]
\centering \includegraphics{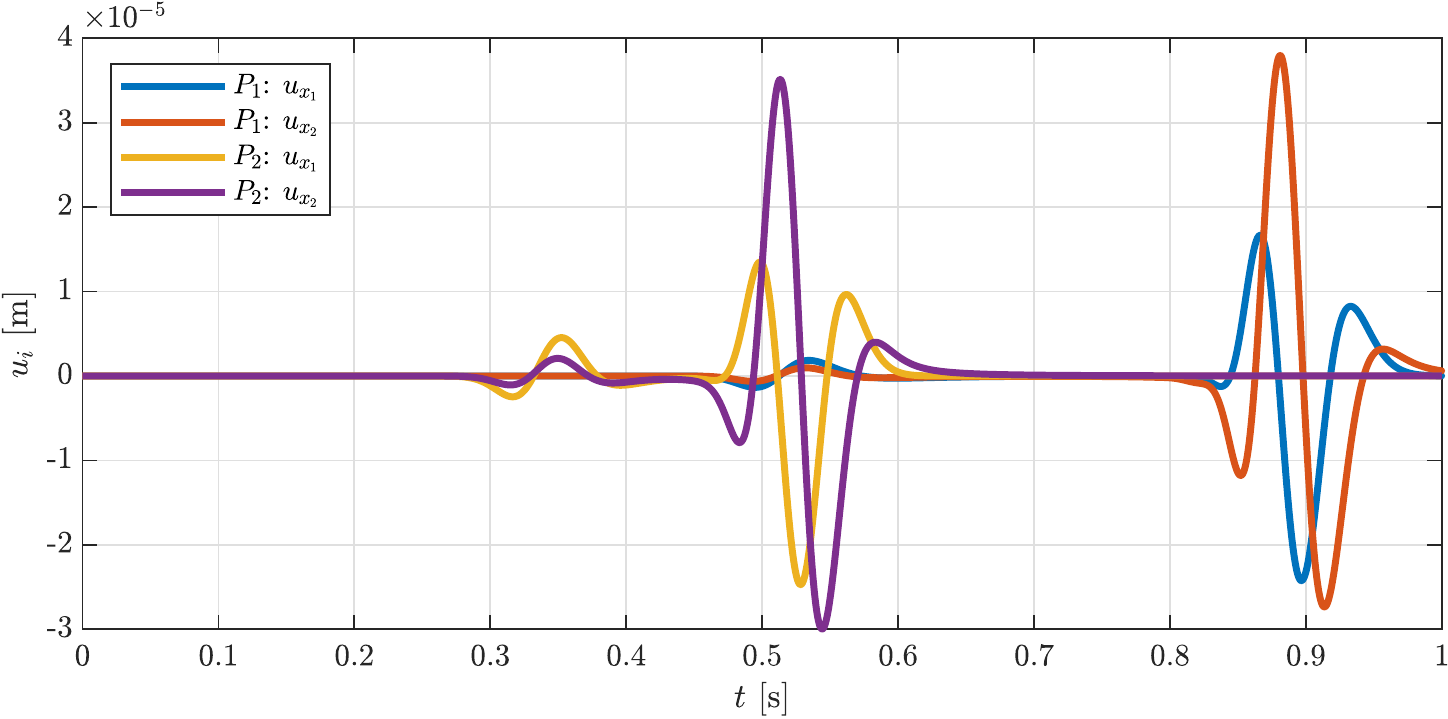} \caption{Displacement response at the observation points $P_{1}$ and $P_{2}$---Lamb's
problem. \label{fig:Lamb_U}}
\end{figure}

The displacement response due to the excitation by means of the point
force following a Ricker wavelet are evaluated at two observation
points $P_{1}$ and $P_{2}$ located at the top surface 640$\,$m
and 1280$\,$m away from the excitation source. Contour plots of the
entire domain for selected time instance, where we clearly observe
the different wave types propagating through the structure are provided
in Fig.~\ref{fig:DispLamb}, while the displacement histories for
the two observation points are given in Fig.~\ref{fig:Lamb_U}.

\subsubsection{Analysis of the computational costs}

\label{sec:CostsLamb} Lamb's problem is only investigated for the
spatial discretization discussed in the previous subsection, i.e.,
an element size of 5$\,$m is selected leading to 821,762 DOFs. Utilizing
a coarser discretization would lead to physically not meaningful results
and is therefore, not conducted. On the other hand, a finer resolution
would yield quickly several million DOFs, which is not feasible for
the performance assessment done in this paper. For studying much larger
systems, it is worth switching from a direct solver to an iterative
solver which is out of the scope of this contribution but is a part
of ongoing research activities. For this example, the normalized computational
time with respect to Newmark's constant average acceleration method
is 1 for the present scheme of order $\mathcal{O}(1,1)$ (which is
mathematically identical), 4.55 for order $\mathcal{O}(2,2)$, 5.53
for order $\mathcal{O}(3,3)$, and 9.96 for order $\mathcal{O}(4,4)$.
These values are moderately higher compared to the ones obtained in
Sect.~\ref{sec:CostsRod}. Still, we have to bear in mind that a
significant increase in accuracy is achieved with elevating the order
of the time-stepping scheme. It is obvious that the numerical overhead
of a high-order technique is dependent on the selected equation solvers
and the numerical model, i.e., sparsity and bandwidth of the matrix,
among others. These factors will be thoroughly investigated in forthcoming
communications. 

\subsubsection{Analysis of the accuracy of the time-integrator}

\label{sec:AccuracyLamb} 
\begin{figure}[b!]
\centering \includegraphics{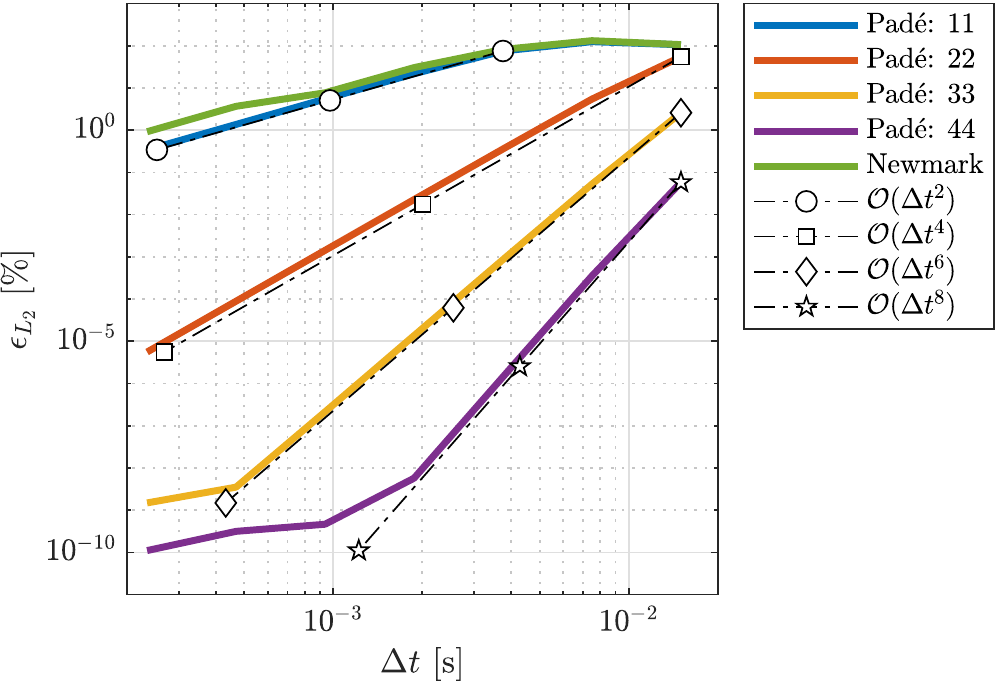} \caption{Displacement error in the $L_{2}$-norm for Lamb's problem. The dash-dotted
lines indicate the optimal rates of convergence corresponding to slopes
of 2 (circle), 4 (square), 6 (diamond), and 8 (pentagram), respectively.
\label{fig:ErrorLamb}}
\end{figure}

In order to assess the accuracy the displacement error in the $L_{2}$-norm
($\epsilon_{L_{2}}$)---see Eq.~\eqref{eq:ErrorL2}---is computed.
In Fig.~\ref{fig:ErrorLamb}, the error is depicted for the chosen
spatial discretization. Similar to the previous numerical examples,
we observe optimal convergence until an error plateau is reached.
As theoretically predicted, the novel time integrator of order $\mathcal{O}(1,1)$
and Newmark's constant average acceleration method yield identical
results. Due to the higher rates of convergence a similar accuracy
can be easily reached using much larger time step sizes when the novel
approach is employed. Considering an error threshold of roughly 1\%,
which is acceptable in most engineering applications, a time increment
of $\Delta t\,{=}\,4.01\times10^{-4}\,$s is required for Newmark's
constant average acceleration method, while for the high-order schemes
significantly larger values are acceptable. Considering the proposed
time integrator of order $\mathcal{O}(2,2)$, $\Delta t\,{=}\,4.87\times10^{-3}\,$s
is sufficient, while the time increment can be further increased for
orders $\mathcal{O}(3,3)$ and $\mathcal{O}(4,4)$ to $\Delta t\,{=}\,1.27\times10^{-2}\,$s
and $\Delta t\,{=}\,2.0\times10^{-2}\,$s, respectively. That is to
say, the time step size can be increased by factors of 12.14, 31.67,
and 49.88 for the novel time-integrator. These values need to be put
into perspective with the increased costs per times step as found
in the previous subsection. There, it was observed that the computational
effort increases by a factor of roughly 4.5, 5.5, and 10 for the fourth-,
sixth- and eighth-order accurate schemes. Therefore, the attainable
speed-up is observed to be 2.67, 5.73, and 5.01. Again, please keep
in mind that this value is problem dependent and might also be related
to the properties of the numerical method that is employed for the
spatial discretization. 
\begin{figure}[t!]
\centering \includegraphics{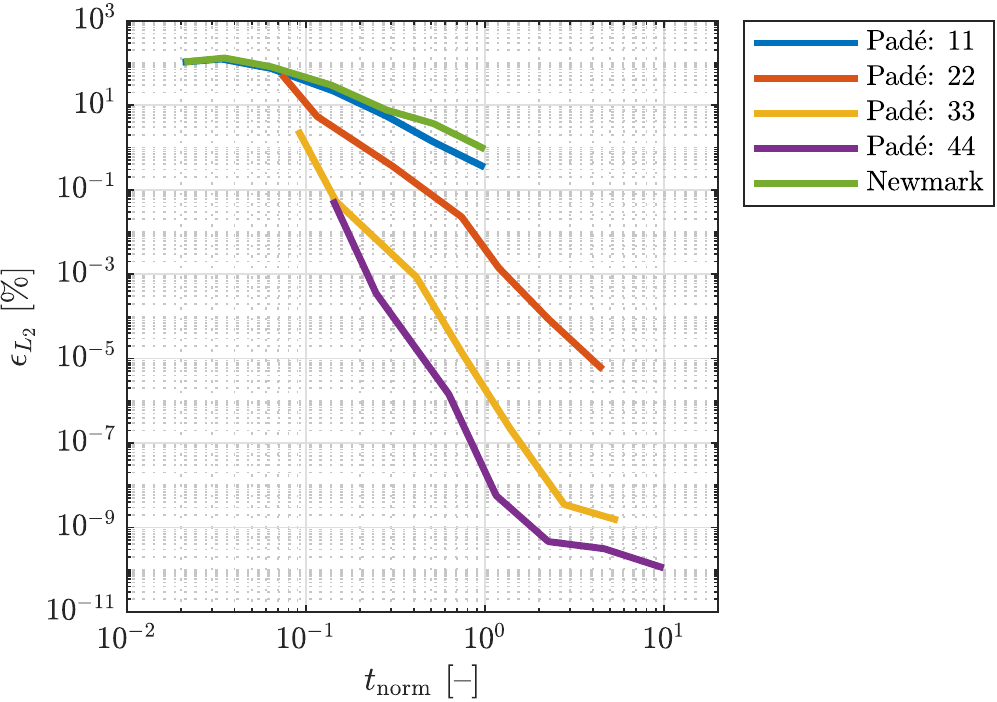} \caption{Displacement error in the $L_{2}$-norm for the different configurations
over the normalized computational time---Lamb's problem. \label{fig:ErrorTcpuLamb}}
\end{figure}

The conclusions that can be drawn from the findings reported in the
current section are summarized in Fig.~\ref{fig:ErrorTcpuLamb}.
Here, the attainable accuracy is plotted against the normalized computational
time. As a reference, the computational time corresponding to Newmark's
method with a time step size of $\Delta t_{\mathrm{max}}\,{=}\,2.34375\times10^{-4}\,$s
is taken. The results again illustrate that employing a high-order
time integration scheme leads to significant savings in terms of computational
time for a prescribed error threshold or to significantly more accurate
results for the same computational time.

\subsection{Wave propagation in a mountainous region}

\label{sec:Mountain} 
\begin{figure}[!t]
\centering \subfloat[Time-domain signal $f(t)$]{\includegraphics{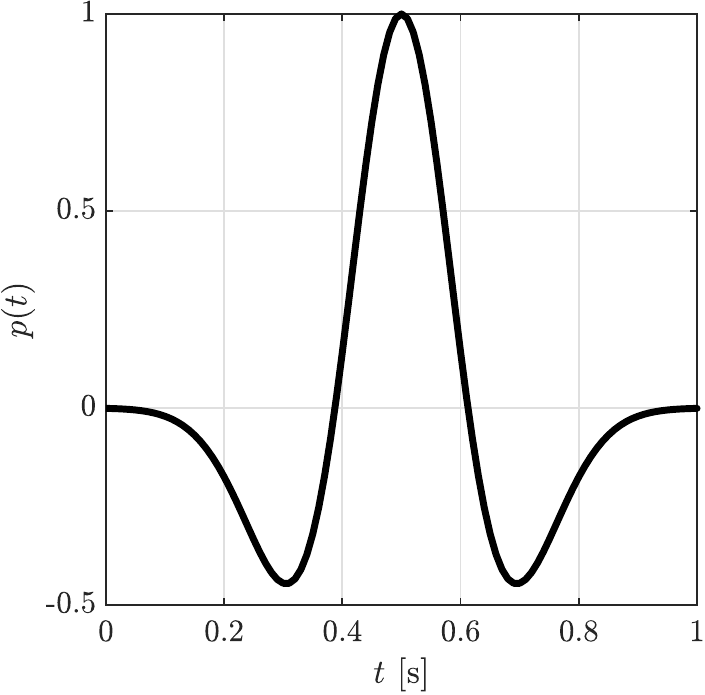}

}\hfill{}\subfloat[Frequency-domain signal $\hat{F}(f)$]{\includegraphics{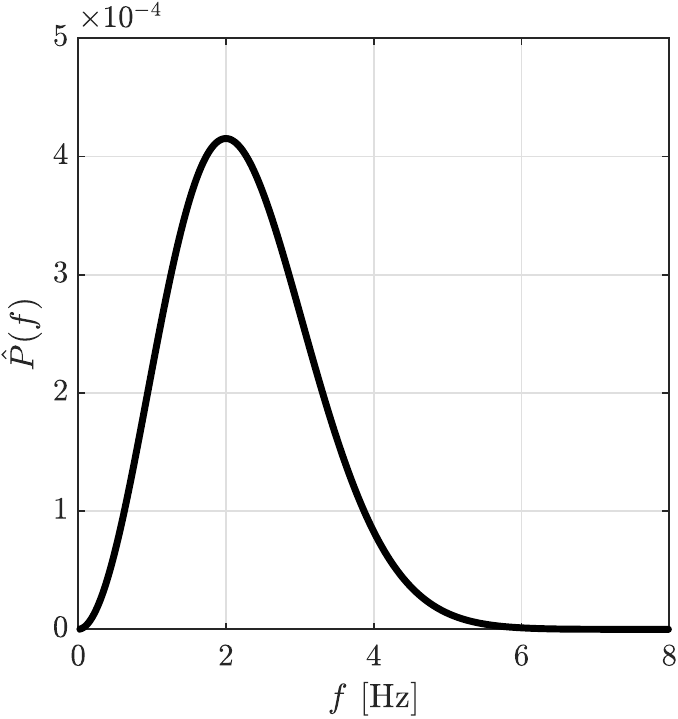}

}\caption{Excitation signal in the time and frequency domain---Mountainous
region. \label{fig:ExcitationSignalMawson}}
\end{figure}

As a final example, a three-dimensional mountainous
region is chosen\footnote{The simulation has been run on a workstation with the following specifications:
Precision 7920-Tower Workstation; Intel(R) Xeon(R) CPU E5-2637 v4
@ 3.50$\,$GHz; 512$\,$GB DDR4 (8$\times$64$\,$GB), 2400 MHz; NVIDIA
Quadro K1200.}. Here, we are simulating the wave propagation in a mountainous region
of size 6.94$\,$km$\times$7.8$\,$km to showcase the performance
of the proposed time integration scheme for fully three-dimensional
models. The initial STL file\footnote{Upon request, the authors are open to share the model files for reproduction
of the published results and further analyses.} of the region is obtained using the online service Terrain2STL at
http://jthatch.com/Terrain2STL/. The material properties are assumed
as: Young's modulus $E\,{=}\,55\,$GPa, Poisson's ratio $\nu\,{=}\,0.2$,
and the mass density $\rho\,{=}\,2400\,$kg/m$^{3}$. Thus, the P-wave
and S-wave velocities are $c_{\mathrm{L}}\,{=}\,5046\,$m/s and $c_{\mathrm{T}}\,{=}\,3090\,$m/s,
respectively. In terms of boundary conditions, it is assumed that
the displacements normal to the side and bottom faces are fixed. The
excitation is realized by means of a concentrated force in $x_{2}$-direction
applied at the centroid of the bottom surface. The excitation follows
a Ricker wavelet with a center frequency of $f_{\mathrm{ex}}\,{=}\,2\,$Hz,
while the time parameter is set to $t_{0}\,{=}\,\sfrac{1}{f_{\mathrm{ex}}}$
(see Fig.~\ref{fig:ExcitationSignalMawson} for the time- and frequency-domain
signal). The maximum frequency of interest $f_{\mathrm{max}}$ is
determined at the threshold when the amplitude in the frequency-spectrum
is constantly below 1\% of its maximum value. In our specific example,
$f_{\mathrm{max}}$ is 5.53$\,$Hz and therefore, a suitable time
steps size to start the investigation is $\Delta t_{\mathrm{max}}\,{=}\,\sfrac{1}{2f_{\mathrm{max}}}\,{\approx}\,9\times10^{-2}\,$s.

\begin{figure}[!t]
\centering \subfloat[STL data]{\includegraphics[width=0.475\textwidth]{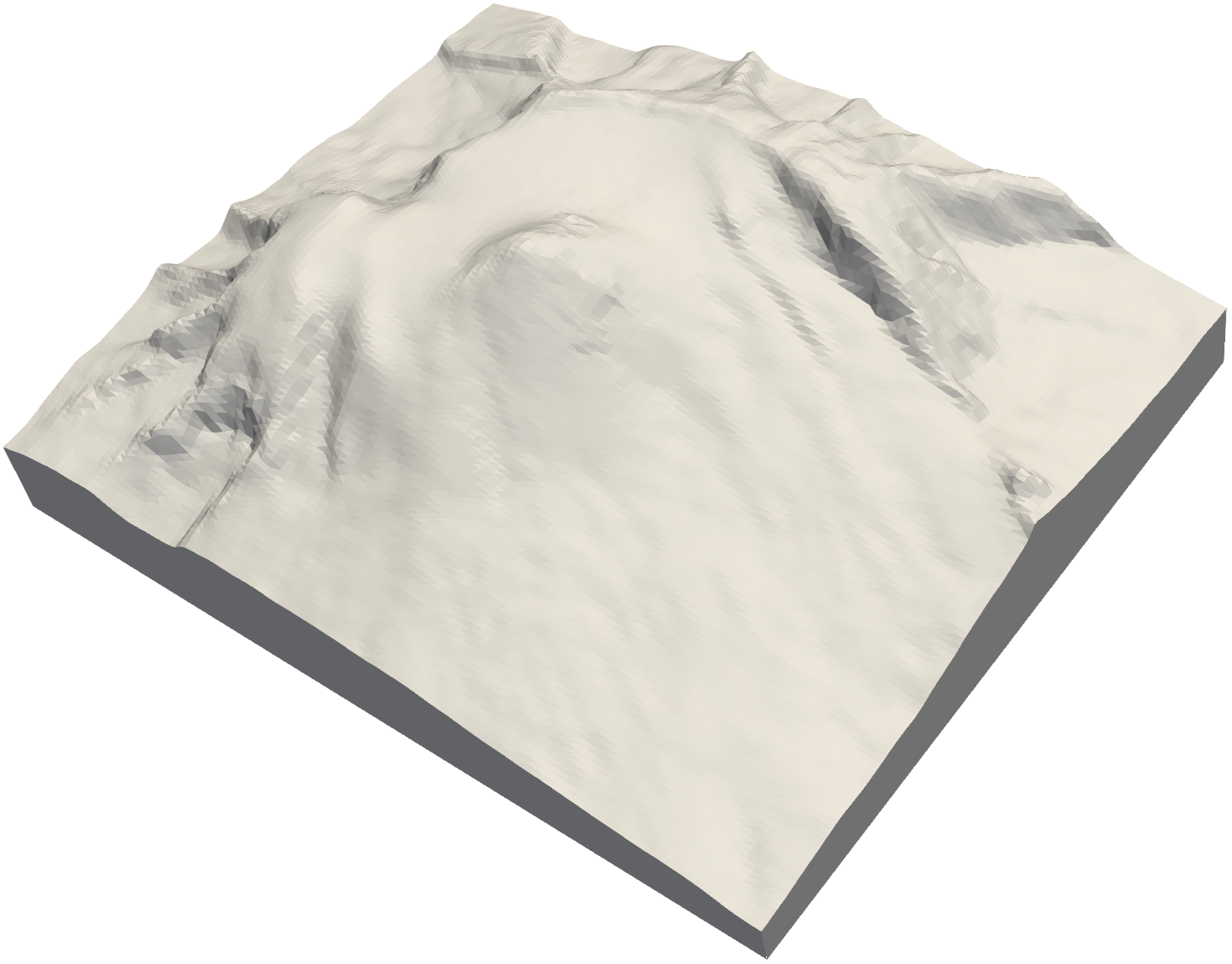}

}\hfill{}\subfloat[Octree discretization]{\includegraphics[width=0.475\textwidth]{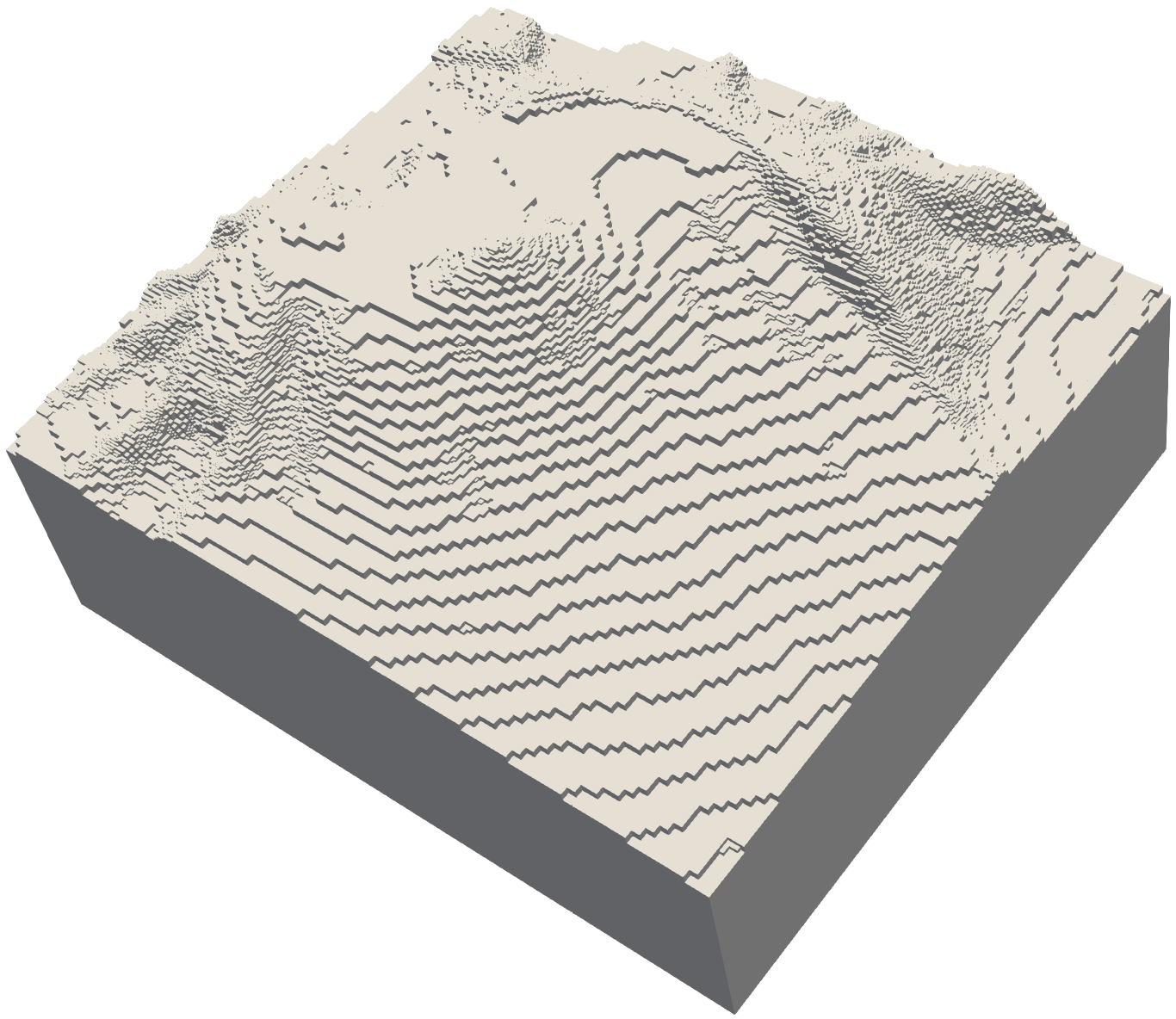}

}\\
 \subfloat[Cut view of the discretization]{\includegraphics[width=0.6\textwidth]{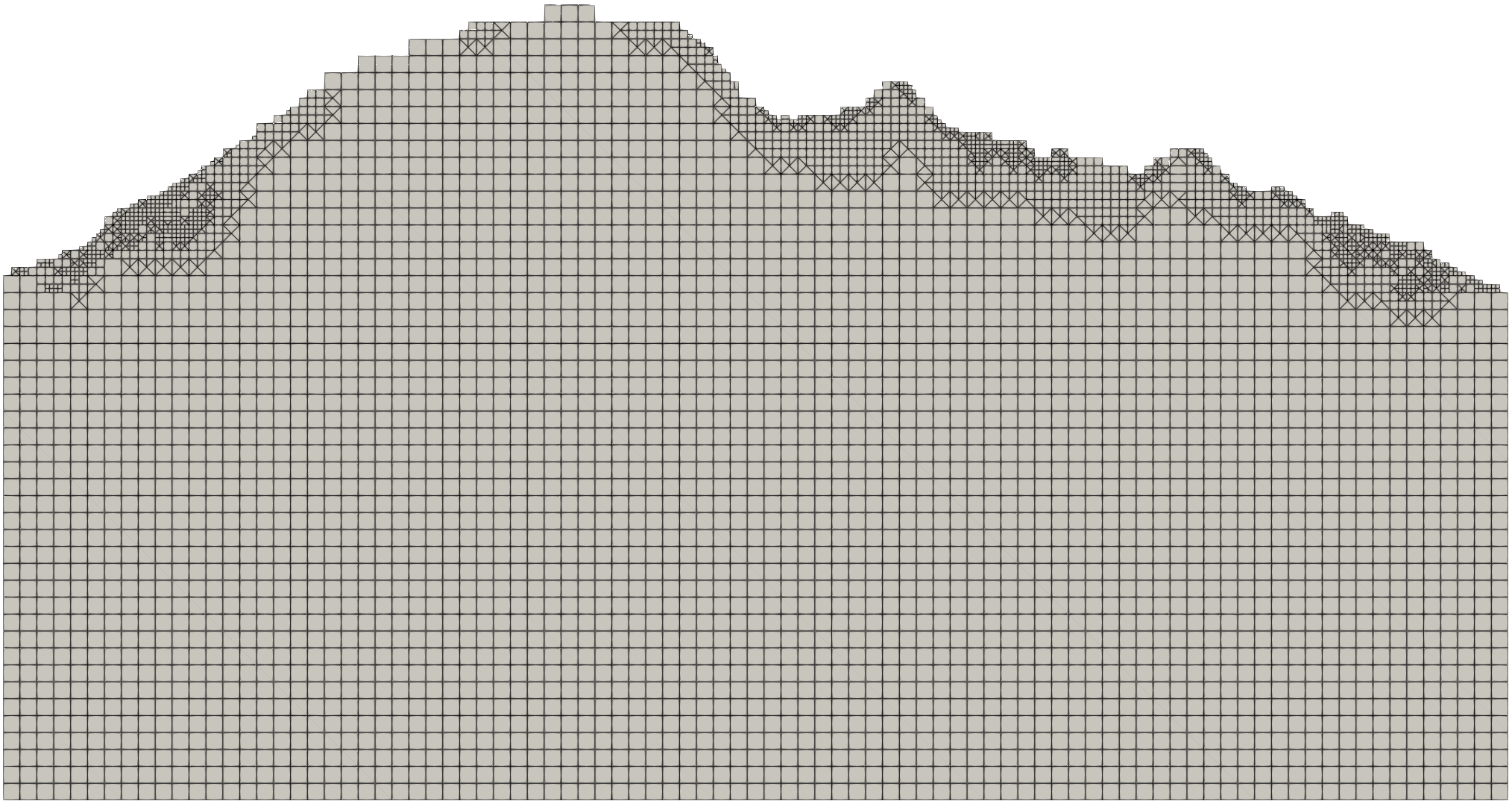}

}\caption{Three-dimensional model of a mountainous region. \label{fig:ModelMawsonPeak}}
\end{figure}

The displacement response has been saved at the peak of the mountainous
region (3,224$\,$m, 4,160$\,$m, 3,952$\,$m) and at other points
with 75\% (3.224$\,$m, 4,160$\,$m, 2,912$\,$m), 50\% (3.224$\,$m,
4,160$\,$m, 1,976$\,$m), and 25\% (3,224$\,$m, 4,160$\,$m, 0,936$\,$m)
of the height. Since a highly complex octree mesh has been used the
nodes nearest to those locations have been picked. The original geometry
of the region and its discretization are depicted in Fig.~\ref{fig:ModelMawsonPeak}.

The simulations in this section are based the scaled boundary finite
element method (SBFEM), which is a semi-analytical non-standard discretization
technique. For a comprehensive introduction and thorough assessment
of the performance of SBFEM, the interested reader is referred to
the monograph by Song \citep{BookSong2018}. At this point, we only
want to mention that using the SBFEM it is straightforward to set
up octree meshes based on STL files without encountering the ``hanging
nodes'' problem \citep{ArticleSaptura2017,ArticleLiu2017}. Here,
each element is modeled as a scaled boundary polyhedron and therefore,
a conformal/compatible coupling between elements of different sizes
is \emph{a priori} guaranteed, i.e., there is no need for a special
treatment as required in conventional FEA. In the current analysis,
a simple octree mesh without trimming has been used for the sake of
simplicity, but it is also possible to achieve geometry conforming
meshes \citep{ArticleTalebi2016,ArticleZhang2019} to improve the
quality of the solution. Overall, the mesh consists of 173,786 octree
cells (featuring linear shape functions on the boundary) with a ratio
of 1/4 for the smallest (cell size: 26$\,$m) to the largest cell
(cell size: 104$\,$m). Thus, we have 251,779 nodes in the discretization
and consequently, 755,337 DOFs. In order to be able to obtain an accurate
reference solution in a relatively short time, the central difference
method in conjunction with a lumped (diagonal) mass matrix is applied.
The time increment is set to $\Delta t_{\mathrm{CDM}}\,{=}\,1\times10^{-4}\,$s.
More details on the particular implementation of the CDM are discussed
in Ref.~\citep{ArticleZhang2021}.

\begin{figure}[p]
\centering \subfloat[$t\,{=}\,0.6058\,$s\label{subfig:Mawson_1}]{\includegraphics[width=0.45\textwidth]{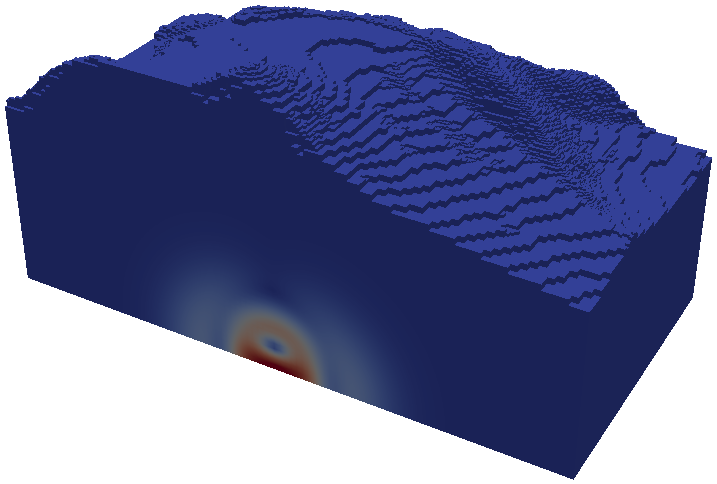}

}\hfill{}\subfloat[$t\,{=}\,0.8078\,$s\label{subfig:Mawson_2}]{\includegraphics[width=0.45\textwidth]{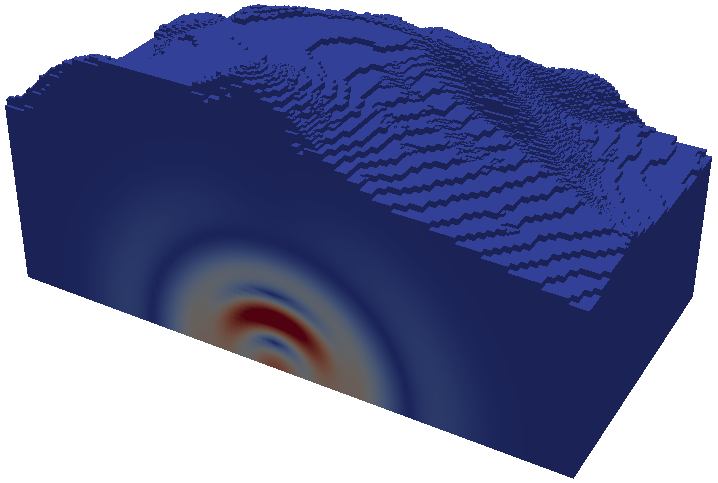}

}\\
 \subfloat[$t\,{=}\,1.0097\,$s\label{subfig:Mawson_3}]{\includegraphics[width=0.45\textwidth]{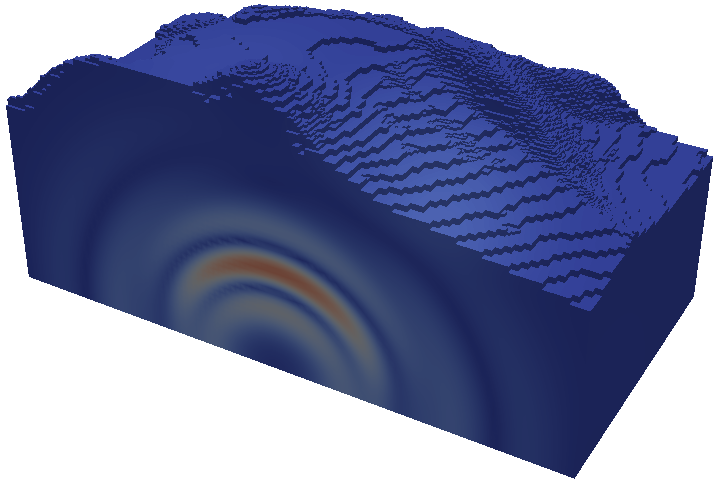}

}\hfill{}\subfloat[$t\,{=}\,1.2117\,$s\label{subfig:Mawson_4}]{\includegraphics[width=0.45\textwidth]{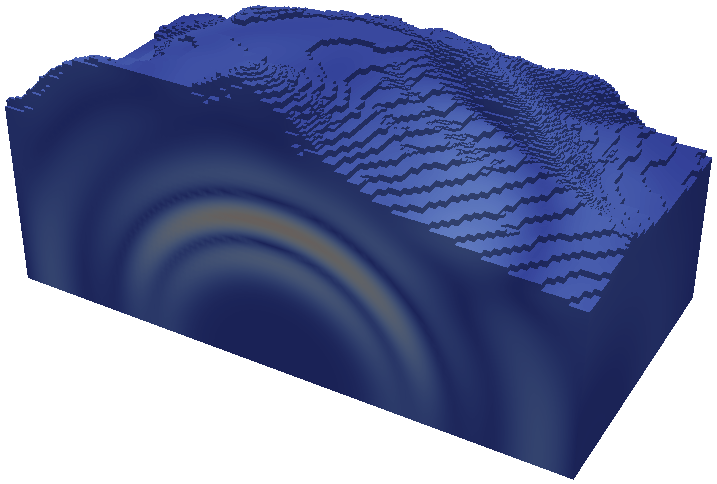}

}\\
 \subfloat[$t\,{=}\,1.4136\,$s\label{subfig:Mawson_5}]{\includegraphics[width=0.45\textwidth]{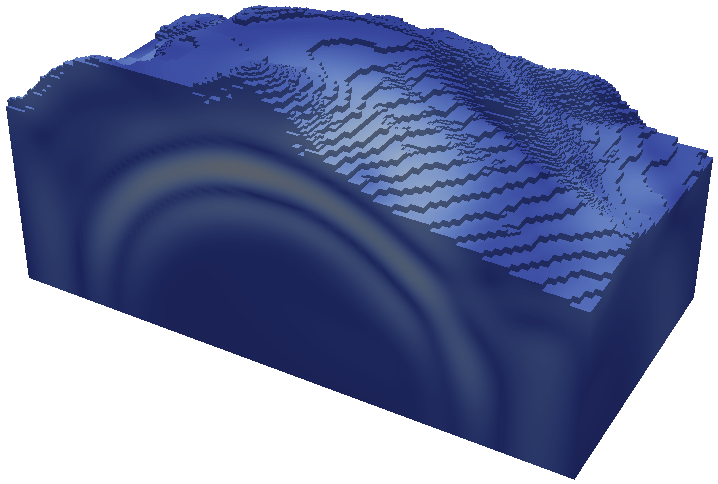}

}\hfill{}\subfloat[$t\,{=}\,1.6156\,$s\label{subfig:Mawson_6}]{\includegraphics[width=0.45\textwidth]{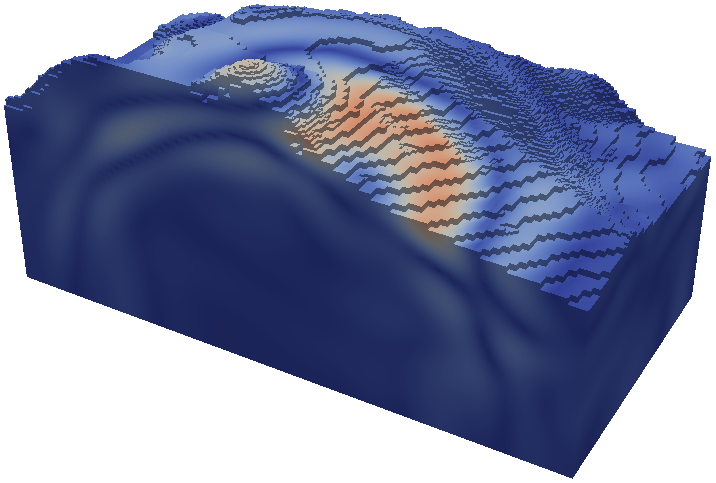}

}\\
 \subfloat[$t\,{=}\,1.8175\,$s\label{subfig:Mawson_7}]{\includegraphics[width=0.45\textwidth]{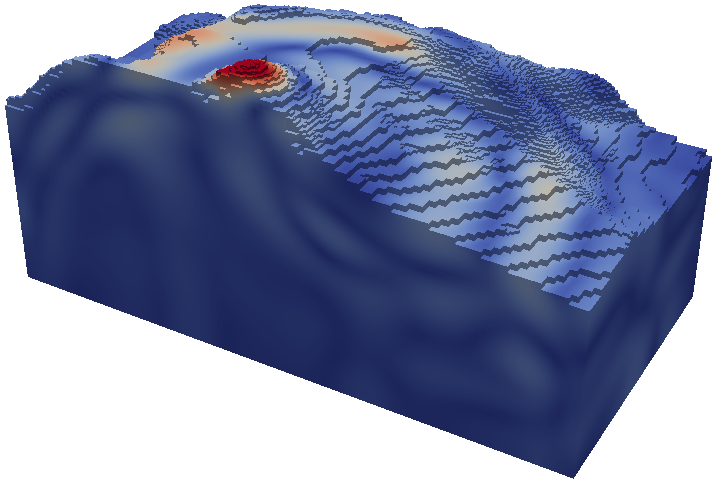}

}\hfill{}\subfloat[$t\,{=}\,2.0194\,$s\label{subfig:Mawson_8}]{\includegraphics[width=0.45\textwidth]{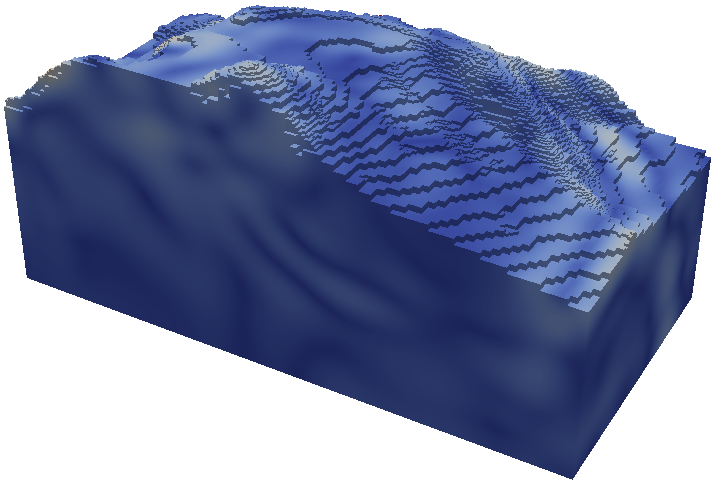}

}\caption{Displacement field ($u_{\mathrm{mag}}$) for the mountainous region
at different time instances ($\Delta t\,{=}\,0.0018\,$s; Padé scheme
of order $\mathcal{O}$(1,1). For a better visibility of the propagating
waves, the maximum color value of the displacement is set to $5\times10^{-5}\,$m.\label{fig:DispMawson}}
\end{figure}

In this three-dimensional example, we will compute the transient response
using the novel time-stepping scheme of order $\mathcal{O}$(4,4)
with the time step size $\Delta t_{\mathrm{max}}$. Additionally,
simulations employing Newmark's constant average acceleration method
will be executed with time step sizes of $\Delta t_{\mathrm{max}}$
and $\sfrac{\Delta t_{\mathrm{max}}}{50}$. Judging for the findings
reported in the previous sections, a clear difference in the displacement
history is expected for a time increment $\Delta t_{\mathrm{max}}$.
However, when using the smaller time step size for Newmark's method
a significantly improved agreement is expected. Using this approach,
we can gain some insights into the performance of the time integrator
for complex three-dimensional problems. Extensive parametric studies
for such examples are still computationally very expensive and prohibitive
at this point. Snapshots of the wave propagation in the mountainous
region for different time instances are depicted in Fig.~\ref{fig:DispMawson}.
These contour plots illustrate the propagation of the primary wave
through the mountainous region. In further, more sophisticated, analyses,
numerical models of faults, different rock layers, and unbounded domains
can be straightforwardly added.

\begin{figure}[!p]
\centering \subfloat{\includegraphics{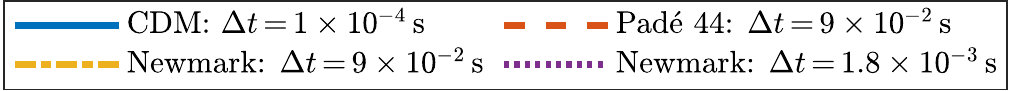}}\\
 \setcounter{subfigure}{0} \subfloat[Displacement response ($u_{\mathrm{{x_{2}}}}$) at point $P_{1}$(peak)
\label{subfig:P1_Mountain}]{\includegraphics{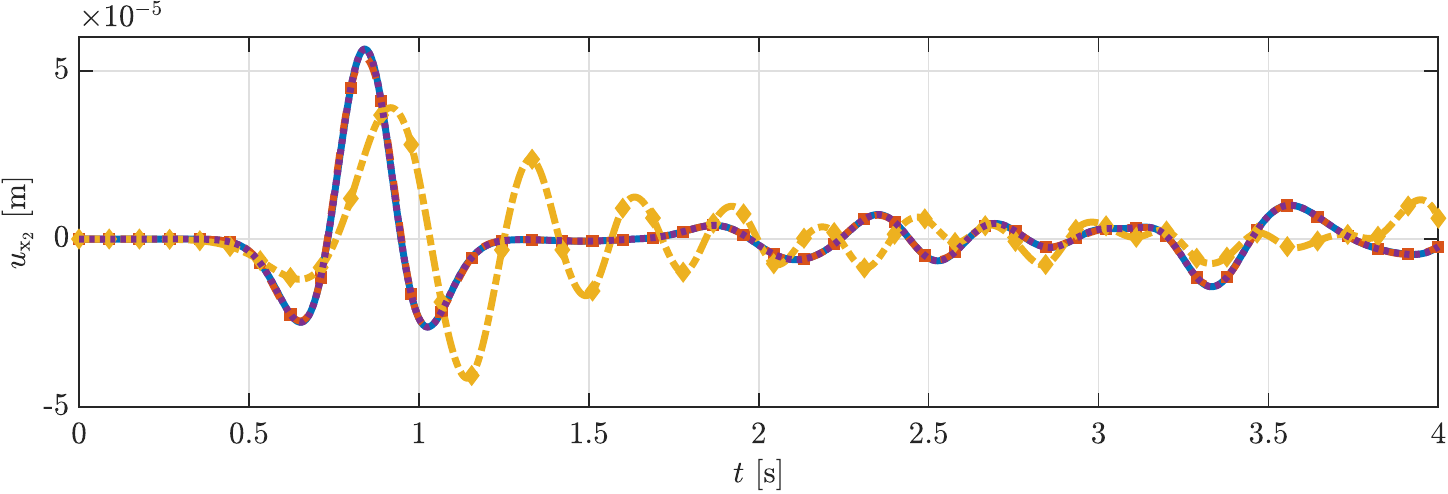}

}\hfill{}\subfloat[Displacement response ($u_{\mathrm{{x_{2}}}}$) at point $P_{2}$(75\%
of the peak height)\label{subfig:P2_Mountain}]{\includegraphics{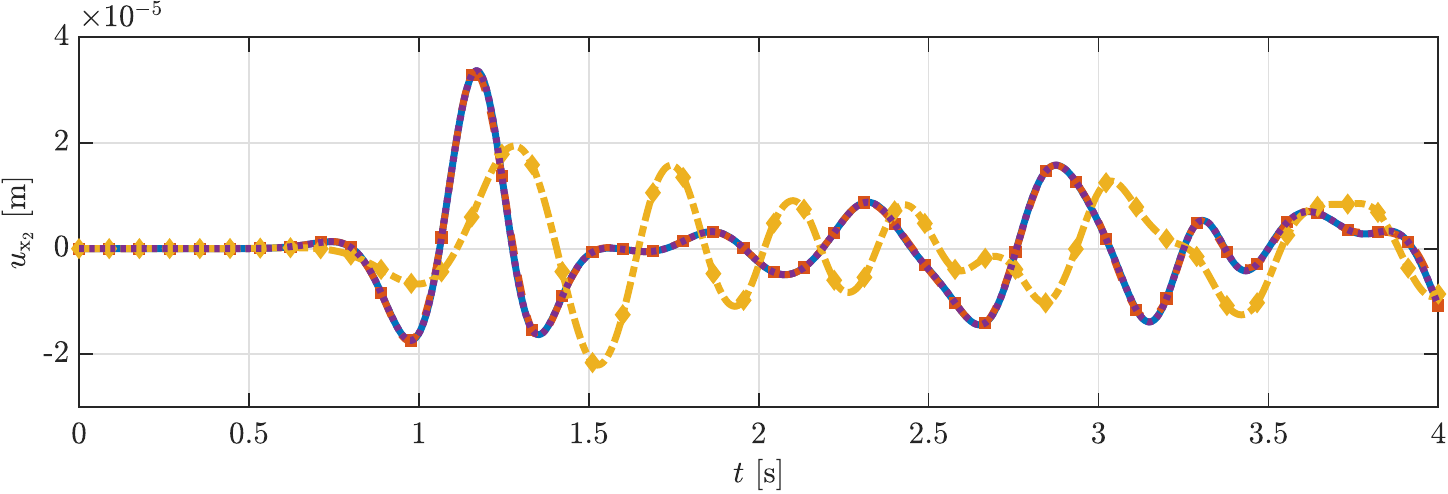}

}\\
 \subfloat[Displacement response ($u_{\mathrm{{x_{2}}}}$) at point $P_{3}$(50\%
of the peak height)\label{subfig:P3_Mountain}]{\includegraphics{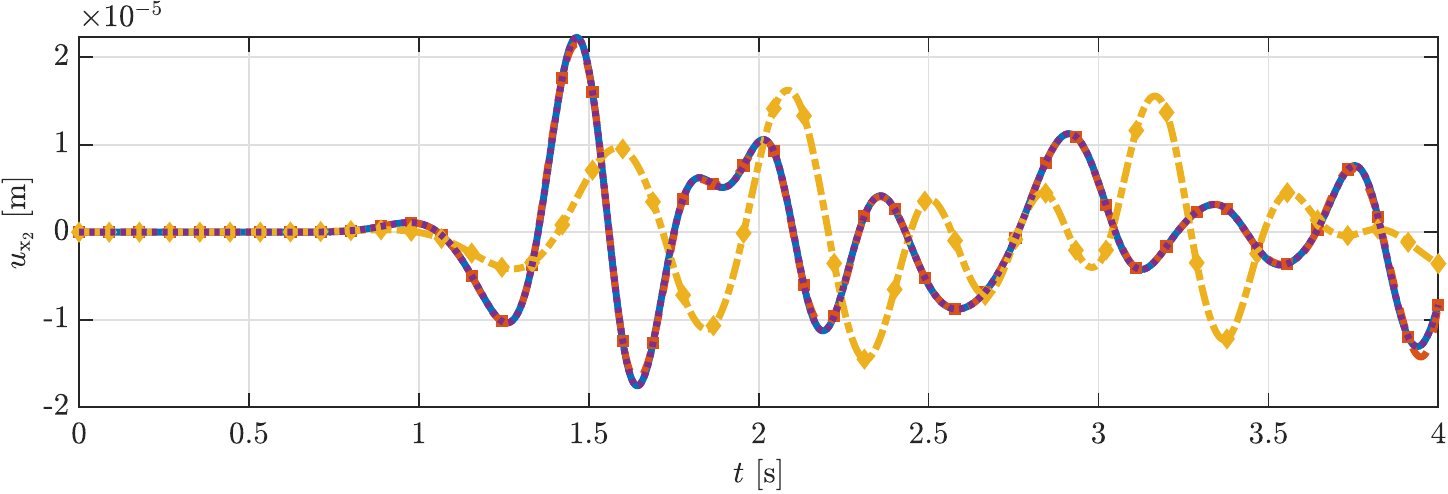}

}\hfill{}\subfloat[Displacement response ($u_{\mathrm{{x_{2}}}}$) at point $P_{3}$(25\%
of the peak height)\label{subfig:P4_Mountain}]{\includegraphics{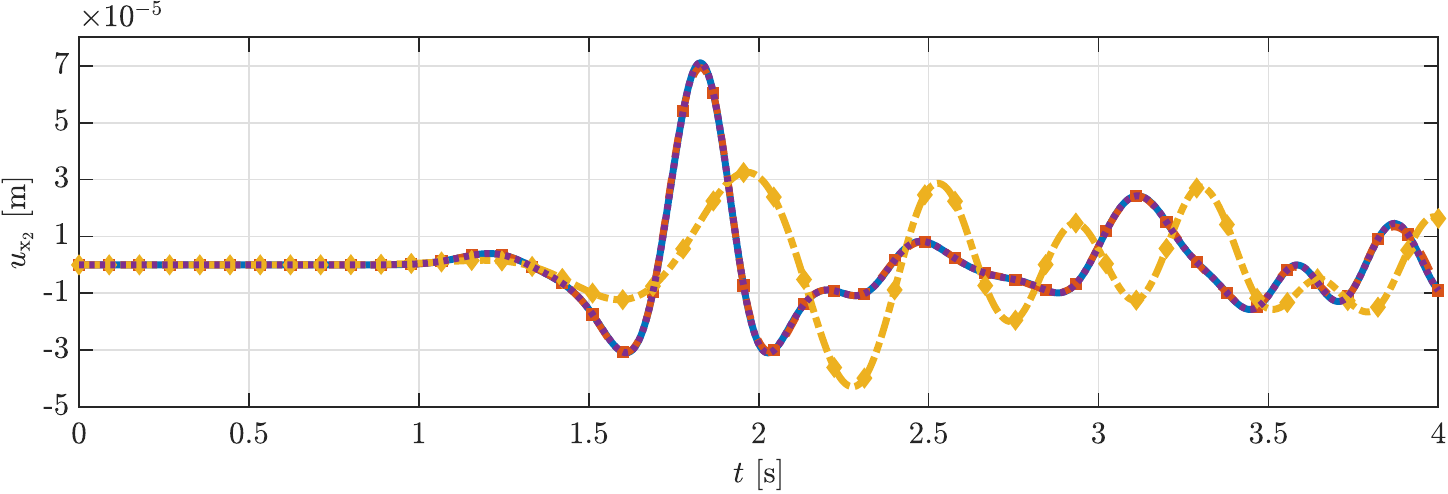}

}\caption{Displacement history at the four observation points. \label{fig:Mawson_U}}
\end{figure}

The displacement field $u_{x_{2}}$ at the four observation points
is depicted in Fig.~\ref{fig:Mawson_U}. Here, we observe clear differences
in the response if the same time step size is used for both Newmark's
constant average acceleration method and the novel eighth-order scheme.
The results obtained employing the proposed time-stepping scheme are
very accurate and closely match the reference results. Since a very
large time step has been chosen, we included markers to indicate the
actual time points that have been computed. To obtain a smooth curve,we
interpolated between these few points using $C_{1}$-contnuous cubic
splines. Therefore, the small observed differences are not actually
related to the time-integration itself but to the required interpolation.
In contrast, if the time increment is significantly decreased for
Newmark's method (in this case by a factor of 50), also an excellent
agreement of the results is noted compared to the reference solution.
This result emphasizes again that highly accurate solution are obtained
with the novel method for very large time step sizes. Considering
the computational time, we normalize the values with respect to Newmark's
constant average acceleration method with the time increment $\Delta t_{\mathrm{max}}$.
The novel eighth-order time-stepping algorithm only requires 4.05
times the computational time, while decreasing the time step size
by a factor of 50 an increase in the consumed time by a 24.1 is observed.
Note that for the current example, a large portion of the solution
time is invested in the initial factorization of effective system
matrices.

\section{Conclusions and outlook on future research activities}

\label{sec:Conclusion} In the present article, a high-order implicit
time-stepping scheme based on a Padé series expansion of the matrix
exponential function has been proposed. The salient features of this
time integrator are: (i) accuracy of order $2p$, where $p$ is the
polynomial degree of the diagonal Padé expansion; (ii) negligible
period elongation for high-order schemes; (iii) no amplitude decay
(no numerical damping); and (iv) high efficiency in the solution of
the semi-discrete equations of motion. The proposed method is especially
suited for applications where a high accuracy of the results is of
utmost importance and long-term analyses are required. Although, it
is generally argued that high-frequency numerical damping is a wanted
property in time integrators, there are also applications where any
damping (physical or numerical) is unwanted and would deteriorate
the accuracy of the analysis. If, however, spurious waves are excited
and damping is indispensable, preliminary results of the authors suggest
that the addition of a very low amount of stiffness-proportional damping
(i.e., $\beta_{\mathrm{R}}\,{<<}\,1$ and $\alpha_{\mathrm{R}}\,{=}\,0$)
is adequate and might even out perform time-stepping methods offering
high-frequency numerical damping. Note that the addition of damping
does not impair the efficiency of the novel algorithm since a system
of equations has to be solved in each time step anyway due to the
implicit nature of the method. \\

In future research activities, the authors plan to extend the proposed
scheme to include controllable numerical damping where the spectral
radius of the amplification matrix is prescribed \emph{a priori}.
Thus, we will be able to construct a high-order family of time-stepping
methods that offer\emph{ A-} and \emph{L-}stable algorithms suitable
for different areas of application. Additionally, to allow for large-scale
transient analyses, the direct solver needs to be exchanged by a suitable
iterative solver including tailored pre-conditioners. This also opens
the pathway to a targeted parallelization of the algorithm and its
use on high-performance applications. 

\section*{Dedication}

This paper is dedicated to the late Dr.~John~P.~Wolf, who was a
mentor and close friend to the first author of this paper, Chongmin~Song.
The presented research is highly influenced by the scientific contributions
of Dr.~John~P.~Wolf and Dr.~Antonio~Paranesso~\citep{Wolf1991,Wolf1992}
regarding the development of lumped-parameter models for dynamic soil-structure
interaction analysis by means of Padé and partial fraction approximations.
Therefore, the novel time integration scheme proposed in this article
clearly benefits from the direct exposure to the guidance of Dr.~John~P.~Wolf
for which Chongmin~Song is very grateful. 

\section*{Acknowledgment}

This research was partially supported by the Australian Research Council
(ARC) under its Discovery Project scheme through Grant Numbers DP180101538
and DP200103577. Additionally, the authors are grateful for the help
provided by Dr Junqi Zhang (UNSW Sydney) in setting up the numerical
model of the mountainous region analyzed in Sect.~\ref{sec:Mountain}.

\appendix

\section{Matrix exponential and Taylor series}

\label{sec:Appendix:-Matrix-exponential} For the sake of completeness
and in the purpose of ensuring that the article is largely self-contained,
we will repeat the necessary theory that is directly related to this
paper. For more comprehensive information and a complete theory on
matrix functions the interested reader is referred to the pertinent
literature \citep{Gantmacher1977,Golub1996}.

The matrix exponential function of an arbitrary square matrix $\mathbf{A}$
can be defined using a Taylor series expansion 
\begin{equation}
e^{\mathbf{A}}=\mathbf{I}+\sum\limits _{k=1}^{\infty}\cfrac{1}{k!}\mathbf{A}^{k}\,,\label{eq:MPFunc-1}
\end{equation}
where $\mathbf{I}$ is an identity matrix of the same dimensions as
$\mathbf{A}$.

One important application of the matrix exponential function is to
solve generic IVPs resulting in a system of linear ODEs written as
\begin{equation}
\mathbf{y}'(x)=\mathbf{A}\mathbf{y}(x)\label{eq:MPFunc-7}
\end{equation}
with the IC 
\begin{equation}
\mathbf{y}(0)=\mathbf{y}_{0}\,.\label{eq:MPFunc-8}
\end{equation}
The prime $\square'$ denotes a derivative with respect to the argument
of the function $x$. The general solution to this ODE is given by
\begin{equation}
\mathbf{y}(x)=e^{\mathbf{A}x}\mathbf{y}_{0}\,.\label{eq:MPFunc-9}
\end{equation}
The matrix exponential $e^{\mathbf{A}x}$ is called the fundamental
matrix of the system of ODEs expressed in Eq.~\eqref{eq:MPFunc-7}.
An important mathematical property of the function $e^{\mathbf{A}}$
is that it is commutative with respect to the matrix $\mathbf{A}$
and its inverse $\mathbf{A}^{-1}$, i.e., 
\begin{equation}
\mathbf{A}e^{\mathbf{A}}=e^{\mathbf{A}}\mathbf{A}\qquad\text{and}\qquad\mathbf{A}^{-1}e^{\mathbf{A}}=e^{\mathbf{A}}\mathbf{A}^{-1}\label{eq:MPFunc-21}
\end{equation}
hold, which can be shown using the definition of the matrix exponential
function given in Eq.~\eqref{eq:MPFunc-1} by means of a Taylor series.

As a specific example, let us consider the case of free vibration
of an undamped SDOF system. Equation~\eqref{eq:ODE_SDOF} is written
as 
\begin{equation}
\ddot{u}+\omega_{n}^{2}u=0\label{eq:MPFunc-ex-eq}
\end{equation}
with the initial conditions (see Eq.~\eqref{eq:InitialCondition})
\begin{equation}
u(t=0)=u_{0}\quad\mathrm{and}\quad\dot{u}(t=0)=\dot{u}_{0}\label{eq:InitialCondition-1}
\end{equation}
The matrix $\mathbf{A}$ in Eq.~\eqref{eq:matrixA} is defined as
\begin{equation}
\mathbf{A}=\begin{bmatrix}0 & -\omega_{n}^{2}\\
+1 & 0
\end{bmatrix}\,.\label{eq:MPFunc-ex2-1}
\end{equation}
Here, $\mathbf{A}$ has two imaginary eigenvalues $\pm\mathrm{i}$.
The matrix exponential function $e^{\mathbf{A}t}$ is expressed according
to Eq.~\eqref{eq:MPFunc-1} as a infinite sum of power functions
\begin{equation}
e^{\mathbf{A}t}=\begin{bmatrix}1 & 0\\
0 & 1
\end{bmatrix}+t\begin{bmatrix}0 & -\omega_{n}^{2}\\
+1 & 0
\end{bmatrix}+\dfrac{t^{2}}{2!}\begin{bmatrix}-\omega_{n}^{2} & 0\\
0 & -\omega_{n}^{2}
\end{bmatrix}+\dfrac{t^{3}}{3!}\ \begin{bmatrix}0 & \omega_{n}^{4}\\
-\omega_{n}^{2} & 0
\end{bmatrix}+\dfrac{t^{4}}{4!}\begin{bmatrix}\omega_{n}^{4} & 0\\
0 & \omega_{n}^{4}
\end{bmatrix}+\cdots\label{eq:MPFunc-ex-EAt}
\end{equation}
This expression can be further simplified using the Taylor series
expansions of the sine- and cosine-functions 
\begin{subequations}
\begin{alignat}{2}
\cos(x) & =\sum\limits _{k=0}^{\infty}\dfrac{(-1)^{k}}{(2k)!}x^{2k} &  & =1-\dfrac{x^{2}}{2!}+\dfrac{x^{4}}{4!}-\ldots\\
\sin(x) & =\sum\limits _{k=0}^{\infty}\dfrac{(-1)^{k}}{(2k+1)!}x^{2k+1} &  & =x-\dfrac{x^{3}}{3!}+\dfrac{x^{5}}{5!}-\ldots
\end{alignat}
\end{subequations}
 Hence, the matrix exponential function in Eq.~\eqref{eq:MPFunc-ex-EAt}
can be given as 
\begin{equation}
e^{\mathbf{A}t}=\begin{bmatrix}\cos(\omega_{n}t) & -\omega_{n}\sin(\omega_{n}t)\\
\dfrac{1}{\omega_{n}}\sin(\omega_{n}t) & \cos(\omega_{n}t)
\end{bmatrix}\,.\label{eq:MPFunc-ex-3}
\end{equation}
Following Eq.~\eqref{eq:MPFunc-9}, the solution of Eq.~\eqref{eq:MPFunc-ex-eq}
with Eq.~\eqref{eq:InitialCondition-1} is obtained as
\begin{equation}
\left\{ \begin{array}{c}
\dot{u}\\
u
\end{array}\right\} =\begin{bmatrix}\cos(\omega_{n}t) & -\omega_{n}\sin(\omega_{n}t)\\
\dfrac{1}{\omega_{n}}\sin(\omega_{n}t) & \cos(\omega_{n}t)
\end{bmatrix}\left\{ \begin{array}{c}
\dot{u_{0}}\\
u_{0}
\end{array}\right\} \label{eq:MPFunc-ex-sln}
\end{equation}
which is consistent with the solution in Eq.~\eqref{eq:FrorcedVibrationsDamped}.

\section{Padé expansion of the exponential function}

\label{sec:Pade-approximation} In a Padé expansion, a function is
approximated as a ratio of two polynomials, i.e., a rational function.
The Padé approximation of a function $y(x)$ of order $\mathcal{O}(L,M)$
is written as 
\begin{equation}
y_{L/M}(x)=\cfrac{P_{L}(x)}{Q_{M}(x)}=\cfrac{p_{0}+p_{1}x+\ldots+p_{L}x^{L}}{q_{0}+q_{1}x+\ldots+q_{M}x^{M}}\,,\label{eq:Pade}
\end{equation}
where $L$ and $M$ are the degrees of polynomials in the numerator
and denominator, respectively. The constants $p_{k}$ ($k=0,1,\ldots,L)$
and $q_{k}$ ($k=0,1,\ldots,M)$ can be determined by equating the
Padé series with a Taylor expansion. In the present study, only the
\emph{diagonal} Padé approximations of the exponential function $e^{x}$
is of interest. Hence, the order of the approximation is $\mathcal{O}(M,M)$
or in short $M$, which is expressed as 
\begin{equation}
e_{M}^{x}=\cfrac{P_{M}(x)}{Q_{M}(x)}\,,\label{eq:ExpPade}
\end{equation}
where the polynomials in the numerator and denominator are given by
\begin{subequations}
\label{eq:Exp_M} 
\begin{align}
P_{M}(x) & =\sum\limits _{m=0}^{M}\cfrac{\left(2M-m\right)!}{m!\left(M-m\right)!}\;x^{m}\,,\label{eq:Exp_PM}\\
Q_{M}(x) & =P_{M}(-x)=\sum\limits _{m=0}^{M}\cfrac{\left(2M-m\right)!}{m!\left(M-m\right)!}\;(-x)^{m}\,.\label{eq:Exp_QM}
\end{align}
\end{subequations}
 From Eqs.~\eqref{eq:Exp_M}, it is obvious that $Q_{M}(x)=P_{M}(-x)$
holds. The accuracy of the approximation is of the order of $2M$.
The denominator polynomial $Q_{M}(x)$ of degree $M$ can be factorized
as 
\begin{equation}
Q_{M}(x)=\prod\limits _{i=1}^{M}\left(r_{i}-x\right)=\left(r_{1}-x\right)\left(r_{2}-x\right)\ldots\left(r_{M}-x\right)\,,\label{eq:Exp_QM_factor}
\end{equation}
where $r_{i}$ denote the roots of the polynomial, which are either
real or pairs of complex conjugates numbers. Note that a fraction
can be decomposed into partial fractions if the denominator is a factorized
polynomial, e.g., 
\begin{equation}
\cfrac{1}{\left(r_{1}-x\right)\left(r_{2}-x\right)}=\cfrac{1}{r_{2}-r_{1}}\left(\cfrac{1}{r_{1}-x}-\cfrac{1}{r_{2}-x}\right)\label{eq:PartialFraction}
\end{equation}
with the condition that $r_{1}\ne r_{2}$ holds. 

\section{Lagrange interpolation}

\label{sec:LagrangeInterpolation} A simple and efficient way of interpolating
functions by means of polynomials of any order $p$ is provided by
Lagrange interpolation, where the explicit calculation of the polynomial
coefficients is avoided. This polynomial interpolation scheme is implemented
by introducing a set of Lagrangian interpolation polynomials defined
using a particular set of points (often also called nodes) with the
coordinates $\xi_{k}\in\{1,2,\ldots,p+1\}$. It is easy to show that
a Lagrangian interpolation polynomial corresponding to the node $l$
is given as 
\begin{equation}
L_{l}(\xi)=\prod\limits _{k=1,\;k\ne l}^{p+1}\cfrac{\xi-\xi_{k}}{\xi_{l}-\xi_{k}}\,.\label{eq:LagrangeInterp}
\end{equation}
By construction, a Lagrangian interpolation polynomial holds the Kronecker-delta
property, i.e., it is equal to zero at all points, except the $l$th
point where it is equal to unity 
\begin{equation}
L_{l}(\xi_{k})=\delta_{lk}\,.\label{eq:Kronecker}
\end{equation}
In principle, the nodal positions can be chosen arbitrarily but this
does not guarantee stable and converging results. Especially, in the
context of the spectral element method (SEM) it has been shown that
GLL-points hold favorable approximation properties \citep{BookKarniadakis2005,BookPozrikidis2014}.
Here, two points are placed at the beginning and end of the interpolation
interval, while the other points are distributed in a non-equidistant
but symmetric fashion in the interior. Often the interpolation points
are given in the interval $[+1,-1]$ such that 
\begin{subequations}
\begin{align}
\xi_{1} & =-1\,,\\
\xi_{p+1} & =+1\,.
\end{align}
\end{subequations}
 The inner points are located at the GLL-points, which satisfy the
following condition 
\begin{equation}
\cfrac{\mathrm{d}}{\mathrm{d}\xi}\;\mathfrak{L}_{p}(\xi_{l})=0\quad\forall\;l=2,3,\ldots,p\,.
\end{equation}
Here, $\mathfrak{L}_{p}$ denotes the Legendre polynomial of order
$p$ and its derivative is generally referred to as Lobatto polynomial.
Using GLL-points the maximum of each Lagrangian interpolation polynomial
is located at its corresponding node and the value is unity. An arbitrary
function $F(\xi)$ is now approximated as 
\begin{equation}
F(\xi)\approx\hat{F}(\xi)=\sum\limits _{k=1}^{p+1}a_{l}L_{l}(\xi)\,,
\end{equation}
where $a_{l}$ are the values of the original function at the GLL-points
\begin{equation}
a_{l}=F(\xi_{l})\,.
\end{equation}
For the purpose of approximating the forcing term in the equations
of motion in terms of the dimensionless coordinate $s$ defined in
the interval $[0,1]$ a linear mapping needs to be introduced 
\begin{equation}
\xi=2s-1\,.
\end{equation}
Using the above expression, a function of $\xi$ can be rewritten
in terms of $s$. The interpolation polynomials can be naturally also
directly defined in terms of $s$ by mapping the GLL-points. 
\addcontentsline{toc}{section}{References} 
 \bibliographystyle{elsarticle-num}
\bibliography{bibliography/Literature}

\end{document}